\documentclass{amsart}

\usepackage[german, english]{babel}
\usepackage{amssymb,amsmath,amsthm,amscd}

\input{diagrams}

\advance\oddsidemargin by -0.1cm
\advance\evensidemargin by -0.7cm
\newtheorem{theorem}{Theorem}
\newtheorem{proposition}{Proposition}
\newtheorem{corollary}{Corollary}
\newtheorem{lemma}{Lemma}
\newtheorem{definition}{Definition}

\newcommand{\modulus}[1]{\left\lvert #1 \right\rvert}
\newcommand{\norm}[1]{\left\| #1 \right\|}
\newcommand{\mklm}[1]{\left\{ #1 \right\}}
\newcommand{\eklm}[1]{\left\langle #1 \right\rangle}
\newcommand{\mnd}{\!\!\!}

\newcommand{\Z}{{\mathbb Z}}
\newcommand{\C}{{\mathbb C}}
\newcommand{\R}{{\mathbb R}}

\newcommand{\D}{{\mathcal D}}
\newcommand{\E}{{\mathcal E}}
\newcommand{\F}{{\mathcal F}}
\renewcommand{\H}{{\mathcal H}}

\newcommand{\M}{{\mathcal M}}
\newcommand{\T}{{\rm T}}
\renewcommand{\P}{{\mathbb P}}
\newcommand{\Q}{{\mathbb Q}}
\renewcommand{\O}{{\mathcal O}}
\newcommand{\X}{{\mathcal X}}

\newcommand{\1}{{\bf 1}}
\renewcommand{\epsilon}{\varepsilon}
\renewcommand{\phi}{\varphi}
\renewcommand{\rho}{\varrho}
\newcommand{\HM}{{\mathfrak H}}

\newcommand{\Cinft}{{\rm C^{\infty}}}
\newcommand{\Ctest}{{\rm C^{\infty}_0}}
\newcommand{\Sob}{{\rm H}}
\renewcommand{\L}{{\rm L}}

\renewcommand{\O}{{\mathrm O}}
\newcommand{\SU}{\mathrm{SU}}
\newcommand{\U}{{\mathrm U}}

\newcommand{\Spin}{\mathrm{Spin}}

\renewcommand{\div}{\text{div}\,}
\renewcommand{\det}{\text{det}\,}
\renewcommand{\Im}{\mathrm{Im}\,}
\renewcommand{\Re}{\mathrm{Re}\,}
\newcommand{\vol}{\text{vol}\,}
\newcommand{\dist}{\text{dist}\,}

\DeclareMathOperator{\Ker}{Ker}
\DeclareMathOperator{\supp}{supp}

\DeclareMathOperator{\tr}{tr}
\DeclareMathOperator{\gd}{\partial}
\DeclareMathOperator{\gdq}{\overline{\partial}}
\DeclareMathOperator{\Ric}{Ric}
\DeclareMathOperator{\grad}{grad}
\DeclareMathOperator{\II}{II}
\DeclareMathOperator{\Dirac}{D}
\DeclareMathOperator{\Tr}{Tr}
\DeclareMathOperator{\arsinh}{arsinh}

\begin{document}

\author{Pablo Ramacher}
\title[Hypersurfaces in Eguchi--Hanson space]{Geometric and analytic properties of families of hypersurfaces in Eguchi--Hanson space}
\address{Pablo Ramacher, Humboldt--Universit\"at zu Berlin, Institut f\"ur Reine Mathematik, Sitz: Rudower Chaussee 25, D--10099 Berlin, Germany}
\subjclass{53C27,58J50,83C20}
\keywords{Eguchi--Hanson metric, Hopf tori, families of hypersurfaces, Dirac operator and Laplace operator, spectrum, Einstein--Dirac system}
\email{ramacher@mathematik.hu-berlin.de}
\thanks{Supported by the SFB 288 of the DFG}

\begin{abstract}
We study the geometry of families of hypersurfaces in Eguchi--Hanson space that arise as complex line bundles over curves in $S^2$ and are three--dimensional, non--compact Riemannian manifolds, which  are foliated in Hopf tori  for closed cur\-ves. They are negatively curved, asymptotically flat spaces, and we compute the complete three--dimensional curvature tensor as well as the second fundamental form, giving  also some results concerning their  geodesic flow. We show  the non--existence of  $\L^p$--harmonic functions on these hypersurfaces for every $p \geq 1$ and arbitrary curves, and determine the  infima of the spectra of the Laplace and of the square of the Dirac operator in the case of closed curves. We also show that, in this case, zero lies in the spectrum of the Dirac operator. For circles we compute the $\L^2$--kernel of the Dirac operator in the sense of spectral theory and show that it is trivial. We consider further the Einstein--Dirac system on these spaces and construct explicit examples of $T$--Killing spinors on them.
\end{abstract}

\maketitle

\tableofcontents

\section{Introduction}

In this paper we shall study certain families of hypersurfaces in Eguchi--Hanson space that arise as complex line bundles over curves on $S^2\simeq\C\cup\mklm{\infty}$. They are three-dimensional, open, asymptotically flat  Riemannian manifolds of non--positive scalar curvature which, in case of a closed curve, are foliated in  Hopf tori. We describe their geometry in detail, computing the complete three--dimensional curvature tensor as well as the second fundamental form, and give  also some results on the structure of the geodesic flow.
Since an explicit description of the geometric properties of these hypersurfaces is possible, we are able to make precise statements about the spectra of the scalar Laplacian and the Dirac operator and also about the existence of solutions of spinorial field equations. In particular, we show that there are no $\L^p$--harmonic functions for every $p \geq 1$ and arbitrary curves, and that for curves arising by  M\"{o}bius transforms from closed curves the  spectra of the scalar Laplacian and the square of the Dirac operator come arbitrarily close to zero, implying that zero lies in the spectrum of the considered operators. In the mentioned case, it also turns out that zero lies in the spectrum of the Dirac operator. In case that the considered curves are generalized circles in $\C$ that arise by M\"{o}bius transforms from circles in $\C$ with center at the origin the $\L^2$--kernel of the Dirac operator in the sense of spectral theory can be computed explicitly and we show that it is trivial. 
 As it turns out, these hypersurfaces  do not admit solutions to the Einstein--Dirac system; such solutions can only be obtained by deformation into a singular situation. Nevertheless, we can construct explicit examples of $T$--Killing spinors, which are solutions of  a generalized Killing equation for spinors.

Hopf tori have been extensively studied, see e.g. \cite{willmore}, and where first considered by Pinkall \cite{pinkall}. If $\pi: S^3 \rightarrow S^2$ denotes the Hopf fibration, the inverse image of any closed curve in $S^2$ will be an immersed torus in $S^3$, which is called a \emph{Hopf torus}. Using Hopf tori Pinkall showed that every compact Riemann surface of genus one can be conformally embedded as a flat torus into the unit sphere $S^3$. As a further application, and using elastic curves in $S^2$, he constructed new examples of compact embedded Willmore surfaces in $\R^3$, which are extremal surfaces for the Willmore functional $\int \HM \,  dA$, where $\HM$ denotes the mean curvature.

 The Eguchi--Hanson metric is a four--dimensional metric, which can be constructed in  the total space of the fibration $p:T^\ast\P^1(\C)\rightarrow \P^1(\C)\simeq S^2$, and since  its holonomy is contained in $\SU(2)$, it is Ricci flat and self--dual. Both the Hopf fibration $\pi$ and the projection $p$ are compatible with the action of $\U(2)$ in $\C\cup \mklm{\infty}$, and, like  the standard metric in $S^3$, the Eguchi--Hanson metric  is invariant under this action. Therefore, its restriction to  the three--dimensional projective space $\P^3(\R)$, which is  immersed  in $\T^\ast\P^1(\C)$ as the set of all cotangential vectors of unit length, corresponds exactly to the standard metric in $S^3$. For this reason the projection $p$ is a geometric extension of the Hopf fibration, and the preimage of any closed curve on $S^2$ under the projection $p$ gives rise to a three--dimensional non--compact Riemannian manifold foliated in Hopf tori. Its end is of topological type $T^2 \times (0,\infty)/ \mklm{\pm 1}$, where $T^2$ is the two--dimensional torus. 
Nevertheless, the corresponding Willmore functional turns out to be unbounded, so that the considered hypersurfaces are not accessible to integral geometry.
The interest in Eguchi--Hanson space itself originates from a result of Schoen and Yau \cite{schoen-yau}, who proved that a complete asymptotically Euclidean four--manifold whose Ricci tensor vanishes is necessarily flat. For Ricci flat asymptotically locally Euclidean  K\"{a}hler metrics  this turns out not to be true, the first example of such a metric being given by the Eguchi--Hanson metric \cite{eguchi-hanson}. 

\medskip

We give now a description of the main results of this  work. The Sections \ref{sec:1}, \ref{sec:2}, \ref{sec:3} and \ref{sec:8} are concerned with the geometry of the hypersurfaces studied, the Sections \ref{sec:4}, \ref{sec:6} and \ref{sec:7} with the spectra of the Dirac and the Laplace operator, while Section \ref{sec:5} is devoted to the study of spinorial field equations. The Eguchi--Hanson metric is described in Section \ref{sec:1}: it depends on a real parameter $t>0$, thus giving rise to a one--parameter family of Riemannian metrics $g_t$. These metrics become degenerate along the zero section in case that $t=0$. For any curve $\Gamma(s)=r(s) e^{i\phi_\Gamma(s)}$ in $\C\cup \mklm{\infty}\simeq\P^1(\C)$ we consider  its preimage $M^3_\Gamma :=p^{-1}(\Gamma)$
and obtain a family of hypersurfaces $(M^3_\Gamma,h_t)$, where we assume that $\Gamma(s)$ is parametrized by arc length and $h_t$ denotes the induced Riemannian metric. Each of these hypersurfaces is a complex line bundle over $\Gamma$, and introducing the polar coordinates $\rho$ and $\phi$ in each fiber, we obtain a parametrization of $M^3_\Gamma$ outside the zero section by the coordinates $s,\, \rho,\, \phi$, see Section \ref{sec:2}. Since the coefficients of $h_t$ do not depend on $\phi$, the corresponding $S^1$--symmetry is an isometry. We determine the inner geometry of the hypersurfaces and in Theorem \ref{thm:2.1}, page \pageref{thm:2.1}, the complete Ricci tensor is computed with respect to an orthonormal frame, 
one eigenvalue being positive, one negative and the third one becoming negative at infinity, yielding, for the scalar curvature, the expression
\begin{displaymath}
S=-\frac{\rho^4(r^2+1)^2}{(\rho^4(r^2+1)^2+t^4)^{3/2}}.
\end{displaymath}
 It is negative and tends  to zero for large $\rho$ and $r$ with the order $1/\rho^2(r^2+1)$. For $ t \not=0$, $S$ remains regular at $\rho=0$, i.e., the scalar curvature vanishes on the zero section. In Section \ref{sec:3} we turn to the study of the Levi--Civita connection of the Eguchi--Hanson space and determine the second fundamental form of the hypersurfaces $M^3_\Gamma$ with respect to the above orthonormal frame, thus obtaining
\begin{displaymath}
\II^\ast =\left ( \begin{array}{ccc} 0 & 0 & 0 \\ 0 & 0 & \frac{K}{4\rho}\sqrt{\frac K {r^2+1}} \\ 0 &\frac{K}{4\rho}\sqrt{\frac K {r^2+1}}  & \HM
\end{array} \right ),
\end{displaymath}
see Theorem \ref{thm:3.1} on  page \pageref{thm:3.1}, and Corollary \ref{cor:3.1} on  page \pageref{cor:3.1}, where $K$ is the function $K=2\rho^2(r^2+1)^2/\sqrt{\rho^4(r^2+1)^2+t^4}$ and  $\HM$ denotes the mean curvature.  It is given by the geodesic curvature $k_g$ of $\Gamma$ as a curve in $S^2$ according to the formula
\begin{align*}
\HM=\sqrt{\frac 2 {\sqrt{\rho^4(r^2+1)^2+t^4}}}\cdot k_g.
\end{align*}
This appears to be natural, since the geometry of the vector bundle $T^\ast \P^1(\C)$ is determined by the elliptic geometry of $\P^1(\C)\simeq S^2$. The above formula also implies that $M^3_\Gamma$ is a minimal surface if and only if $k_g=0$, i.e., if $\Gamma$ is a great circle in $S^2$. Further, since the function $\rho^2(r^2+1)$ corresponds to the distance in $\T^\ast\P^1(\C)$, both the scalar curvature $S$ and the mean curvature $\HM$, as well as the components of the Ricci tensor and of the second fundamental form, are manifestly invariant under the action of the isometry group $\U(2)$.
Section \ref{sec:8} contains some results concerning the geodesic flow of the hypersurfaces $M^3_\Gamma$. So, in case $\Gamma$ is a circle in $\C$ with center at the origin, we are able to compute the distance of a point in $M^3_\Gamma$ to the curve $\Gamma\subset M^3_\Gamma$, i.e., to the zero section, see Proposition \ref{prop:8.1} on page \pageref{prop:8.1}, and, in this way, to calculate the exponential growth of $M^3_\Gamma$ explicitly. 

In Section \ref{sec:4} the vanishing of the $\L^p$--kernel, $p\geq 1$, of the scalar Laplacian on the hypersurfaces  $(M^3_\Gamma,h_t)$ is proved for every $t\geq 0$ and every curve $\Gamma$ by showing the existence of a canonical exhaustion function on the considered hypersurfaces, see Proposition \ref{prop:4.1} and Corollary \ref{cor:4.1} on  page \pageref{prop:4.1}. The result then follows from the work of Greene and Wu \cite{greene-wu}, who studied integrals of certain generalized subharmonic functions on connected non--compact Riemannian manifolds admitting such a function, and showed that these integrals cannot be bounded. For the smallest spectral value of the scalar Laplacian we obtain, in Section \ref{sec:7}, the first estimate 
\begin{displaymath}
\mu_0(M^3_\Gamma) \leq t^{-2}, \qquad t>0,
\end{displaymath}
where $\Gamma$ is a closed curve, see Corollary \ref{cor:7.1} on page \pageref{cor:7.1}, since 
by general theory lower bounds for the Ricci tensor of open complete manifolds imply upper bounds for the smallest spectral value of the Laplace operator \cite{eichhorn}. By using the $\min$--$\max$ principle we are then able to determine the infimum of the spectrum of the closure of the Laplacian $\overline{\Delta}$ on $(M^3_\Gamma, h_t)$, obtaining
\begin{displaymath}
\inf \sigma (\overline{\Delta})< \delta
\end{displaymath}
for every $\delta>0$ and arbitrary $t>0$ and closed curves $\Gamma$. Since, by Corollary \ref{cor:4.1},  zero can be no $\L^2$--eigenvalue, we therefore get that zero lies in $\sigma_{\mathrm{ess}}(\overline \Delta)$, the essential spectrum of $\overline \Delta$. A result of Brooks \cite{brooks} then implies that in this case the hypersurfaces $M^3_\Gamma$ must be of subexponential growth, generalizing the previously obtained result.
Section \ref{sec:5} is devoted to the study of spinorial field equations. In \cite{friedrich-kim} Friedrich and Kim showed that in dimension 3 the existence of a solution to the Einstein Dirac system is equivalent  to the existence of a so--called Weak Killing or WK spinor. For the existence of such a spinor geometric integrability conditions that are independent of the considered spin structure are known, and we show that, for $t>0$, these conditions can never be fulfilled, implying that there cannot be any solutions to the Einstein Dirac system on the hypersurfaces $(M^3_\Gamma,h_t)$ for any $t>0$ and any curve $\Gamma$, see Proposition \ref{prop:5.2} on page \pageref{prop:5.2}. Nevertheless, such solutions can be constructed explicitly with respect to the trivial spin structure in case that $t=0$, the manifolds considered then being no longer complete. As remarked above, the Eguchi--Hanson metric is self--dual and, due to  this, there is a parallel spinor on Eguchi--Hanson space. By restricting this spinor to the hypersurfaces $M^3_\Gamma\subset T^\ast\P^1(\C)$  we show in Proposition \ref{prop:5.3} that there exists a $T$--Killing spinor on $M^3_\Gamma$ if and only if $M^3_\Gamma$ is a minimal surface. The spectrum of the Dirac operator $\Dirac$ is studied in Section \ref{sec:6}. There we show, by estimating the Rayleigh quotient from above and using again the $\min$--$\max$ principle, that the infimum of the spectrum of $\overline{\Dirac}^2$ on $(M^3_\Gamma,h_t)$ becomes arbitrarily small, 
\begin{displaymath}
\inf \sigma (\overline{\Dirac}^2)< \delta,
\end{displaymath}
where $\delta >0$ and $\Gamma$ is a closed curve, see Theorem \ref{thm:6.1} on page \pageref{thm:6.1}, $t>0$ being arbitrary; here the involved spin structure is again the trivial one. In this case it also follows that $0 \in \sigma(\overline D)$,  by explicit construction of an approximating sequence. 
In case that $\Gamma$ is a circle in $\C$ with center at the origin, an isometric $S^1 \times S^1$--action is given and the $\L^2$--kernel of the Dirac operator and of its closure decompose into the unitary representations of this action according to the spectral decomposition of the corresponding generators $i\gd_\phi$, $i\gd_s$, and with respect to the trivial spin structure one obtains 
\begin{displaymath}
\Ker_{\L^2}(\Dirac_{|(M^3_\Gamma\setminus \Gamma)})=\bigoplus\limits_{\beta=-1,-2,\dots} \H_0 \otimes H_\beta,
\end{displaymath}
while on $M^3_\Gamma$ the $\L^2$--kernels of the Dirac operator and its closure turn out to be trivial. Thus, in this case, $0 \in \sigma_{\mathrm{ess}}^{\L^2}(\overline \Dirac)$. Since $M^3_\Gamma$ and $M^3_{A\Gamma}$ are isometric for every $A \in \U(2)$, statements for a particular curve in $\C\cup \mklm{\infty}$ can be generalized to curves that arise from it by M\"{o}bius transforms.

\section{The Eguchi--Hanson space $(T^\ast\P^1(\C),g_t)$}
\label{sec:1}

Let $G$ be a finite nontrivial subgroup of $\U(m)$ that acts freely on $\C^m \setminus \mklm{0}$. Then $\C^m/G$ carries an isolated quotient singularity at zero and any resolution $(M,\pi)$ of $\C^m/G$ is a non--compact complex manifold. A K\"ahler metric $g$ on $M$ is said to be \emph{asymptotic} to the Euclidean metric $h$ on $\C^m/G$ if there is a smooth surjective map $f:M\rightarrow  \C^m/G$ such that $f^{-1}(0)$ is a connected, simply connected, finite union of compact submanifolds of $M$ and $f$ induces a diffeomorphism $M/f^{-1}(0)\simeq (\C^m/\mklm{0})/G$. Under this diffeomorphism $f_\ast(g)$ should satisfy
\begin{equation}
\label{eq:1.1}
f_\ast(g)=h+\O(r^{-4}), \qquad \nabla f_\ast(g)=\O(r^{-5}), \qquad \nabla^2 f_\ast(g)=\O(r^{-6})
\end{equation}
for large $r$, where $r$ is the distance from the origin  and $\nabla$ is the flat connection in $\C^m/G$. Such a metric is called an \emph{asymptotically locally Euclidean} or \emph{ALE metric}. Notice that the topological type of the end is given by a quotient of the Euclidean space. In the following we will mainly be concerned with the case of $m=2$.

In \cite{schoen-yau} Schoen and Yau proved that a complete asymptotically Euclidean four--manifold whose Ricci--tensor vanishes is necessarily flat. Nevertheless, a similar statement for  Ricci--flat ALE K\"ahler--metrics does not hold, since, as mentioned, the topology of the end differs from the topology of Euclidean space.  An important class of Ricci--flat K\"ahler metrics which give rise to ALE spaces is given by the so--called hyperk\"ahler structures. 
In the case of an oriented four--dimensional smooth manifold $X$ a hyperk\"ahler structure is a metric whose holonomy is contained in $\SU(2)$. A manifold with such a structure is Ricci--flat and self--dual, and its metric is K\"ahler with respect to each of the three anticommuting  complex structures. Alternatively, a hyperk\"ahler structure on $X$ may be defined to be a triple of smooth, closed $2$--forms $\sigma_1,\, \sigma_2,\, \sigma_3$ on $X$ that can be represented locally according to
\begin{equation}
\label{eq:1.2}
\sigma_1=l_1 \wedge l_4+l_2 \wedge l_3, \qquad \sigma_2=l_1 \wedge l_3 -l_2 \wedge l_4, \qquad \sigma_3=l_1 \wedge l_2+l_3 \wedge l_4,
\end{equation}
where $(l_1, \dots, l_4)$ is a local oriented frame of $1$--forms on $X$. 
The systematic construction of ALE metrics with holonomy $\SU(2)$ as hyperk\"ahler quotients was initiated by Hitchin \cite{hitchin} and carried over by Kronheimer \cite{kronheimer1, kronheimer2}, who studied the spaces $\C^2/G$ for general polyhedra groups $G\subset \SU(2)$ and showed the existence of hyperk\"ahler metrics on the resolution $M$ for the considered groups $G$, giving a complete classification. For cyclic groups these metrics are explicitly known. 

\medskip

The first example of a hyperk\"ahler ALE four--manifold was found by Eguchi and Hanson \cite{eguchi-hanson}. We  will now briefly proceed to describe its construction. Let $m=2$ and $G=\Z_2$ and consider the mapping
\begin{displaymath}
\Phi: \C^2 \longrightarrow \C^3, \qquad \Phi(z_1,z_2)=(z^2_1,z^2_2,z_1z_2).
\end{displaymath}
The image of $\C^2$ under $\Phi$ is 
\begin{displaymath}
X\,:=\, \Im \Phi\,=\, \mklm{(w_1,w_2,w_3)\in \C^3:w_1w_2=w_3^2},
\end{displaymath}
and $\Phi$ induces a bijection $\Phi: \C^2/\mklm{\pm 1}\rightarrow X$ so that $X$ becomes analytically equivalent to $\C^2/\mklm{\pm 1}$. The canonical bundle over $\P^1(\C)$,
\begin{gather*}
H\,:=\,\mklm{(l,v) \in \P^1(\C) \times \C^2: v \in l},
\end{gather*}
can be described explicitly as follows. If one introduces the homogeneous coordinates $[\alpha:\beta]$ in $P^1(\C)$, then the total space $H$ consists of all equivalence classes of triples $[\alpha,\beta,\gamma]$ with respect to the equivalence relation $(\alpha,\beta,\gamma) \sim (A \alpha,A \beta, {\gamma}/{A}) ,$ where $A \in \C^\ast$, i.\ e.
\begin{gather*}
H \, =\, \mklm{(\alpha, \beta, \gamma) \in (\C^2\setminus \mklm{0})\times \C}\big / \sim.
\end{gather*}
The one--dimensional complex tangential bundle  $\T\P^1(\C)$ is biholomorphic to the square of the dual of the canonical bundle \cite{milnor}
\begin{displaymath}
\T\P^1(\C) \, =\, H^\ast \otimes H^\ast \, ,
\end{displaymath} 
from which one obtains, for the cotangential bundle $T^\ast\P^1(\C)$, the description
\begin{displaymath}
\T^\ast\P^1(\C)=H^2=\mklm{(\alpha,\beta,\gamma) \in (\C^2 \setminus \mklm{0}) \times \C} \big / \sim_1,
\end{displaymath}
with the equivalence relation $(\alpha,\beta,\gamma) \sim_1 (A \alpha, A \beta,{\gamma}/{A^2})$.
Notice that  $H^2$ is simply-connected.
We define now the mapping
\begin{displaymath}
\pi: H^2 \longrightarrow X, \quad \pi([\alpha,\beta,\gamma])=(\alpha^2 \gamma, \beta^2 \gamma, \alpha \beta \gamma).
\end{displaymath}
The preimage of the point $(0,0,0)$ under $\pi$ is the zero section of the bundle $H^2$. Away from this set $\pi:H^2\setminus \P^1(\C) \longrightarrow X \setminus \mklm{(0,0,0)}$ is bijective,  and hence  $(H^2,\pi)$ represents a resolution of the singularity of $\C^2/\mklm{\pm 1}$ at zero. Summing up one obtains the diagram
\begin{diagram}
                         &                    & H^2  &                          &              \\
                         & \ldTo^{\pi_1}      &      & \rdTo(2,1)^\pi           & X            \\
\C^2/\mklm{\pm 1} \qquad &                    &      & \ruTo(4,1)_{\simeq}^\Phi &              \\
\end{diagram}
where the mapping $\pi_1$ is given by the formula $\pi_1([\alpha,\beta,\gamma])=[\alpha\sqrt{\gamma},\beta\sqrt{\gamma}]$.
The closed holomorphic 2--form $dz_1\wedge dz_2$ and the function $u_1:=\modulus{z_1}^2+\modulus{z_2}^2$ on $\C^2$ are invariant under reflections at the origin, descend to $\C^2/\mklm{\pm1}$ and, thus, lift to forms on $H^2$, which we will denote by $dz_1\wedge dz_2$ and $u_1$ as well. 

We come now to the description of the Eguchi--Hanson metric.
Following \cite{lebrun} we consider, on the complex manifold $H^2$, the family of real--valued functions $f_t \in \E^{0,0}(H^2)$ depending on the parameter $t$, 
\begin{gather*}
f_t:=\sqrt{u_1^2+t^4}+t^2 \log \frac{u_1}{\sqrt{u_1^2+t^4}+t^2}, \quad t>0. 
\end{gather*}
Here the function $u_1:H^2 \rightarrow \R$ is explicitly given by $u_1([\alpha,\beta,\gamma])=(\modulus{\alpha}^2+\modulus{\beta}^2)|\gamma|$, from which it follows that, away from the exceptional curve, i.e., the zero section, $u_1$ is a smooth function, and the same holds for $f_t$. 
For $t>0$ the associated form 
\begin{displaymath}
\omega_t:=i\gd \gdq f_t \in \E^{1,1}(H^2)
\end{displaymath}
is regular even in the exceptional curve and thus defines a K\"ahler form on $H^2$.
For using homogeneous coordinates we can define  a complex analytic structure on $H^2$ as follows. Let $U_\alpha=\mklm{[\alpha,\beta,\gamma]: \alpha \not=0},\,U_\beta=\mklm{[\alpha,\beta,\gamma]: \beta \not=0}$ be open subsets in $H^2$ and define the homeomorphisms
\begin{align*}
h_\alpha: U_\alpha \rightarrow \C^2, \quad [\alpha,\beta,\gamma]=[1,{\beta}/{\alpha},\gamma \alpha^2]&\longmapsto ({\beta}/{\alpha},\gamma \alpha^2),\\
h_\beta: U_\alpha \rightarrow \C^2, \quad [\alpha,\beta,\gamma]=[{\alpha}/{\beta},1,\gamma \beta^2]&\longmapsto ({\alpha}/{\beta},\gamma \beta^2).
\end{align*}
Since $h_\alpha \circ h_\beta^{-1}:h_\beta(U_\alpha \cap U_\beta) \rightarrow h_\alpha (U_\alpha \cap U_\beta)$ is a biholomorphic mapping, this gives a complex analytic structure on $H^2$. We can therefore choose the functions $\beta$ and $\gamma$ as local coordinates  in $U_\alpha$ by setting $\alpha$ equal to $1$,  so that $\gd:\E^{(p,q)}(H^2)\rightarrow \E^{(p+1,q)}(H^2)$ and $\gdq:\E^{(p,q)}(H^2)\rightarrow \E^{(p,q+1)}(H^2)$ are given by 
\begin{displaymath}
\gd=\frac{\gd}{\gd \beta} d\beta+\frac{\gd}{\gd \gamma} d\gamma, \qquad \gdq=\frac{\gd}{\gd \overline \beta} d\overline \beta+\frac{\gd}{\gd \overline \gamma} d\overline\gamma
\end{displaymath}
on $U_\alpha$. The regularity of $\omega_t$ for $t>0$ then follows by noting that the derivatives of $f_t$ with respect to $\gamma, \beta, \overline \gamma$ and $\overline \beta$ become regular (note that $u_1=(1+\beta\overline \beta)\sqrt{\gamma\overline\gamma}$ ). For example, in $\C$, one has 
\begin{displaymath}
\gd \gdq \left (\sqrt{|z|^2+t^4}+t^2 \log \frac{|z|}{\sqrt{|z|^2+t^4}+t^2}\right)=\frac{2t^6+2t^2 \bar z z +(2t^4+\bar z z) \sqrt{t^4+\bar zz}}{4(t^4+\bar zz)(t^2+\sqrt{t^4+\bar zz})^2} dz \wedge d\bar z.
\end{displaymath}
In case that $t=0$ on has $f_0=u_1$, and $\omega_0$ becomes degenerate along the zero section.
On $H^2$ the K\"ahler form $\omega_t$ induces a Riemannian metric through the formula 
\begin{displaymath}
g_t(X,Y):=\omega_t(X,JY), \qquad X,Y \in \X(H^2),
\end{displaymath}
where $J$ denotes the complex structure of $H^2$. For  $t\not=0$ $(H^2,g_t)$ becomes a complete Riemannian manifold.

 The complex manifold $M^4:=H^2 \setminus \P^1(\C)$ is an open dense subset of $H^2$, so it suffices for the study of  the geometric properties of $H^2$ to consider $\omega_t$ as well as the other relevant geometric objects just on $M^4$. Further, since  $\pi_1$ maps $M^4$ bijectively  onto $\C^2\setminus \mklm{0}/\mklm{\pm1}$,  $\omega_t$ can be explicitly computed on $M^4$ with respect to the coordinates $z_1,z_2$. For $u_1=z_1\bar z_1+z_2\bar z_2 \in \E^{0,0}(\R^4)$ as a function on $\C^2$ one has therefore
\begin{gather*}
du_1=\sum\limits_{i=1,2} \bigg (\frac{\gd u_1}{\gd z_i} dz_i+ \frac{\gd u_1}{\gd \bar z_i} d\bar z_i\bigg ) \in \E^{1,0}(\R^4) \oplus \E^{0,1}(\R^4), \\
\gdq u_1=z_1 d\bar z_1+z_2 d\bar z_2, \qquad \gd u_1=\bar z_1 dz_1+\bar z_2 dz_2,
\end{gather*}
see e.\ g.\ \cite{wells}, which yields $\gdq f_t=(\sqrt{u_1^2+t^4}/u_1) \gdq u_1$ and thus
\begin{equation*}
\begin{split}
\omega_t&=-i\frac{t^4}{u_1^2\sqrt{u_1^2+t^4}} \Big\{\modulus{z_1}^2 dz_1 \wedge d\bar z_1+\modulus{z_2}^2 dz_2 \wedge d\bar z_2+z_1\bar z_2 \,dz_2\wedge d\bar z_1  \\ & +\bar z_1 z_2 \,dz_1\wedge d\bar z_2\Big\} +i \frac{\sqrt{u_1^2+t^4}}{u_1}\mklm{dz_1\wedge d\bar z_1+dz_2\wedge d\bar z_2}.
\end{split}
\end{equation*}
The form $dz_i\wedge d\bar z_j$ is expressed with respect to the coordinates $z_1=x_1+iy_1, \,z_2=x_2+iy_2$ by
\begin{displaymath}
dz_i\wedge d\bar z_j=dx_i\wedge dx_j+dy_i \wedge d y_j-i(dx_i \wedge d y_j+d x_j \wedge d y_i),
\end{displaymath}
and the action of $J$ is given by $J(\gd_{x_i})=\gd_{y_i},\, J(\gd_{y_i})=-\gd_{x_i}$. Computation of $g_t$ restricted to $M^4$ then gives 
\begin{displaymath}
g_t=\left(\begin{array}{cccc} G_1&0&-G_4&-G_3 \\ 0&G_1&G_3&-G_4 \\ -G_4&G_3&G_2&0 \\  -G_3&-G_4&0&G_2 \end{array} \right),
\end{displaymath}
where 
\begin{align*}
G_1&=G-H(x_1^2+y_1^2), & G_2&=G-H(x_2^2+y^2_2), \\
G_3&=H(x_1y_2-y_1x_2), & G_4&=H(x_1x_2+y_1y_2), \\
\end{align*}
and $G,\,H:M^4\rightarrow \R$ are the smooth functions
\begin{displaymath}
G:=\frac{2\sqrt{u_1^2+t^4}}{u_1}, \qquad H:=\frac{2t^4}{u_1^2\sqrt{u_1^2+t^4}}. 
\end{displaymath}
For later use we define the smooth function
\begin{displaymath}
 K:M^4 \rightarrow \R,\qquad K:=G-Hu_1=4G^{-1}.
\end{displaymath}
From this it becomes evident that $g_t$ satisfies condition \eqref{eq:1.1}. 
One further computes the volume form $dM^4=\omega_t \wedge \omega_t$ to be $V(dz_1\wedge d\bar z_1\wedge dz_2 \wedge \bar dz _2$), where $V\equiv 2$. By the general theory of K\"ahler manifolds \cite{kobayashi-nomizu} the Ricci--form $\Ric=i\gd \gdq \log V$ then vanishes and it follows immediately that the Riemannian curvature tensor with respect to the decomposition $\bigwedge^2(M^4)=\bigwedge^2_+(M^4)\oplus \bigwedge^2_-(M^4)$ is given by
\begin{displaymath}
R=\left ( \begin{array}{cc} W_+ & 0 \\ 0 & W_- \end{array} \right ) +\left ( \begin{array}{cc} 0 & B^\ast \\ B & 0 \end{array} \right ) - \frac{S}{12}=\left ( \begin{array}{cc} W_+ & 0 \\ 0 & 0 \end{array} \right ),
\end{displaymath}
where $W_-$ and $W_+$ are the negative and positive part of the Weyl tensor respectively, $B$ is the trace--free part of the Ricci tensor and $S$ denotes the scalar curvature. The condition $B=0$ implies that $(H^2,g_t)$ is an Einstein space and the vanishing of $W_-$ means that $(H^2,g_t)$ is self--dual; the latter is equivalent to the statement that the bundle $\bigwedge^2_-(M^4)$ is flat which in turn implies that there exist three parallel forms on $\bigwedge^2_-(M^4)$. These forms can be chosen as $\omega_t$ and the two closed $2$--forms $\sigma_2,\,\sigma_3$ defined by $\sigma_2+i\sigma_3=dz_1\wedge dz_2$. One can show that the triple $(\omega_1,\sigma_2,\sigma_3)$ may locally be written in the form \eqref{eq:1.1} and thus forms a hyperk\"ahler structure on $M^4$ and hence on $H^2$. 

We consider now the projection 
\begin{displaymath}
p: H^2=T^\ast\P^1(\C) \longrightarrow \P^1(\C)\cong \C \cup \mklm{\infty} \cong S^2,
\end{displaymath}
 which is explicitly given by $[\alpha, \beta,\gamma] \mapsto [\alpha:\beta]\mapsto {\alpha}/{\beta}$. The function $u_1$ is invariant under the standard action of $\U(2)$  on $\C^2$ resp.\ $\C^2\setminus\mklm{0}/\mklm{\pm1}\simeq M^4$, which is given by its matrix representation.
On the other hand, $\U(2)$ acts as a group of holomorphic transformations on $\C\cup \mklm{\infty}$ by the so--called \emph{M\"{o}bius transform}
\begin{displaymath}
\left ( \begin{array}{cc} a & b \\ c & d \end{array} \right ) z=\frac{az+b}{cz+d},
\end{displaymath}
resp.\ on $\P^1(\C)$ by 
\begin{displaymath}
\left ( \begin{array}{cc} a & b \\ c & d \end{array} \right ) [\alpha:\beta]=[a\alpha+b\beta:c\alpha+d\beta].
\end{displaymath}
Taking the mapping $\tilde p:=p \circ \pi^{-1}_1: \C^2 \setminus\mklm{0}/\mklm{\pm1} \rightarrow \P^1(\C),\,$ which is given by $\tilde p[z_1,z_2] = [z_1:z_2]$, one therefore sees that 
\begin{displaymath}
\tilde p \left (\left ( \begin{array}{cc} a & b \\ c & d \end{array} \right )[z_1,z_2] \right)=\left ( \begin{array}{cc} a & b \\ c & d \end{array} \right ) \tilde p[z_1,z_2],
\end{displaymath}
which means that the diagram
\begin{diagram}
   H^2\setminus \P^1(\C) & \rTo ^{\quad \pi_1\quad }_{\quad\simeq\quad}  & \C^2\setminus \mklm{0} /\mklm{\pm1}\\
      \dTo  \,p          & \ldTo_{\tilde p}        &                                    \\
    \P^1(\C)             &                         &                                    \\
\end{diagram}
is compatible with the group action of $\U(2)$. Since the set $\C^2\setminus\mklm{0}/\mklm{\pm 1}$ is mapped by $\U(2)$ onto itself, it follows that, by extending the action of $\U(2)$ to $H^2$, the exceptional curve in $H^2$ must be mapped onto itself, too. The projection $p:H^2\rightarrow \P^1(\C)$ is  therefore also  compatible with the $\U(2)$--action.  Since $u_1$ vanishes on the zero section, it  becomes $\U(2)$--invariant on $H^2$. The K\"{a}hler metric $\omega_t$ and the Riemannian metric $g_t$, which are defined by means of the function $u_1$, are thus also invariant under $\U(2)$ in $H^2$.

\section{Hypersurfaces in ($T^\ast\P^1(\C),g_t)$ and their inner geometry}
\label{sec:2}

We now introduce certain hypersurfaces in $T^\ast\P^1(\C)$ and to this end consider  for any curve $\Gamma(s)=u(s)+iv(s)=r(s) e^{i\phi_\Gamma(s)}$ in $\C \cup \mklm{\infty}$ its preimage 
$$M^3_\Gamma :=p^{-1}(\Gamma)=\mklm{[\Gamma(s),1,\gamma]\in H^2:\gamma \in \C},$$ 
obtaining  a real three--dimensional hypersurface  in $H^2$. Let $h_t$ be the Riemannian metric on $M^3_\Gamma$ induced by $g_t$. The three--manifold $M^3_\Gamma$ is open and in case of a closed curve its end is  of topological type $T^2\times (0,\infty)/\mklm{\pm 1}$, where $T^2=S^1\times S^1$ is the two--dimensional torus. The hypersurfaces $M^3_\Gamma$ are asymptotically flat, but no ALE spaces, since their end is not modeled  on the end of $\R^3/G$
.  Note that $M^3_\Gamma$ is a one--dimensional complex vector bundle over $\Gamma$.

Since $p:H^2\rightarrow \P^1(\C)$ is compatible with the action of $\U(2)$, and since  $g_t$ and hence  $h_t$ are invariant under this action, 
$M^3_\Gamma$  is mapped  isometrically onto $M^3_{A\, \Gamma}$, where $A \in \U(2)$. Remember  that under  M\"{o}bius transforms generalized circles in $\C$ are mapped again into generalized circles.

We will now compute the inner geometry of the hypersurfaces $M^3_\Gamma$ and
assume from now on that $\Gamma(s)$ is parametrized by arc length. Using the projection $p$ one obtains a  parametrization  $\Psi:[0,L_\Gamma) \times (0,\infty) \times [0,2 \pi) \rightarrow \pi_1(M^3_\Gamma\cap M^4), (s,\rho,\phi) \mapsto [x_1,y_1,x_2,y_2]$ of the hypersurfaces $M^3_\Gamma$ outside the zero section
\begin{align*}
M^3_\Gamma \cap M^4&\simeq\mklm{\Big[\rho(u(s)\cos \phi-v(s) \sin \phi), \,\rho(v(s)\cos \phi+u(s) \sin \phi), \,\rho \cos \phi, \,\rho \sin \phi\Big]},
\end{align*}
where $s$ is the length parameter of $\Gamma$ and $\sqrt{\gamma}=\rho\cdot e^{i\phi}\in \C^\ast$ denotes the parameter of the fiber over $\Gamma$. All the following calculations will be performed in $M^3_\Gamma \, \cap M^4$, which is dense in $M^3_\Gamma$. 
The vector fields on $M^3_\Gamma$ induced by the parametrization $\Psi$ read
\begin{align*}
\gd_s=\Psi_\ast\left(\gd_s\right)
&=\Big( \rho(\dot u \cos \phi-\dot v \sin \phi), \, \rho(\dot v \cos \phi + \dot u \sin \phi),\,0,\, 0 \Big), \\
\gd_\rho=\Psi_\ast\left(\gd_\rho\right)&=\Big(u \cos \phi - v \sin \phi, \,v \cos \phi + u \sin \phi, \,\,\cos \phi, \,\sin \phi\Big),  \\
\gd_\phi=\Psi_\ast\left(\gd_\phi\right)&=\Big(-\rho( u \sin \phi+v \cos \phi), \,- \rho(v \sin \phi - u \cos \phi),\,-\rho \sin \phi, \,\rho \,\cos \phi\Big).
\end{align*}
Further, one has
\begin{equation*}
\big|\dot\Gamma(s)\big|^2=\dot u^2+\dot v^2=\dot r^2+r^2\dot \phi_\Gamma^2=1,
\end{equation*}
since $s$ is the arc length parameter of $\Gamma$. Note also that $u\dot u+v \dot v=r \dot r$ and $u\dot v - v\dot u=r^2 \dot \phi_\Gamma$. Moreover, outside the zero section the following identities hold:
\begin{align*}
\nonumber(x_1^2+y_1^2)_{\big|M^3_\Gamma}=\rho^2 r^2,& & &(x_2^2+y_2^2)_{\big|M^3_\Gamma}=\rho^2, \\
\nonumber(x_1y_2-y_1x_2)_{\big|M^3_\Gamma}=-v\rho^2,& & &(x_1x_2+y_1y_2)_{\big|M^3_\Gamma}=u \rho^2,
\end{align*}
and ${u_1}_{|M^3_\Gamma}=\rho^2(r^2+1)$.
Let $h_t$ be the Riemannian metric on $M^3_\Gamma$ induced by $g_t$. In the case of a closed curve $\Gamma$ the hypersurface $(M^3_\Gamma,h_t)$ is a complete Riemannian manifold for $t\not=0$.
Making use of the above relations one obtains the following proposition.

\begin{proposition}
On $M^3_\Gamma \cap M^4$, the coefficients of  the induced Riemannian metric $ h_t$ with respect to the local coordinate frame $\mklm{\gd_s,\gd_\rho,\gd_\phi}$ are given by  
\begin{align*}
h_{11}&=(K+H\rho^2) \rho^2, &
h_{12}&=r \dot r \rho K, \\
h_{13}&=\dot \phi_\Gamma r^2\rho^2 K, &
h_{22}&=(r^2+1)K, \\
h_{23}&=0, & h_{33}&=(r^2+1)\rho^2 K.
\end{align*}
\endproof
\end{proposition}
The function $K$ is given on $M^3_{\Gamma}\cap M^4$ by the formula
\begin{displaymath}
K=\frac{2\rho^2(r^2+1)}{\sqrt{\rho^4(r^2+1)^2+t^4}},
\end{displaymath}
and the functions $G$ and $H$ by
\begin{displaymath}
G:=\frac{2\sqrt{\rho^4(r^2+1)^2+t^4}}{\rho^2(r^2+1)}, \qquad H:=\frac{2t^4}{\rho^4(r^2+1)^2\sqrt{\rho^4(r^2+1)^2+t^4}}. 
\end{displaymath}

In order to compute the relevant geometric quantities of $M^3_\Gamma$  it turns out to be convenient to work within the framework of Cartan. For this purpose we determine an orthonormal frame with respect to $h_t$ by the ansatz
\begin{subequations}
\label{eq:2.1}
\begin{gather}
Y_1=\frac{1}{\sqrt{h_{22}}} \gd_\rho, \qquad \qquad Y_2=\frac{1}{\sqrt{h_{33}}} \gd_\phi,\qquad
Y_3=D \gd_s+ E \gd_\rho+ F \gd_\phi.
\end{gather}
The vector fields $Y_1$ and $Y_2$ are normalized to length $1$; since $h_{23}=0$, they are orthogonal to each other. From the condition $h_t(Y_1,Y_3)=h_t(Y_2,Y_3)=0$ together with $h_t(Y_3,Y_3)=1$ one obtains
\begin{equation}
D=h_{22} \Sigma, \qquad E= -h_{12}\Sigma, \qquad F=-\frac{h_{13}}{\rho^2} \Sigma. 
\end{equation}
\end{subequations}
Here we have introduced the function $\Sigma:= \rho \sqrt{(\det h_t)^{-1}}$ and one computes
\begin{align*}
\det h_t&= h_{22}(h_{11}h_{33}-h_{12}^2\rho^2-h_{13}^2)
        =4(r^2+1) K \rho^4=\frac{8 \rho^6(r^2+1)^2}{\sqrt{\rho^4(r^2+1)^2+t^4}}  > 0.
\end{align*}
The vector fields $\mklm{Y_1,Y_2,Y_3}$ are defined on $M^3_\Gamma\cap M^4$ and outside the exceptional curve do represent  a global section in the frame bundle of $M^3_\Gamma$.
Note that  since $h_{22}$ is positive,  $D$ is always positive.
The local base of $1$--forms $\mklm{\omega^1,\omega^2,\omega^3}$ dual to the orthonormal frame is then given by  
\begin{equation}
\label{eq:2.2}
\omega^1=\sqrt{h_{22}}\Big(d\rho-\frac{E}{D}ds\Big),\qquad 
\omega^2=\sqrt{h_{33}}\Big(d\phi-\frac{F}{D}ds\Big),\qquad
\omega^3=\frac{1}{D}ds,
\end{equation}
and the connection forms $\omega_{ji}=h_t(\nabla Y_j, Y_i)$ of the  Levi--Civita--connection $\nabla$ on $M^3_\Gamma$ as well as  the components 
of the Riemannian curvature tensor   $R^j_{kli}=R_{klij}$
are uniquely determined by Cartan's structure equations
\begin{align}
\label{eq:2.3a}
d\omega^i&=\sum\limits_{j=1}^3\omega_{ij} \wedge \omega^j, \\
\label{eq:2.3b}
d\omega_{ij}&=\sum\limits_{k=1}^3\omega_{ik} \wedge \omega_{kj} + \frac{1}{2} \sum\limits_{k,l=1}^3R^j_{kli} \omega^k \wedge \omega^l.
\end{align}
We determine now the  connection forms of the considered hypersurfaces. 
\begin{proposition}
\label{prop:2.2}
With respect to the orthonormal frame \eqref{eq:2.1}, the forms $\omega_{ji}$ of the Levi--Civita connection on $M^3_\Gamma\cap M^4$ are given by
\begin{align*}
\omega_{12}&=\frac{1}{\sqrt{h_{22}}}\left (\frac{1}{\rho}+\frac{1}{2} (\log K)_{,\rho}\right) \omega^2=\left ( 1 + \frac{\rho}{2}(\log K)_{,\rho} \right)\left (d\phi - \frac{F}{D} ds\right), \\
\omega_{13}&=-\frac{1}{\sqrt{h_{22}}}(\log D)_{,\rho} \,\omega^3= -\frac{1}{D\sqrt{h_{22}}}(\log D)_{,\rho} \,ds,\\
\omega_{23}&=0.
\end{align*}
\end{proposition}
\begin{proof}
Using the orthonormal frame $\mklm{Y_1,Y_2,Y_3}$ on $M^3_\Gamma\cap M^4$, the components of the Levi--Civita--connection can be obtained  via the formulas 
\begin{equation}
\label{eq:2.4}
2h_t(\nabla _{Y_i} {Y_j},{Y_k})=h_t([{Y_i},{Y_j}],{Y_k})-h_t([{Y_j},{Y_k}],{Y_i})+h_t([{Y_k},{Y_i}],{Y_j}),
\end{equation}
resulting from the Koszul--formula.
A direct  computation of the commutators yields  
\begin{align*}
[Y_1,Y_2]&=-\frac{1}{\sqrt{h_{22}}}\left (\frac{1}{\rho}+\frac{1}{2}\frac{K_{,\rho}}{K} \right) Y_2,\\ 
[Y_1,Y_3]&=\left (E_{,\rho}+\frac{D}{2h_{22}} \Big( 2 r\dot r K+ (r^2+1)K_{,s}\Big) +\frac{E}{2h_{22}}(r^2+1) K_{,\rho}-\frac{E}{D} D_{,\rho} \right ) Y_1 \\
&+\rho\left (F_{,\rho}-\frac{F}{D} D_{,\rho}\right ) Y_2+\frac{1}{\sqrt{h_{22}}} \frac{D_{,\rho}}{D} Y_3, \\
[Y_2,Y_3]&=\left(\frac{D}{2h_{33}} \rho^2\Big(2r \dot r K+(r^2+1) K_{,s}\Big)+\frac{E}{2 h_{33}}(r^2+1)(2\rho K+ \rho^2 K_{,\rho})\right ) Y_2.
\end{align*}
Now, a short calculation gives
\begin{equation}
\label{eq:2.5}
(r^2+1) K_{,s}=r \dot r \rho K_{,\rho},
\end{equation}
by which one further calculates
\begin{gather*}
E_{,\rho}-\frac{E}{D} D_{,\rho}=\Sigma \left ( \frac{h_{12}}{h_{22}} (r^2+1) K_{,\rho}- r \dot r K - r \dot r \rho K_{,\rho} \right ) = -\Sigma r \dot r K, \\
{D\Big(2r \dot r K + (r^2+1) K_{,s}\Big)+E(r^2+1) K_{,\rho}}=  2(r^2+1) r \dot r K^2 \Sigma,
\end{gather*}
which shows that the first coefficient of $[Y_1,Y_3]$  vanishes.
Similarly, it can be seen  that the second coefficient is also zero, since 
\begin{align*}
F_{,\rho}-\frac{F}{D} D_{,\rho} &=-h_{13,\rho} \frac{1}{\rho^2} \Sigma - h_{13} \left ( -2 \frac{1}{\rho^3}\Sigma+\frac{1}{\rho^2}\Sigma _{,\rho} \right )\\
&+h_{13} \frac{1}{\rho^2} \Sigma \frac{1}{(r^2+1) K \Sigma}  ( r^2+1) \left ( K_{,\rho} \Sigma + K \Sigma _{,\rho} \right) \\ &= \Sigma \left ( -r^2 \dot \phi_\Gamma(2 \rho K + \rho^2 K_{,\rho})\frac{1}{\rho^2}+2 r^2 \dot \phi_\Gamma\left ( \frac{1}{\rho} K + \frac{1}{2} K_{,\rho} \right ) \right )=0, 
\end{align*}
and by using \eqref{eq:2.5} again one sees that the commutator $[Y_2,Y_3]$ vanishes completely.   In the equations \eqref{eq:2.4} therefore only the terms 
\begin{align*}
h_t([Y_2,Y_1],Y_2)= \frac{1}{\sqrt{h_{22}}}\left (\frac{1}{\rho}+\frac{1}{2} (\log K)_{,\rho}\right), \quad h_t([Y_3,Y_1],Y_3)= -\frac{1}{\sqrt{h_{22}}}(\log D)_{,\rho},
\end{align*}
are non--trivial, and for the forms $\omega_{ij}$ this gives the stated expressions.
\end{proof}
Summing up, one obtains that the structure equations  \eqref{eq:2.3a} read
\begin{align*}
d\omega^1&=\omega_{12}\wedge \omega^2+\omega_{13}\wedge \omega^3=0, \\
d\omega^2&=\omega_{21}\wedge \omega^1+\omega_{23}\wedge \omega^3= \frac{1}{\sqrt{h_{22}}}\left (\frac{1}{\rho}+\frac{1}{2} (\log K)_{,\rho}\right) \omega^1\wedge \omega^2, \\
d\omega^3&=\omega_{31}\wedge \omega^1+\omega_{32}\wedge \omega^2=-\frac{1}{\sqrt{h_{22}}}(\log D)_{,\rho}\,\omega^1\wedge \omega^3.
\end{align*}

We are now able to compute the components of the Riemannian curvature tensor as well as the Ricci tensor and the scalar curvature of the hypersurface $M^3_\Gamma$.

\begin{proposition}
With respect to the section \eqref{eq:2.1}, the components of the Riemannian curvature tensor $R$ of the hypersurfaces $(M^3_\Gamma,h_t)$ are given by
\begin{align*}
R_{1212}&=\frac{1}{2h_{22}}\left (\frac{1}{\rho} (\log K)_{,\rho} + ( \log K)_{,\rho\rho}\right ), \\
R_{2323}&=-\frac{1}{2h_{22}}\left (\frac{2}{\rho} +(\log K)_{,\rho}\right ) (\log D)_{,\rho},\\
R_{1313}&=\frac{1}{2h_{22}} \left (-2(\log D)_{,\rho\rho}+(\log K)_{,\rho} (\log D)_{,\rho}+ 2(\log D)^2_{,\rho}\right ),
\end{align*}
while $R_{1213}$, $R_{2312}$ and $R_{2313}$ vanish.
\end{proposition}
\begin{proof}
We calculate the components of $R$ by using the structure equations \eqref{eq:2.3b} of the hypersurface $M^3_\Gamma$.
By Proposition \ref{prop:2.2} the 2--forms $\omega_{13} \wedge \omega_{32}$ and $\omega_{12} \wedge \omega_{23}$ vanish and
\begin{align*}
\omega_{21}\wedge \omega_{13}&= \frac{1}{h_{22}}\left (\frac{1}{\rho} +\frac{1}{2}(\log K)_{,\rho}\right ) (\log D)_{,\rho} \,\omega^2\wedge \omega^3.
\end{align*}
Further, the differentials $d\omega_{ij}$ of the connection forms are given by
\begin{align*}
d\omega_{12}&=\left (\frac{\rho}{2} (\log K)_{,\rho} \right )_{,\rho} d\rho \wedge d\phi+\left (\frac{\rho}{2} (\log K)_{,\rho} \right )_{,s} ds \wedge d\phi\\ &-\left ( \frac{F}{D}\left( 1+\frac{\rho}{2} (\log K)_{,\rho}\right ) \right )_{,\rho} d\rho \wedge ds \\
&=\left (\frac{1}{2} (\log K)_{,\rho} + \frac{\rho}{2} ( \log K)_{,\rho\rho}\right )\frac{1}{h_{22}\rho} \,\omega^1\wedge \omega^2 \\ &-\left (\left ( \frac{1}{2} (\log K)_{,\rho}+\frac{\rho}{2} (\log K)_{,\rho\rho}\right )\frac{E}{\sqrt{h_{33}}} +\frac{\rho}{2}(\log K)_{,s\rho}\frac{D}{\sqrt{h_{33}}}\right ) \omega^2 \wedge \omega^3, \\ 
d\omega_{13}&=-\left (\frac{1}{D\sqrt{h_{22}}} (\log D)_{,\rho} \right )_{,\rho} d\rho \wedge ds\\
&=\left (-\frac{1}{{h_{22}}} (\log D)_{,\rho\rho}+ \frac{1}{2{h_{22}}} (\log K)_{,\rho} (\log D)_{,\rho}+ \frac{1}{{h_{22}}} (\log D)^2_{,\rho}\right )\omega^1\wedge \omega^3, \\
d\omega_{23}&=0,
\end{align*}
and one obtains the stated formulas for the components $R_{ijkl}$ of the curvature tensor by using \eqref{eq:2.3b}. Notice that
\begin{eqnarray*}
\lefteqn{\left ((\log K)_{,\rho}+\rho (\log K)_{,\rho\rho}\right )E +\rho(\log K)_{,s\rho}D} \\
& &=\left ((\log K)_{,\rho}+\rho (\log K)_{,\rho\rho}\right ) (-r\dot r \rho K \Sigma)+\rho \left ( \frac {\rho r \dot r}{r^2+1} ( \log K)_{,\rho} \right )_{,\rho}(r^2+1) K \Sigma=0, 
\end{eqnarray*}
implying that  $R_{2312}$ vanishes.
\end{proof}

\begin{theorem}
\label{thm:2.1}
The components $R_{ij}$ of the Ricci tensor $\Ric$ of the Riemannian $C^\infty$--manifolds $(M^3_\Gamma,h_t)$ are given with respect to the orthonormal frame \eqref{eq:2.1}   by
\begin{equation*}
\Ric=\frac{1}{2\Big(\rho^4(r^2+1)^2 +t^4\Big)^{3/2}}\left (\begin{array}{ccc} 2t^4 & 0 & 0 \\ 0 & 2t^4-\rho^4(r^2+1)^2 & 0 \\ 0 & 0 & -4t^4-\rho^4(r^2+1)^2 \end{array} \right ),
\end{equation*}
and the scalar curvature is
\begin{equation*}
S=R_{11}+R_{22}+R_{33}=-\frac{\rho^4(r^2+1)^2}{\Big(\rho^4(r^2+1)^2 +t^4\Big)^{3/2}}.
\end{equation*}
\end{theorem}
\begin{proof}
 One computes 
\begin{align*}
(\log K)_{,\rho}& = \frac{2t^4}{\Big(\rho^4(r^2+1)^2+t^4\Big)\rho}, \\
(\log K)_{,\rho\rho}& = -\frac{2t^4}{\Big(\rho^4(r^2+1)^2+t^4\Big)\rho^2}-\frac{8t^4 \rho^2 (r^2+1)^2}{\Big(\rho^4(r^2+1)^2 +t^4\Big)^2}, \\
(\log \Sigma)_{,\rho}& =-\frac{1}{\rho}-\frac{t^4}{\Big(\rho^4(r^2+1)^2+t^4\Big)\rho}= -\frac{1}{\rho}-\frac{1}{2}(\log K)_{,\rho}, \\
(\log D)_{,\rho}&= (\log K)_{,\rho}+(\log \Sigma)_{,\rho}=-\frac{1}{\rho}+\frac{1}{2}(\log K)_{,\rho}, 
\end{align*}
obtaining thus, with the previous proposition, for the components $R_{ij}=\sum_{k=1}^3 R_{ikkj}$ of the Ricci tensor that 
\begin{align*}
R_{11}&= \frac{1}{2h_{22}}\left [ -\frac{1}{\rho} ( \log K)_{,\rho} - (\log K )_{,\rho\rho}-(\log K)_{,\rho}\left ( -\frac{1}{\rho}+\frac{1}{2} (\log K )_{, \rho} \right )\right. \\
 &\left. -2 \left ( \frac{1}{\rho^2} -\frac{1}{\rho}(\log K)_{,\rho} +\frac{1}{4} (\log K )^2_{,\rho} \right ) +2\left ( \frac{1}{\rho^2} +\frac{1}{2} (\log K )_{,\rho\rho}\right ) \right ]\\
 &=\frac{1}{2h_{22}} \left ( -(\log K)^2_{,\rho} +2 \frac{1}{\rho} (\log K ) _{,\rho}\right ),\\
R_{22}&= \frac{1}{2h_{22}}\left [ -\frac{1}{\rho} ( \log K)_{,\rho} - (\log K )_{,\rho\rho}+\frac{1}{\rho^2} \Big( 2+ \rho \,(\log K)_{,\rho}\Big) \left ( -{1}+\frac{\rho}{2} (\log K )_{, \rho} \right )\right ] \\
 &=\frac{1}{2h_{22}} \left ( -\frac{1}{\rho} ( \log K)_{,\rho} - (\log K )_{,\rho\rho}-\frac{2}{\rho^2}+ \frac{1}{2} ( \log K)^2_{,\rho} \right ),\\
R_{33}&=\frac{1}{2h_{22}}\left [-(\log K)_{,\rho}\left ( -\frac{1}{\rho}+\frac{1}{2} (\log K )_{, \rho} \right )-2 \left ( \frac{1}{\rho^2} -\frac{1}{\rho}(\log K)_{,\rho} +\frac{1}{4} (\log K )^2_{,\rho} \right ) \right. \\
& \left.+ 2\left ( \frac{1}{\rho^2} +\frac{1}{2} (\log K )_{,\rho\rho}\right )+\frac{1}{\rho} \Big( 2+ \rho (\log K)_{,\rho}\Big) \left ( -\frac{1}{\rho}+\frac{1}{2} (\log K )_{, \rho} \right ) \right ]\\
& =\frac{1}{2h_{22}} \left ( ( \log K ) _{,\rho\rho} - \frac{1}{2}( \log K ) ^2 _{,\rho} + \frac{3}{\rho}(\log K )_{,\rho}-\frac{2}{\rho^2} \right ),
\end{align*}
the remaining coefficients being equal to zero. In the same way as the components of the Riemannian curvature tensor turn out to be bounded when $\rho\to 0$, the components of $\Ric$ and thus $S$ stay bounded, too. Explicitly one has  
\begin{align*}
R_{11}&=\frac{1}{\rho^2 (r^2+1)^2 \Big( \rho^4 (r^2+1)^2 +t^4\Big)^{3/2}} \left ( \frac{t^4}{\rho^2} \Big( \rho^4(r^2+1)^2 + t^4 \Big)-\frac{t^8}{\rho^2} \right )\\
&=\frac{t^4}{\Big(\rho^4(r^2+1)^2 +t^4\Big)^{3/2}},
\intertext{}
R_{22}&=\frac{1}{2\rho^2 (r^2+1)^2 \Big( \rho^4 (r^2+1)^2 +t^4\Big)^{3/2}} \left ( -\frac{t^4}{\rho^2}\Big(\rho^4(r^2+1)+t^4\Big) + \frac{t^8}{\rho^2}  \right. \\
&\left. +\frac{t^4}{\rho^2}\Big(\rho^4(r^2+1)^2 +t^4\Big) + 4 t^4 \rho ^2 (r^2+1)^2- \frac{1}{\rho^2} \Big( \rho^8(r^2+1)^4 +t^8 +2 \rho^4(r^2+1)^4\Big) \right)\\&=\frac{2t^4-\rho^4 (r^2+1)^2}{2\Big(\rho^4(r^2+1)^2 +t^4\Big)^{3/2}},\\
R_{33}&=\frac{1}{2\rho^2 (r^2+1)^2 \Big( \rho^4 (r^2+1)^2 +t^4\Big)^{3/2}} \left (-\frac{t^4}{\rho^2}\Big(\rho^4(r^2+1)^2 +t^4\Big) -4t^4\rho^2(r^2+1)^2 \right. \\
&-\left. \frac{t^8}{\rho^2}+\frac{3t^4}{\rho^2} \Big(\rho^4(r^2+1)^2 +t^4\Big)-\frac{1}{\rho^2}\Big(\rho^8(r^2+1)^4+t^8+2 \rho^4(r^2+1)^2t^4\Big) \right )\\
&=\frac{-4t^4-\rho^4(r^2+1)^2}{2\Big(\rho^4(r^2+1)^2 +t^4\Big)^{3/2}},
\end{align*}
showing that the divergent terms cancel out and the assertion follows.
\end{proof}

Hence, the scalar curvature $S$ is  negative  and tends  as ${1}/{\rho^2(r^2+1)}$ to zero as $\rho$ and $r$ go to infinity. For  $t\not= 0$ all components of the Riemannian and Ricci tensor as well as $S$ remain regular at $\rho=0$ and are therefore defined everywhere  on the Riemannian manifolds $M^3_\Gamma$. For  $t=0$ the scalar curvature degenerates at $\rho=0$ in concordance with the fact that the hypersurfaces $M_\Gamma^3$ are  no longer complete in this case.

\section{The second fundamental form of the hypersurfaces $M^3_\Gamma$}
\label{sec:3}

We proceed now  studying the  second fundamental form  of the hypersurfaces $M^3_\Gamma$.
In order to do so, we need  the Levi--Civita connection $\nabla^{H^2}$ of the Eguchi--Hanson space $(H^2, g_t)$. It can be obtained from the Koszul formula, which reads for commuting vector fields as follows:
\begin{displaymath}
2 g_t(\nabla^{H^2}_X Y, Z)=X(g_t(Y,Z))+Y(g_t(Y,Z))-Z(g_t(X,Y)).
\end{displaymath}
In the following we will denote the coordinates $x_1,y_1,x_2,y_2$ of the dense complex manifold $M^4\subset H^2$ by $x_1,x_2,x_3,x_4$ so that the components of $\nabla^{H^2}$ on $M^4$ are given by
\begin{gather}
\label{eq:3.0}
\Gamma_{ijk}=\frac{1}{2} \left ( \frac{\gd g_{ik}}{\gd x_j}+\frac{\gd g_{jk}}{\gd x_i}-\frac{\gd g_{ij}}{\gd x_k}\right ),\qquad
\Gamma_{ij}^k=g^{kl}\Gamma_{ijl}.
\end{gather}
Because of the symmetry
\begin{equation}
\label{eq:3.1}
g_{kl}(x_m)=g_{\bar k\bar l}(x_{\overline m}), \qquad \bar k:=k+2 \mod 4,
\end{equation}
of the covariant coefficients of the metric one obtains, for the Christoffel symbols, the relations
\begin{equation}
\label{eq:3.2}
\Gamma_{klm}(x_n)=\Gamma_{\bar k\bar l\overline m}(x_{\bar n}).
\end{equation}
Now the contravariant coefficients $g^{ij}$ of $g_t$ are given by the matrix
\begin{displaymath}
g_t^{-1}=\frac{1}{\sqrt{\det g_t}} \left (\begin{array}{cccc} G_2 & 0& G_4 & G_3 \\ 0 & G_2 &-G_3 & G_4 \\ G_4 & -G_3 & G_1 & 0 \\ G_3 & G_4 & 0 & G_1 \end{array} \right ),  
\end{displaymath}
where $\det g_t=(G_1G_2-G_3^2-G_4^2)^2=16$.
Further, the derivatives of the functions  $G$ and $H$ are
\begin{displaymath}
G_{,x_i}=-2x_iH, \qquad H_{,x_i}=-2x_iI,
\end{displaymath}
with $I:={2t^4}(3u_1^2+t^4)/{u_1^3(u_1^2+t^4)^{3/2}}$. A straightforward calculation yields the Christoffel symbols of the first kind. 

\begin{proposition}
The Christoffel symbols of the first kind of the Eguchi--Hanson space $(H^2,g_t)$  are given on $M^4$ by
\begin{alignat*}{2}
\Gamma_{111}&=-x_1(2H-I(x_1^2+x_2^2)),&\Gamma_{114}&=-2I\bigg(x_1x_2x_3+\frac{1}{2}x_4(x^2_2-x_1^2)\bigg),\\ \Gamma_{112}&=x_2(2H-I(x_1^2+x_2^2)),& \Gamma_{131}&=-x_3(H-I(x_1^2+x^2)), \\ 
\Gamma_{113}&=2I\bigg(x_1x_2x_4+\frac{1}{2}x_3(x_1^2-x_2^2)\bigg),\quad&\Gamma_{132}&=x_4(H-I(x_1^2+x_2^2));
\end{alignat*}
the remaining ones can be obtained from these by taking into account the symmetry $\Gamma_{ijk}=\Gamma_{jik}$ as well as the relations \eqref{eq:3.2} together with the additional symmetries
\begin{align*}
\Gamma_{22i}&=-\Gamma_{11i}, & \Gamma_{44i}&=-\Gamma_{33i}, &
\Gamma_{23i}&=\Gamma_{14i}, & \Gamma_{24i}&=-\Gamma_{13i} 
\end{align*}
and
\begin{alignat*}{4}
\Gamma_{12i}&=\Gamma_{11(i-1)},\qquad & \Gamma_{14i}&=\Gamma_{13(i-1)},\qquad & \Gamma_{34i}&=\Gamma_{33(i-1)}&\qquad&\text{for $i$ even,}\\
\Gamma_{12i}&=-\Gamma_{11(i+1)}, & \Gamma_{14i}&=-\Gamma_{13(i+1)},& \Gamma_{34i}&=-\Gamma_{33(i+1)}&\qquad&\text{for $i$ odd}.
 \end{alignat*}
\endproof
\end{proposition}

The Christoffel symbols of the second kind are derived from these formulas as indicated in \eqref{eq:3.0}. By the symmetries of the Levi--Civita connection $\nabla ^{H^2}$ it is sufficient to compute only six of them explicitly. So one has
\begin{align*}
\Gamma_{11}^1%&=\frac{1}{4}\left [ G_2\Gamma_{111}+G_4\Gamma_{113}+G_3\Gamma_{114}\right ]\\
 &=-\frac{1}{4} x_1G(2H-I(x_1^2+x_2^2))+\frac{1}{4}H\bigg [x_1(x_3^2+x_4^2)(2H-I(x_1^2+x_2^2))\\
&+2I\bigg[(x_1x_3+x_2x_4)\bigg(x_1x_2x_4+\frac{1}{2}x_3(x_1^2-x_2^2)\bigg)\\&-(x_1x_4-x_2x_3)\bigg(x_1x_2x_3+\frac{1}{2}x_4(x_2^2-x_1^2)\bigg)\bigg]\bigg]\\
&=x_1\bigg[-\frac{1}{4}G(2H-I(x_1^2+x_2^2))+\frac{1}{2}H^2(x_3^2+x_4^2)+\frac{1}{4}HI\Big[-(x_1^2+x_2^2)(x^2_3+x^2_4)\\
&+(x_3^2+x_4^2)(x_1^2-x_2^2)+2(x_2^2x_4^2+x_2^2x_3^2)\Big] \bigg]=x_1 A_1,
\end{align*}
where we have introduced
$
A_1=\frac{1}{2}H^2(x_3^2+x_4^2)-\frac{1}{4}G(2H-I(x_1^2+x_2^2)).
$
In a similar way one obtains
\begin{align*}
\Gamma_{11}^2%&=\frac{1}{4}\left [ G_2\Gamma_{112}-G_3\Gamma_{113}+G_4\Gamma_{114}\right ]\\
 &=-x_2 A_1.
\end{align*}
Further one checks that
\begin{align*}
\Gamma_{11}^3%&=\frac{1}{4}(G_4\Gamma_{111}-G_3\Gamma_{112}+G_1\Gamma_{113})\\
&=\frac{1}{4}G\Gamma_{113}+\frac{1}{4}H\Big [-\Gamma_{113}(x_1^2+x_2^2)\\
&-(2H-I(x_1^2+x_2^2))\big [x_1(x_1x_3+x_2x_4)+x_2(x_1x_4-x_2x_3) \big] \Big]\\
&=\frac{1}{4}G\Gamma_{113}-\frac{1}{2}H^2(2x_1x_2x_4+x_3(x_1^2-x_2^2))\\
&+\frac{1}{4}H\Big[-\Gamma_{113}(x_1^2+x_2^2)+I(x_1^2+x_2^2)(2x_1x_2x_4+x_3(x_1^2-x_2^2))\Big]\\
&=\frac{1}{2}C(x_1x_2x_4+\frac{1}{2}x_3(x_1^2-x_2^2)),
\end{align*}
as well as
\begin{align*}
\Gamma_{11}^4%&=\frac{1}{4}(G_3\Gamma_{111}+G_4\Gamma_{112}+G_1\Gamma_{114})\\
&=-\frac{1}{2}C(x_1x_2x_3+\frac{1}{2}x_4(x_2^2-x_1^2)),
\end{align*}
where $C=(IG-2H^2)$.
Finally, one calculates
\begin{align*}
\Gamma_{13}^1%&=\frac{1}{4}(G_2\Gamma_{131}+G_4\Gamma_{133}+G_3\Gamma_{134})\\
&=\frac{1}{4}G\Gamma_{131}+\frac{1}{4}H\Big [-\Gamma_{131}(x_3^2+x_4^2)\\
&+(H-I(x_3^2+x_4^2))\big [-x_1(x_2x_4+x_1x_3)+x_2(x_1x_4-x_2x_3) \big] \Big]\\
&=-\frac{1}{4}x_3\Big[ G(H-I(x_1^2+x_2^2))+H^2(-x_3^2-x_4^2+x_1^2+x_2^2)\Big]\\
&+\frac{1}{4}HI\Big [x_3(x_1^2+x_2^2)(x_3^2+x_4^2)-x_3(x_3^2+x_4^2)(x_1^2+x_2^2)\Big ]=-x_3B_1,\\
\Gamma_{13}^2&=x_4B_1,
\end{align*}
with $B_1$ given by
$
B_1=\frac{1}{4} \big [G(H-I(x_1^2+x_2^2))+H^2(x_1^2+x_2^2-x_3^2-x_4^2)\big ].
$
Taking into account the relations \eqref{eq:3.1}, \eqref{eq:3.2} one thus obtains that the $\Gamma^i_{jk}$ are given as follows.
\begin{proposition}
The components $\Gamma^i_{jk}$ of the Levi--Civita connection of the Eguchi--Hanson space $(H^2,g_t)$ are given on $M^4$ by 
\begin{alignat*}{2}
\Gamma_{11}^1&=x_1A_1,&\Gamma_{33}^3&=x_3A_2,\\
\Gamma_{11}^2&=-x_2A_1,&\Gamma_{33}^4&=-x_4A_2,\\
\Gamma_{11}^3&=\frac{C}{2}\bigg(x_1x_2x_4+\frac{1}{2}x_3(x_1^2-x_2^2)\bigg),\qquad &\Gamma_{33}^1&=\frac{C}{2}\bigg(x_2x_3x_4+\frac{1}{2}x_1(x_3^2-x_4^2)\bigg),\\
\Gamma_{11}^4&=-\frac{C}{2}\bigg (x_1x_2x_3+\frac{1}{2}x_4(x_2^2-x_1^2)\bigg),\qquad &\Gamma_{33}^2&=-\frac{C}{2}\bigg(x_1x_3x_4+\frac{1}{2}x_2(x_4^2-x_3^2)\bigg),\\
\Gamma_{13}^1&=-x_3B_1,&\Gamma_{13}^3&=-x_1B_2,\\
\Gamma_{13}^2&=x_4B_1,&\Gamma_{13}^4&=x_2B_2,\\
\end{alignat*}
where 
\begin{align*}
A_1&=\frac{1}{2}H^2(x_3^2+x_4^2)-\frac{1}{4}G(2H-I(x_1^2+x_2^2)),\\
B_1&=\frac{1}{4} \Big [G(H-I(x_1^2+x_2^2))+H^2(x_1^2+x_2^2-x_3^2-x_4^2)\Big ],\\
A_2&=\frac{1}{2}H^2(x_1^2+x_2^2)-\frac{1}{4}G(2H-I(x_3^2+x_4^2)),\\
B_2&=\frac{1}{4} \Big [G(H-I(x_3^2+x_4^2))+H^2(x_1^2+x_2^2-x_3^2-x_4^2)\Big ]
\end{align*}
and $C=(IG-2H^2)$. All remaining $\Gamma^i_{jk}$ can be obtained from the above by using $\Gamma_{ij}^k=\Gamma_{ji}^k$ as well as the relations
\begin{align*}
\Gamma_{22}^i&=-\Gamma_{11}^i, & \Gamma_{44}^i&=-\Gamma_{33}^i,&
\Gamma_{23}^i&=\Gamma_{14}^i, & \Gamma_{24}^i&=-\Gamma_{13}^i 
\end{align*}
and
\begin{alignat*}{4}
\Gamma_{12}^i&=\Gamma_{11}^{(i-1)},\qquad & \Gamma_{14}^i&=\Gamma_{13}^{(i-1)},\qquad & \Gamma_{34}^i&=\Gamma_{33}^{(i-1)}&\qquad&\text{for $i$ even,}\\
\Gamma_{12}^i&=-\Gamma_{11}^{(i+1)}, & \Gamma_{14}^i&=-\Gamma_{13}^{(i+1)},& \Gamma_{34}^i&=-\Gamma_{33}^{(i+1)} &\qquad&\text{for $i$ odd}.
\end{alignat*}
\endproof
\end{proposition}

In order to describe the outer geometry of the hypersurfaces $M^3_\Gamma$, we first determine a field of unit  normal vectors $N:M^3_\Gamma \rightarrow (TM^3_\Gamma)^\perp$ on $M^3_\Gamma$. Up to orientation such a field is given by the conditions
\begin{displaymath}
g_t(Y_i,N)=0, \quad i=1,2,3, \qquad g_t(N,N)=1,
\end{displaymath}
which are equivalent to the system of equations
\begin{gather*} 
N_1(u \cos \phi-v \sin \phi)K+N_2(v \cos \phi + u \sin \phi) K  
+N_3K \cos \phi+N_4K \sin \phi=0,\\[.3em]
-N_1(u \sin \phi +v \cos \phi)K-N_2(v \sin \phi -u \cos \phi) K
-N_3K \sin \phi +N_4 K \cos \phi=0,\\[.3em]
\Big\{N_1(\dot u \cos \phi-\dot v \sin \phi)+N_2(\dot v \cos \phi+\dot u \sin \phi)\Big\}(K+H\rho^2)-N_3(\cos \phi(v\dot v+u\dot u)\\+ \sin \phi(\dot u v -\dot v u)) H\rho^2-N_4(\sin \phi(v\dot v +u \dot u)-\cos \phi ( \dot u v - \dot vu)) H \rho^2=0,\\[.3em]
(N_1^2+N_2^2)(K+H\rho^2)+(N_3^2+N_4^2)(K+H\rho^2r^2)+2(N_4N_1-N_3N_2)v\rho^2 H\\
-2(N_3N_1+N_4N_2)u \rho^2H=0.
\end{gather*}
By solving these equations with respect to the components $N_i$ of the unit normal vectors one obtains the following proposition.
\begin{proposition}
On $M^3_\Gamma\cap M^4$ a field of unit normal vectors is given by 
\begin{displaymath}
N=\frac{1}{2} \sqrt{\frac {K}{r^2+1}} (w_1,-w_2,-rw_3,rw_4),
\end{displaymath}
where the functions $w_i$  are
\begin{align*}
 w_1&=\dot v \cos \phi + \dot u \sin \phi, & w_2&=\dot u \cos \phi - \dot v \sin \phi,\\
w_3&=\dot r \sin \phi  +r\dot \phi_\Gamma \cos \phi,& w_4&=\dot r \cos \phi- r\dot \phi_\Gamma \sin \phi.
\end{align*}
\endproof
\end{proposition}

We note that $w_2$ and  $w_1$ can be viewed as the real and imaginary part of $\dot \Gamma(s)e^{i\phi}$, $w_4$ and $w_3$ as the real and  imaginary part of $\dot \Gamma(s)e^{i\phi}e^{-i\phi_\Gamma}$ respectively.
%We remark further, that $(Y_1,Y_2,Y_3,N)$ represents a so--called adapted frame of the Levi--flat real hypersurface $M^3_\Gamma$.
By construction the hypersurfaces $M^3_\Gamma$ are imbedded in $H^2$. If $N$ denotes the field of unit normal vectors determined above, the second fundamental form of $M^3_\Gamma$ is defined by
\begin{gather}
\label{eq:3.3}
{\II:\X(M^3_\Gamma\cap M^4) \times \X(M^3_\Gamma\cap M^4) \longrightarrow \F(M^3_\Gamma),}\qquad \II(X,Y)=g_t(X,\nabla^{H^2}_Y N).
\end{gather}
It is symmetric and bilinear. In the following we will write the coordinates $s, \rho, \phi$ as $\eta_1,\eta_2,\eta_3$, and denote the components of $\II$ with respect to the induced frame of coordinate vector fields by $\II_{ij}$. For shortness, we will simply write $\nabla$ for $\nabla^{H^2}$ in the remaining of this section. Explicitly,
\begin{gather}
\label{eq:3.4}
\begin{split}
\nabla_{\gd_{\eta_j}} N(p)=\Big(\gd_{\eta_j} N_i(p) +N_i(p)\, \nabla_{\gd_{\eta_j}}\Big)\,\gd_{x_i}|_p,\\
\nabla_{\gd_{\eta_j}}\gd_{x_i}|_p=dx_k(\gd_{\eta_j})(p)\,\Gamma^l_{ki}(p)\,\gd_{ x_l}|_p,
\end{split}
\end{gather}
where  $p \in M^3_\Gamma\cap M^4$. On $M^3_\Gamma\cap M^4$ the coordinates $x_1,x_2,x_3,x_4$ can be expressed  by the coordinates  $s,\phi,\rho$ according to
\begin{align*}
x_1&=\rho \,\zeta_1,& x_2&=\rho \,\zeta_2, &
x_3&=\rho \,\cos \phi, & x_4&=\rho\,\sin \phi,
\end{align*}
where we have defined
\begin{displaymath}
\zeta_1=u\cos\phi-v \sin \phi, \qquad \zeta_2=v\cos\phi+u\sin\phi.
\end{displaymath}
Thus one obtains that  on  $M^3_\Gamma\cap M^4$ the polynomials appearing in the expressions for the $\Gamma^{i}_{jk}$ are given by 
\begin{align*}
\bigg [x_1x_2x_4&+\frac{1}{2}x_3(x_1^2-x_2^2)\bigg]|_{M^3_\Gamma}= 
 \rho^3\bigg[-uv\sin\phi+(u^2-v^2)\frac{1}{2} \cos \phi\bigg ],\\
\bigg[x_1x_2x_3&+\frac{1}{2}x_4(x_2^2-x_1^2)\bigg]|_{M^3_\Gamma}=\rho^3\bigg[uv\cos\phi+(u^2-v^2)\frac{1}{2} \sin \phi\bigg ],\\
\bigg[x_2x_3x_4&+\frac{1}{2}x_1(x_3^2-x_4^2)\bigg]|_{M^3_\Gamma}=\frac{\rho^3}{2}[v\sin \phi+u\cos \phi],\\
\bigg[x_1x_3x_4&+\frac{1}{2}x_2(x_4^2-x_3^2)\bigg]|_{M^3_\Gamma}=\frac{\rho^3}{2}[u\sin \phi-v\cos \phi].
\end{align*}

We compute now the covariant derivatives $\nabla _{\gd_{\eta_i}}N$. To this end we first note the relations
\begin{align}
\label{eq:3.5}
\nonumber w_3\cos \phi -w_4\sin \phi=r\dot \phi_\Gamma,& & &w_4\cos \phi +w_3\sin \phi =\dot r,\\
w_2\cos \phi +w_1\sin \phi=\dot u,& & &w_1\cos \phi-w_2\sin \phi=\dot v,\\
\nonumber \zeta_1w_2+w_1\zeta_2=r\dot r,& & &\zeta_1w_1-\zeta_2w_2=r^2\dot \phi_\Gamma.
\end{align}
Because of
\begin{align*}
r\dot u&=r(\dot r \cos \phi_\Gamma-r\dot \phi_\Gamma \sin \phi_\Gamma)=\dot r u -r \dot \phi_\Gamma v,\\
r\dot v&=r(\dot r \sin \phi_\Gamma+r\dot \phi_\Gamma \cos \phi_\Gamma)=\dot r v +r \dot \phi_\Gamma u
\end{align*}
one has further
\begin{subequations}
\label{eq:3.6}
\begin{align}
rw_1&=\cos \phi(\dot r v+r\dot \phi_\gamma u)+\sin \phi(\dot r u - r \dot \phi_\Gamma v)=uw_3+vw_4, \\
rw_2&=\cos \phi(\dot r u-r\dot \phi_\gamma v)-\sin \phi(\dot r v + r \dot \phi_\Gamma u)=-vw_3+uw_4,
\intertext{and thus}
rw_3&=\cos \phi_\Gamma (uw_3+vw_4)-\sin \phi_\Gamma(-vw_3+uw_4) =uw_1-vw_2,\\
rw_4&=\sin \phi_\Gamma (uw_3+vw_4)+\cos \phi_\Gamma(-vw_3+uw_4)=w_1 v+w_2u.
\end{align}
\end{subequations}
\begin{proposition}
\label{prop:3.3}
\begin{align*}
\nabla_{\gd_{s}}N&=\sum\limits^{4}_{i=1}(N_{i,s}+\Phi(\dot r N_i+r\dot \phi_\Gamma N_{i,\phi}))\gd_{x_i}= N_{,s}+\Phi(\dot r N +r \dot \phi_\Gamma N_{,\phi}),
\end{align*}
where $\Phi=-\frac{1}{4} HK \rho^2r.$
\end{proposition}
\begin{proof}
By using the symmetries of the  Christoffel symbols $\Gamma^k_{ij}$ one has
\begin{align*}
\nabla_{\gd_{s}}N%&=\left(N_{i,s}+N_k dx_j\left (\gd_{s}\right)\Gamma^i_{jk} \right) \gd_{x_i}\\
&=\sum\limits_{i=1}^4 \Big[N_{i,s}+\sum\limits_{k=1}^4 N_k\left ( \Gamma^i_{1k} \rho w_2+\Gamma_{2k}^i \rho w_1\right) \Big] \gd_{x_i}\\
&=\sum\limits_{i=1}^4\Big [N_{i,s}+\rho \Gamma^i_{11}(N_1w_2-N_2w_1)+\rho \Gamma^i_{12}(N_2w_2+N_1w_1) \\
&+\rho \Gamma^i_{13}(N_3 w_2-N_4w_1)+\rho \Gamma^i_{14}(N_4w_2+N_3w_1)\Big ] \gd_{x_i}\\
&=\sum\limits_{i=1}^4\bigg [N_{i,s}+\rho\frac{1}{2}\sqrt{\frac{K}{r^2+1}} ( 2\Gamma^i_{11} w_1w_2+\Gamma^i_{12}(w_1^2-w_2^2) \\
& +\Gamma^i_{13}r(-w_3w_2-w_4w_1)+\Gamma^i_{14}r(w_4w_2-w_3w_1))\bigg]\gd_{x_i}.
\end{align*}
The first component of $\nabla_{\gd_{s}}N$ reads
\begin{displaymath}
dx_1(\nabla_{\gd_{s}}N)=N_{1,s}+\frac{1}{2} \sqrt{\frac{K}{r^2+1}}\rho^2r(A_1+B_1)(\dot r w_1 + r \dot \phi_\Gamma w_2),
\end{displaymath}
since by the relations \eqref{eq:3.5} one has
\begin{align*}
2\Gamma^1_{11}w_1w_2+\Gamma^1_{12}&(w_1^2-w_2^2)=A_1\rho(2\zeta_1w_1w_2+\zeta_2(w_1^2-w_2^2))\\
&=A_1\rho [w_1(\zeta_1w_2+w_1\zeta_2)+w_2(\zeta_1w_1-\zeta_2w_2)]\\&=A_1\rho r (\dot r w_1+r \dot \phi_\Gamma w_2),\\
\Gamma^1_{13}r(-w_3w_2-w_4&w_1)+\Gamma^1_{14}r(w_4w_2-w_3w_1)
\\&=B_1\rho r [\cos \phi(w_2w_3+w_1w_4)-\sin \phi(w_4w_2-w_3w_1)]\\
&=B_1\rho r [w_2(w_3\cos \phi -w_4\sin \phi ) +w_1(w_4\cos \phi +w_3\sin \phi )]\\&=B_1\rho r (w_2r \dot \phi _\Gamma+w_1 \dot r).
\end{align*}
The second component is given by
\begin{displaymath}
dx_2(\nabla_{\gd_{s}}N)=N_{2,s}+\frac{1}{2} \sqrt{\frac{K}{r^2+1}}\rho^2r(A_1+B_1)(-\dot r w_2 + r \dot \phi_\Gamma w_1),
\end{displaymath}
as can be verified by an analogous calculation.
As far as the third component is concerned, using also the relations \eqref{eq:3.6} one computes
\begin{align*}
2\Gamma^3_{11}&w_1w_2+\Gamma^3_{12}(w_1^2-w_2^2)=\frac{1}{2} \rho^3 C \bigg[ \bigg (-uv \sin \phi + \frac{u^2-v^2}{2}\cos \phi\bigg )2 w_1w_2 \\
&+\bigg ( uv \cos \phi+\frac{u^2-v^2}{2} \sin \phi\bigg )(w_1^2-w_2^2)\bigg]
\end{align*}
\begin{eqnarray*}
&&=\frac{1}{2} \rho^3 C\bigg [uv\Big(w_1(-w_2 \sin \phi +w_1\cos \phi)+w_2(-w_1 \sin \phi-w_2 \cos \phi)\Big)\\
&&+\frac{u^2-v^2}{2}\Big (w_1(w_2 \cos \phi +w_1 \sin \phi)+ w_2(w_1 \cos \phi -w_2\sin \phi)\Big) \bigg]\\ 
&&=\frac{1}{2} \rho^3 C\bigg[ uv (w_1 \dot v-w_2 \dot u)+\frac{u^2-v^2}{2}(w_1\dot u +w_2 \dot v)\bigg ]\\
&&=\frac{1}{2} \rho^3 C\bigg[\frac{vw_1}{2}(\dot v u - \dot u v)+\frac{uw_2}{2} (-\dot u v+\dot v u )+\frac{uw_1}{2}(v\dot v+u \dot u)-\frac{vw_2}{2}(u\dot u+v\dot v) \bigg]\\
&&=\frac{1}{4} \rho^3 C[(vw_1+uw_2)r^2 \dot \phi_\Gamma+(uw_1-vw_2)r\dot r ]=\frac{1}{4} \rho^3r^2 C(r\dot \phi_\Gamma w_4+\dot r w_3)
\end{eqnarray*}
and
\begin{eqnarray*}
\lefteqn{\Gamma^3_{13}r(-w_3w_2-w_4w_1)+\Gamma^3_{14}r(w_4w_2-w_3w_1)}\\
&&=B_2\rho r [-\zeta_1(-w_2w_3-w_1w_4)-\zeta_2(w_4w_2-w_3w_1)]\\
&&=B_2\rho r [w_3(\zeta_1w_2+\zeta_2w_1)+w_4(\zeta_1w_1-\zeta_2w_2)]=B_2\rho r^2 (w_3 \dot r+w_4 r\dot \phi_\Gamma),
\end{eqnarray*}
so that
\begin{displaymath}
dx_3(\nabla_{\gd_{s}}N)=N_{3,s}+\frac{1}{2} \sqrt{\frac{K}{r^2+1}}\rho^2r^2\bigg(\frac{1}{4}\rho^2C+B_2\bigg)(\dot r w_3 + r \dot \phi_\Gamma w_4).
\end{displaymath}
In the same way one verifies for the fourth component that 
\begin{displaymath}
dx_4(\nabla_{\gd_{s}}N)=N_{4,s}+\frac{1}{2} \sqrt{\frac{K}{r^2+1}}\rho^2r^2\bigg(\frac{1}{4}\rho^2C+B_2\bigg)(-\dot r w_4 + r \dot \phi_\Gamma w_3).
\end{displaymath}
The stated expression for $\nabla_{\gd_{s}}N$ then follows by noting that the equalities
\begin{align*}
A_1+B_1&=\frac{1}{4}(H^2\rho^2(r^2+1)-GH)=-\frac{1}{4}HK,\\
\frac{1}{4}C\rho^2+B_2&=\frac{1}{4}(IG-2H^2)\rho^2+\frac{1}{4}(G(H-I\rho^2)-H^2\rho^2(r^2-1))\\
&=\frac{1}{4}H(G-H\rho^2(r^2+1))=\frac{1}{4} HK
\end{align*}
hold and that the derivatives $N_{i,\phi}$ of the components of the normal vector with respect to $\phi$ are given by the components $N_i$ according to
\begin{equation}
\label{eq:3.7}
N_{i,\phi}=N_{i-1} \quad \text{for $i$ even}, \qquad N_{i,\phi}=-N_{i+1} \quad \text{for $i$ odd},
\end{equation}
thus finishing the proof. We remark that, since $\nabla g_t=0$, one has that $g_t(N,\nabla_Y N)=0$ for all vector fields $Y \in \X(M^3_\Gamma$), and a computation indeed shows that the normal part of $\nabla_{\gd_s} N$ vanishes.
\end{proof}
\begin{proposition}
\label{prop:3.4}
\begin{displaymath}
\nabla_{\gd_{\rho}}N=0.
\end{displaymath}
\end{proposition}
\begin{proof}
One computes
\begin{align*}
\nabla_{\gd_{\rho}}N%&=\left(N_{i,\rho}+N_k dx_j\left (\gd_{\rho}\right)\Gamma^i_{jk} \right) \gd_{x_i}\\
&=\sum\limits_{i=1}^4 \left[N_{i,\rho}+\sum\limits_{k=1}^4 N_k\left ( \Gamma^i_{1k} \zeta_1+\Gamma_{2k}^i \zeta_2+\Gamma^i_{3k}\cos \phi +\Gamma^i_{4k} \sin \phi\right) \right] \gd_{x_i}\\
&=\sum\limits_{i=1}^4\Big[N_{i,\rho}+\Gamma^i_{11}(N_1\zeta_1-N_2\zeta_2)+\Gamma^i_{12}(N_2\zeta_1+N_1\zeta_2)\\
&+\Gamma^i_{13}(N_3\zeta_1-N_4\zeta_2+N_1\cos \phi - N_2\sin \phi)\\&+\Gamma^i_{14}(N_4\zeta_1+N_3\zeta_2+N_2\cos \phi + N_1 \sin \phi)\\
&+\Gamma^i_{33}(N_3 \cos \phi - N_4 \sin \phi) +\Gamma^i_{34}(N_4 \cos \phi+ N_3 \sin \phi)\Big] \gd_{x_i}.
\end{align*}
Once again we calculate the components of $\nabla_{\gd_{\rho}}N$ separately.
By using the symmetries of the $\Gamma^i_{jk}$ and \eqref{eq:3.5}, \eqref{eq:3.6} one obtains
\begin{eqnarray*}
\lefteqn{\Gamma^1_{11}(N_1\zeta_1-N_2\zeta_2)+\Gamma^1_{12}(N_2\zeta_1+N_1\zeta_2)}\\
&&=\frac{1}{2}\sqrt{\frac{K}{r^2+1}}A_1\rho [ \zeta_1(w_1\zeta_1+w_2\zeta_2)+\zeta_2(-w_2\zeta_1+w_1\zeta_2)]=\frac{1}{2}\sqrt{\frac{K}{r^2+1}}A_1 \rho r^2 w_1,
\end{eqnarray*}
\begin{eqnarray*}
\lefteqn{\Gamma^1_{13}\big(N_3\zeta_1-N_4\zeta_2+N_1\cos \phi-N_2 \sin \phi\big)+ \Gamma^1_{14}\big(N_4\zeta_1+N_3\zeta_2+N_2\cos \phi+N_1 \sin \phi\big)}\\
&&=-\frac{1}{2}\sqrt{\frac{K}{r^2+1}}B_1\rho\Big[\cos\phi\big(r(-w_3\zeta_1-w_4\zeta_2)+w_1 \cos \phi + w_2 \sin \phi \big )\\&&+\sin \phi \big ( r(w_4\zeta_1-w_3\zeta_2)-w_2 \cos \phi +w_1\sin \phi \big )\Big ]\\
&&=-\frac{1}{2}\sqrt{\frac{K}{r^2+1}}B_1\rho \Big[w_1+rw_4(\zeta_1\sin \phi-\zeta_2 \cos \phi)-rw_3(\zeta_2\sin \phi+\zeta_1 \cos \phi)\Big]\\
&&=-\frac{1}{2}\sqrt{\frac{K}{r^2+1}}B_1\rho [w_1+rw_4(-v)-rw_3u]=-\frac{1}{2}\sqrt{\frac{K}{r^2+1}}B_1\rho(1-r^2)w_1,\\
\lefteqn{\Gamma^1_{33}(N_3\cos \phi - N_4 \sin \phi) + \Gamma^1_{34}(N_4 \cos \phi + N_3 \sin \phi)}\\
&&=\frac{1}{4} \sqrt{\frac{K}{r^2+1}} r\rho ^3 C\bigg [\bigg(\frac{v}2 \sin \phi + \frac{u}2 \cos \phi\bigg) (-w_3 \cos \phi-w_4 \sin \phi)\\
&&+\bigg(\frac{u}2 \sin \phi- \frac{v}2 \cos \phi\bigg) (w_4 \cos \phi -w_3 \sin \phi)\bigg ]\\
&&=\frac{1}{8} \sqrt{\frac{K}{r^2+1}} r\rho ^3 C (-uw_3-vw_4)=-\frac{1}{8} \sqrt{\frac{K}{r^2+1}} r^2\rho ^3 C w_1,
\end{eqnarray*}
and, moreover,
\begin{displaymath}
dx_1(\nabla_{\gd_{\rho}}N)=N_{1,\rho}+\frac{1}{2} \sqrt{\frac{K}{r^2+1}}\rho \left((A_1+B_1)r^2-\left(B_1+\frac{1}4 \rho^2r^2C\right)\right)w_1.
\end{displaymath}
A similar calculation gives for the second component the expression
\begin{displaymath}
dx_2(\nabla_{\gd_{\rho}}N)=N_{2,\rho}+\frac{1}{2} \sqrt{\frac{K}{r^2+1}}\rho \left(-(A_1+B_1)r^2+\left(B_1+\frac{1}4 \rho^2r^2C\right)\right)w_2.
\end{displaymath}
One calculates further
\begin{eqnarray*}
\lefteqn{\Gamma^3_{11}(N_1\zeta_1-N_2\zeta_2)+\Gamma^3_{12}(N_2\zeta_1+N_1\zeta_2)}\\
&&=\frac{\rho^3 C}{4}\sqrt{\frac{K}{r^2+1}} \bigg[\bigg(-uv \sin \phi+ \frac{u^2-v^2}2 \cos \phi\bigg)(w_1\zeta_1+w_2\zeta_2)\\
&&+\bigg(uv \cos \phi +\frac{u^2-v^2}2 \sin \phi\bigg)(-w_2\zeta_1+w_1\zeta_2)\bigg]\\
&&=\frac{\rho^3 C}{4}\sqrt{\frac{K}{r^2+1}}\bigg[-uv \Big(w_1(\zeta_1 \sin \phi-\zeta_2 \cos \phi) +w_2(\zeta_2 \sin \phi+\zeta_1\cos \phi )\Big)\\
&&+\frac{u^2-v^2}2\Big( w_1(\zeta_1 \cos \phi +\zeta_2\sin \phi)+w_2(\zeta_2 \cos \phi-\zeta_1 \sin \phi) \Big )\bigg]\\
&&=\frac{\rho^3 C}{4}\sqrt{\frac{K}{r^2+1}} \bigg[-uv(-vw_1+uw_2)+\frac{u^2-v^2}2(uw_1+vw_2)\bigg]\\
&&=\frac{\rho^3 C}{4}\sqrt{\frac{K}{r^2+1}} \frac{u^2+v^2}2(uw_1-vw_2)=\frac{1}{8}\sqrt{\frac{K}{r^2+1}}\rho^3r^3Cw_3,
\end{eqnarray*}
as well as
\begin{eqnarray*}
\lefteqn{\Gamma^3_{13}(N_3\zeta_1-N_4\zeta_2+N_1\cos \phi-N_2 \sin \phi)+ \Gamma^3_{14}(N_4\zeta_1+N_3\zeta_2+N_2\cos \phi}\\&&+N_1 \sin \phi)=\frac{\rho B_2}{2}\sqrt{\frac{K}{r^2+1}} \Big [ -\zeta_1\big(r(-w_3\zeta_1-w_4\zeta_2)+w_1 \cos \phi + w_2 \sin \phi \big )\\
&&-\zeta_2 \big ( r(w_4\zeta_1-w_3\zeta_2)-w_2 \cos \phi +w_1\sin \phi \big )\Big ]\\
&&=\frac{\rho B_2}{2}\sqrt{\frac{K}{r^2+1}} \Big[w_1(-\zeta_1\cos \phi-\zeta_2 \sin \phi)+w_2(-\zeta_1\sin \phi + \zeta_2 \cos \phi)\\ &&+rw_3(\zeta_1^2+\zeta_2^2)\Big]\\
&&=\frac{\rho B_2}{2}\sqrt{\frac{K}{r^2+1}}(-uw_1+vw_2+r^3w_3)=\frac{1}{2}\sqrt{\frac{K}{r^2+1}}\rho r(r^2-1) B_2w_3,\\
\lefteqn{\Gamma^3_{33}(N_3\cos \phi - N_4 \sin \phi) + \Gamma^3_{34}(N_4 \cos \phi + N_3 \sin \phi)}\\
&&=\frac{1}{2}\sqrt{\frac{K}{r^2+1}}\rho r A_2 \Big[\cos \phi(-w_3 \cos \phi-w_4 \sin \phi)+\sin \phi(w_4 \cos \phi-w_3 \sin \phi)\Big]\\
&&=-\frac{1}{2}\sqrt{\frac{K}{r^2+1}}\rho r A_2 w_3,
\end{eqnarray*}
thus obtaining  for the third component of  $\nabla_{\gd_{\rho}}N$ that
\begin{displaymath}
dx_3(\nabla_{\gd_{\rho}}N)=N_{3,\rho}+\frac{1}{2} \sqrt{\frac{K}{r^2+1}}\rho r\left(-(A_2+B_2)+\left(B_2+\frac{1}4 \rho^2C\right)r^2\right)w_3.
\end{displaymath}
Finally, by an analogous calculation one finds that the fourth component reads
\begin{displaymath}
dx_4(\nabla_{\gd_{\rho}}N)=N_{4,\rho}+\frac{1}{2} \sqrt{\frac{K}{r^2+1}}\rho r\left(A_2+B_2-\left(B_2+\frac{1}4 \rho^2C\right)r^2\right)w_4.
\end{displaymath}
Since
\begin{gather*}
A_2+B_2=\frac{1}{2}H^2\rho^2r^2-\frac{1}4G(H-I\rho^2)+\frac{1}4G(H-I\rho^2)-\frac{1}4H^2\rho^2(r^2+1)\\
=-\frac{1}4H(G-H\rho^2(r^2+1))=-\frac{1}4 HK,\\
\frac{1}4\rho^2r^2C+B_1=\frac{1}4[G(H-I\rho^2r^2)+H^2\rho^2(r^2-1)]+\frac 14 \rho^2r^2(IG-2H^2)=\frac{1}4HK,
\end{gather*}
the desired statement follows by noting that $N_{i,\rho}= (\log K)_{,\rho} N_i/2$
and 
\begin{displaymath}
\frac{1}2 (\log K)_{,\rho}-\frac{1}4 HK \rho (r^2+1)=0.
\end{displaymath}
\end{proof}

It remains to compute the covariant derivative of $N$ with respect to $\gd_{\phi}$. 

\begin{proposition}
\label{prop:3.5}
\begin{displaymath}
\nabla_{\gd_{\phi}}N=\sum\limits^4_{i=1} \Xi N_{i,\phi} \gd_{x_i}=\Xi N_{,\phi},
\end{displaymath}
where $\Xi=1-\frac{1}4 HK \rho^2(r^2+1)$.
\end{proposition}
\begin{proof}
Again,
\begin{eqnarray*}
\lefteqn{\nabla_{\gd_\phi}N%=\left(N_{i,\phi}+N_kdx_j\left(\gd_{\phi}\right)\Gamma_{jk}^i\right)\gd_{x_i}}\\
=\sum\limits^{4}_{i=1} \Big [ N_{i,\phi}+\sum\limits_{k=1}^4 N_k\big[-\rho \zeta_2 \Gamma^i_{1k}+\rho \zeta_1 \Gamma^i_{2k}-\rho \sin \phi \Gamma^i_{3k}+\rho \cos \phi \Gamma^i_{4k} \big] \Big]\gd_{x_i}}\\
&&=\sum\limits^{4}_{i=1}\Big[N_{i,\phi}+\rho \Gamma^i_{11}(-\zeta_2N_1-\zeta_1N_2) + \rho \Gamma^i_{12}(-\zeta_2N_2+\zeta_1N_1)\\
&&+\rho \Gamma^i_{13}(-\zeta_2N_3-N_1\sin \phi  - \zeta_1 N_4-N_2 \cos \phi )
\\&&+\Gamma^i_{14} \rho(-\zeta_2N_4 +N_1 \cos \phi+\zeta_1N_3-N_2 \sin \phi)\\
&&+\rho \Gamma^i_{33}(-N_3\sin \phi-N_4 \cos \phi)+\Gamma^i_{34}\rho(-N_4 \sin \phi + N_3 \cos \phi)\Big]\gd_{x_i}.
\end{eqnarray*}
By the symmetries of the $\Gamma^i_{kj}$ one has for $i$ odd that 
\begin{displaymath}
dx_i(\nabla_{\gd_\phi} N)-N_{i,\phi}=-\rho \left(dx_{i+1}(\nabla_{\gd \rho} N)-N_{i+1,\rho} \right )
\end{displaymath}
and one obtains for the first and third component
\begin{align*}
dx_1(\nabla_{\gd_{\phi}}N)&=N_{1,\phi}+\frac{1}{2} \sqrt{\frac{K}{r^2+1}}\rho^2 \left((A_1+B_1)r^2-\left(B_1+\frac{1}4 \rho^2r^2C\right)\right)w_2,\\
dx_3(\nabla_{\gd_{\phi}}N)&=N_{3,\phi}+\frac{1}{2} \sqrt{\frac{K}{r^2+1}}\rho^2 r\left(-(A_2+B_2)+\left(B_2+\frac{1}4 \rho^2C\right)r^2\right)w_4.
\end{align*}
In an analogous way one has for $i$ even
\begin{displaymath}
dx_i(\nabla_{\gd_\phi} N)-N_{i,\phi}=\rho \left(dx_{i-1}(\nabla_{\gd \rho} N)-N_{i-1,\rho} \right ),
\end{displaymath}
the second and fourth component being given by
\begin{align*}
dx_2(\nabla_{\gd_{\phi}}N)&=N_{2,\phi}+\frac{1}{2} \sqrt{\frac{K}{r^2+1}}\rho^2 \left((A_1+B_1)r^2-\left(B_1+\frac{1}4 \rho^2r^2C\right)\right)w_1,\\
dx_4(\nabla_{\gd_{\phi}}N)&=N_{4,\phi}+\frac{1}{2} \sqrt{\frac{K}{r^2+1}}\rho^2 r\left(-(A_2+B_2)+\left(B_2+\frac{1}4 \rho^2C\right)r^2\right)w_3,
\end{align*}
and the assertion follows. Again, one verifies that $g_t(N,\nabla_{\gd_\phi} N)=0$.
\end{proof}

We are now able to compute the second fundamental form of the hypersurface $M^{3}_\Gamma$.

\begin{theorem}
\label{thm:3.1}
With respect to the coordinate frame $(\gd_ s, \gd_ \rho, \gd_ \phi)$ the components of the second fundamental form of the Riemannian $C^\infty$--manifolds $(M^3_\Gamma,h_t)$  are given by
\begin{equation*}
\II=\frac{\rho}2\sqrt{\frac{K}{r^2+1}} \left (\begin{array}{ccc} G(\dot u \ddot v - \dot v \ddot u)-2H\rho^2 (u\dot v -v\dot u) & 0 & K \\ 0 & 0 & 0 \\ K & 0 & 0 \end{array}\right ).
\end{equation*}
\end{theorem}

\begin{proof}
By the equations \eqref{eq:3.3}, \eqref{eq:3.4} and Proposition \ref{prop:3.3} one has
\begin{eqnarray*}
\lefteqn{\II_{11}=\sum\limits_{i,j=1}^4g_{ij} dx_i(\gd_s) dx_j(\nabla_{\gd_s}N)=\rho\Big[(N_{3,s}+\Phi(\dot r N_3 + r \dot \phi_\Gamma N_{3,\phi}))}\\
&&\cdot(-w_2G_4+w_1G_3)+(N_{4,s}+\Phi(\dot r N_4 + r \dot \phi_\Gamma N_{4,\phi}))(-w_2G_3+w_1G_4)\\
&&+\big ((N_{1,s}+\Phi(\dot r N_1 + r \dot \phi_\Gamma N_{1,\phi}))w_2+(N_{2,s}+\Phi(\dot r N_2 + r \dot \phi_\Gamma N_{2,\phi}))w_1 \big ) G_1\Big ]\\
&&=\rho \Big [H\rho^2 \big [\Phi(-r \dot r(N_3w_4+N_4w_3)+r^2 \dot \phi_\Gamma(N_4 w_4-N_3w_3))\\&&-r(N_{3,s}w_4+N_{4,s}w_3)\big ]+(G-H\rho^2r^2)\big [(N_{1,s}w_2+w_1 N_{2,s})\\&& +\Phi(\dot r (N_1w_2+N_2w_1)+r\dot \phi_\Gamma(-N_2w_2+N_1w_1))\big]\Big ],    
\end{eqnarray*}
where we made use of the relations \eqref{eq:3.7} as well as
\begin{gather*}
-w_2G_4+w_1G_3=(-w_2u-w_1v) \rho ^2 H=-rw_4\rho^2H,\\
-w_2G_3 -w_1 G_4=(w_2v-w_1u) \rho ^2 H= - r w_3\rho^2 H.
\end{gather*}
Since further $N_3w_4+N_4w_3=0$, $N_1w_2+N_2w_1=0$ and 
\begin{align*}
N_4w_4-N_3w_3&=\frac{1}2 \sqrt{\frac{K}{r^2+1}}r(w_3^2+w_4^2)=\frac{1}2 \sqrt{\frac{K}{r^2+1}}r\\
-N_2w_2+N_1w_1&=\frac{1}2 \sqrt{\frac{K}{r^2+1}}(w_1^2+w_2^2)=\frac{1}2 \sqrt{\frac{K}{r^2+1}},
\end{align*}
one obtains for the expression above that 
\begin{align*}
\II_{11}&=H\rho^3\left[\Phi r^3\dot \phi_\Gamma\frac{1}2 \sqrt{\frac{K}{r^2+1}}-r\left(-\left (\frac{1}2 \sqrt{\frac{K}{r^2+1}}r\right)_{,s}(w_3w_4-w_4w_3)\right.\right.\\
&\left.\left.+\frac{1}2 \sqrt{\frac{K}{r^2+1}}r(-w_{3,s}w_4+w_{4,s}w_3)\right ) \right ]+\rho(G-H\rho^2r^2)\left [ \Phi r \dot \phi_\Gamma\frac{1}2 \sqrt{\frac{K}{r^2+1}}\right.\\
&\left.+\left (\frac{1}2 \sqrt{\frac{K}{r^2+1}}r\right)_{,s}(w_1w_2-w_2w_1)+\frac{1}2 \sqrt{\frac{K}{r^2+1}}(w_{1,s}w_2-w_{2,s}w_1)\right]\\
&=\frac{\rho}2 \sqrt{\frac{K}{r^2+1}} \Big [Hr^2\rho^2\big(\Phi \dot \phi_\Gamma r-r\dot \phi \ddot r+\dot r(\dot r \dot \phi_\Gamma+r \ddot \phi_\Gamma) \big )\\
&+(\Phi r \dot \phi_\Gamma+\dot u \ddot v- \dot v \ddot u)(G-H\rho^2r^2)\Big ]\\
&=\frac{\rho}2 \sqrt{\frac{K}{r^2+1}}[G(\Phi r \dot \phi_\gamma+\dot u \ddot v -\dot v \ddot u)-H\rho^2r^2\dot \phi_\Gamma].
\end{align*}
Here we made use of the relation
$-w_{3,s}w_4+w_{4,s}w_3=-\dot r ( \dot r \dot \phi_\Gamma+ r \ddot \phi_\Gamma)+r \dot \phi_\Gamma \ddot r$ and
$w_{1,s}w_2-w_{2,s}w_1=\dot u \ddot v -\dot v \ddot u$
as well as
\begin{eqnarray*}
\dot u \ddot v - \dot v \ddot u=\dot \phi_\Gamma+\dot r ( \dot r \dot \phi_\Gamma+r \ddot \phi_\Gamma)-r\dot \phi_\Gamma\ddot r.
\end{eqnarray*}
Because of $G\Phi r \dot \phi_\Gamma=-\frac{1}4GHK\rho^2 r^2\dot \phi_\Gamma=-H\rho^2r^2\dot \phi_\Gamma$  one finally obtains
\begin{equation*}
\II_{11}=\frac{1}2\sqrt{\frac{K}{r^2+1}}\rho [ G(\dot u \ddot v - \dot v \ddot u)-2H\rho^2 (u\dot v -v\dot u)],
\end{equation*}
since  $r^2\dot \phi_\Gamma=u\dot v -v\dot u$. By proposition \ref{prop:3.4},
\begin{eqnarray*}
{\II_{22}=\sum \limits_{i,j=1}^4g_{ij} dx_i(\gd_\rho) dx_j(\nabla_{\gd_\rho}N)=0},
\end{eqnarray*}
and using  the equalities 
\begin{gather*}
G_4w_3-G_3w_4=(uw_3+vw_4) \rho ^2 H= r \rho ^2 w_1 H,\\
-G_3w_3-G_4w_4=(vw_3-uw_4) \rho ^2 H = - r \rho ^2 w_2 H
\end{gather*}
in addition to the above relations, one also sees that $\II_{33}$ vanishes, since by Proposition \ref{prop:3.5}  and  \eqref{eq:3.7}, 
\begin{eqnarray*}
\lefteqn{\II_{33}=\sum\limits_{i,j=1}^4g_{ij} dx_i(\gd_\phi) dx_j(\nabla_{\gd_\phi}N)}\\
&&=\Xi\Big [ -\rho \zeta_2(G_1N_{1,\phi}-G_4 N_{3,\phi}-G_3N_{4,\phi}) +\rho\zeta_1 ( G_1 N_{2,\phi} + G_3N_{3,\phi}-G_4N_{4,\phi})\\
&&-\rho \sin \phi(-G_4N_{1,\phi}+G_3N_{2,\phi}+G_2N_{3,\phi})\\&&+\rho\cos \phi(-G_3N_{1,\phi}-G_4N_{2,\phi}+G_2N_{4,\phi}) \Big ]\\
&&=\Xi\Big [ -\rho \zeta_2(-G_1N_2+G_4N_4-G_3N_3)+\rho \zeta_1 ( G_1 N_1 - G_3 N_4 -G_4N_3)\\&& -\rho \sin \phi ( G_4 N_2 +G_3N_1 - G_2 N_4)+ \rho \cos \phi( G_3N_2-G_4 N_1+ G_2 N_3)\Big]\\
&&=\frac{1}2\sqrt{\frac{K}{r^2+1}}\Xi \rho\Big [ G_1(\zeta_1w_1-\zeta_2w_2) +G_2 r ( - w_3\cos \phi +w_4\sin \phi )\\&&+\zeta_1 r ( G_4w_3-G_3w_4)+\zeta_2 r(-G_3 W_3-G_4w_4) \\&&+ \cos \phi(-G_4w_1 - G_3 w_2)+\sin \phi(-G_3w_1+G_4w_2)\Big]\\
&&=\frac{1}2\sqrt{\frac{K}{r^2+1}}\Xi\rho\Big [(G-H\rho^2r^2)r^2\dot \phi_\Gamma-(G-H\rho^2)r^2 \dot \phi_\Gamma\\
&&+H\rho^2 r (\zeta_1w_1r-\zeta_2w_2 r-w_3 \cos \phi +w_4 \sin \phi) \Big ]\\
&&=\frac{1}2\sqrt{\frac{K}{r^2+1}}\Xi\rho\Big [ (-H\rho^2r^2+H\rho^2)r^2 \dot \phi_\Gamma+H\rho^2 r ( r^3\dot \phi_\Gamma-r \dot \phi_\Gamma)\Big ]=0.
\end{eqnarray*}
Analogously,
\begin{eqnarray*}
{\II_{12}=g_t(\gd _s, \nabla_{\gd_\rho} N)=\sum\limits_{i,j=1}^4g_{ij} dx_i(\gd_s) dx_j(\nabla_{\gd_\rho}N)}=0,
\end{eqnarray*}
\begin{eqnarray*}
\lefteqn{\II_{13}=g_t(\gd _s, \nabla_{\gd_\phi} N)=\sum\limits_{i,j=1}^4g_{ij} dx_1(\gd_s) dx_j(\nabla_{\gd_\phi}N)}\\
&&=\Xi\rho [N_{3,\phi}(-G_4 w_2 +G_3 w_1) + N_{4,\phi}(-G_3w_2-G_4 w_1) +G_1 ( N_{1, \phi}w_2 + N_{2,\phi} w_1)]\\
&&=\Xi \rho \frac{1}2\sqrt{\frac{K}{r^2+1}}[-r w_4(-r \rho ^2w_4H)-rw_3(-r \rho ^2 w_3 H)+ G_1 ( w_1^2+w_2^2)]\\
&&=\Xi \rho \frac{1}2\sqrt{\frac{K}{r^2+1}}[r^2\rho^2 H+G-H\rho^2r^2]\\
&&=\rho \frac{1}2\sqrt{\frac{K}{r^2+1}}\left(1-\frac{1}4HK\rho^2(r^2+1) \right ) G=\frac{1}2\sqrt{\frac{K}{r^2+1}} K.
\end{eqnarray*}
Finally,
\begin{eqnarray*}
{\II_{32}=\sum\limits_{i,j=1}^4g_{ij} dx_i(\gd_\phi) dx_j(\nabla_{\gd_\rho}N)}=0,
\end{eqnarray*}
and the remaining components are determined by the symmetry of $\II$. 
\end{proof}

In order to compute the invariants of $\II$ we need to express the second fundamental form with respect to the orthonormal frame $Y_1, Y_2, Y_3$. In this case we denote its components by $\II^\ast_{ij}$. Let $\kappa_1,\kappa_2,\kappa_3$ be the eigenvalues of $\II$ in this base, regarded as a symmetric transformation on $TM^3_\Gamma$. The mean curvature, the first elementary symmetric function associated with $\II$, is then given by the sum $\HM=\kappa_1+\kappa_2+\kappa_3=\II^\ast_{11}+\II^\ast_{22}+\II^\ast_{33}$. Now, writing  $Y_i=a_{ij} \gd_{\eta_j}$ one has
\begin{equation}
\label{eq:3.8}
\II^\ast_{ij}=g_t(Y_i,\nabla _{Y_j} N)=\sum\limits_{k,l=1}^3a_{ik}a_{jl} g_t(\gd _{\eta_k}, \nabla _{\gd _{\eta_l}} N )=\sum\limits_{k,l=1}^3a_{ik}a_{jl} \II_{kl},
\end{equation}
i.e., $\II^\ast=A \cdot \II \cdot \,^T A$, where the coefficients of $A=(a_{ij})_{i,j}$ are determined by the equations \eqref{eq:2.1}. As a consequence of the previous theorem we obtain then the following result. 

\begin{corollary}
\label{cor:3.1}
The mean curvature of the hypersurfaces $(M^3_\Gamma,h_t)$ is given by
\begin{align*}
\HM=\sqrt{\frac 2 {\sqrt{\rho^4(r^2+1)^2+t^4}}}\cdot k_g,
\end{align*}
where $k_g:=\big ( (\dot u \ddot v -\dot v \ddot u)(r^2+1)-2(u\dot v-v\dot u)\big)/2$ denotes the geodesic curvature of the curve $\Gamma$ in $S^2$. In particular, $M^3_\Gamma$ is a minimal surface if and only if $k_g=0$, i.e., if $\Gamma$ is a great circle in $S^2$.
\end{corollary}
\begin{proof}
One easily sees that  $\II^\ast_{11}=\II^\ast_{22}=0$ as well as
\begin{align*}
\II^\ast_{33}&=D^2 \II_{11}+2DF \II_{13}=\frac{\rho}2 \sqrt{\frac K {r^2+1}}(r^2+1)K^2\Sigma^2\big [(r^2+1) G(\dot u \ddot v -\dot v \ddot u)\\&-2(u\dot v -v\dot u)(H(r^2+1) \rho^2 +K) \big ]\\&=\sqrt{\frac{1}{2\sqrt{\rho^4(r^2+1)^2+t^4}}}\big[(\dot u \ddot v - \dot v \ddot u)(r^2+1)-2(u\dot v - v\dot u) \big].
\end{align*}
So $\HM=\II^\ast_{33}$. We further remark that $\II^\ast_{12}=\II^\ast_{13}=0$ and
\begin{align*}
\II^\ast_{23}=\frac{\rho}2 \sqrt{\frac K {r^2+1}}\frac {DK}{\sqrt{h_{33}}} =\frac{K}{4\rho}\sqrt{\frac K {r^2+1}}=\frac{(r^2+1) \rho^2}{\sqrt{2}(\rho^4(r^2+1)^2+t^4)^{3/4}}.
\end{align*}
Let us now compute the geodesic curvature of $\Gamma$ regarded as a curve in  $S^2\simeq \C\cup \mklm{\infty}$ using the stereographic projection. With respect to the coordinates $u ,v$ the induced metric on $S^2$ reads 
\begin{equation}
\label{eq:3.9}
g_{S^2}=\frac 1 {(r^2+1)^2} \left   ( \begin{array}{cc} 1 & 0 \\ 0 & 1 \end{array} \right).
\end{equation}
Writing $v_1,v_2$ for $u,v$ the geodesic curvature of $\Gamma$ is then given by (see e.g. \cite{kreysig})
\begin{equation*}
k_g=\sqrt{\det g_{S^2}} \left | \begin{array}{cc} \dot v_1 & \dot v_2 \\ \ddot v_1 + \sum\limits_{i,j=1}^2\Gamma^1_{ij} \dot v_i \dot v_j & \ddot v_2 +\sum\limits_{i,j=1}^2 \Gamma^2_{ij} \dot v_i \dot v_j\end{array}  \right | /(g_{ij} \dot v_i \dot v_j)^{3/2},
\end{equation*}
the Christoffel symbols $\Gamma^i_{jk}$ being obtained from the   formulas \eqref{eq:3.0}, where  now the $g_{ij}$ denote the components of $g_{S^2}$ and the $x_i$ should be replaced by the corresponding $v_i$. Note that $\dot v_1^2+ \dot v_2^2=1$. We put $L:=(r^2+1)^2$, $M:=-4/(r^2+1)^3$ and obtain
\begin{align*}
\Gamma^1_{11}&=\Gamma^2_{12}=\Gamma^2_{21}=\frac {LM}{2} v_1, & \Gamma^2_{11}&=-\frac {LM}{2} v_2,\\
\Gamma^1_{12}&=\Gamma^1_{21}=\Gamma^2_{22}=\frac {LM}{2} v_2, & \Gamma^1_{22}&=-\frac {LM}{2} v_1;
\end{align*}
a direct calculation then yields
\begin{displaymath}
\left | \begin{array}{cc} \dot v_1 & \dot v_2 \\ \ddot v_1 + \sum\limits_{i,j=1}^2\Gamma^1_{ij} \dot v_i \dot v_j & \ddot v_2 + \sum\limits_{i,j=1}^2\Gamma^2_{ij} \dot v_i \dot v_j\end{array}  \right |=\dot v_1 \ddot v_2-\dot v_2 \ddot v_1 + \frac M {2L}(v_1 \dot v_2-v_2 \dot v_1),
\end{displaymath}
and by noting that $LM/2=-2/(r^2+1)$ and  $g_{ij} \dot v_i \dot v_j=\sqrt{\det g_{S^2}}=1/L$ one finally has (up to a sign)
\begin{equation*}
k_g=\frac 1{2} \Big ( (\dot u \ddot v -\dot v \ddot u)(r^2+1)-2(u\dot v-v\dot u)\Big),
\end{equation*}
and thus the assertion.
\end{proof}

The third elementary symmetric function associated with $\II$ is the Gauss curvature; it is given by $\kappa_1 \cdot \kappa_2 \cdot \kappa_3=\det \II^\ast$ and equal to zero; the second one is the so--called second order homogeneous curvature. 

\begin{corollary}
The second order homogeneous curvature of the hypersurfaces $(M^3_\Gamma,h_t)$ is 
\begin{equation*}
\kappa_1 \kappa_2+\kappa_1 \kappa_3+\kappa_2 \kappa_3=-\frac{\rho^4(r^2+1)^2}{8(\rho^4(r^2+1)^2+t^4)^{3/2}}=S/8.
\end{equation*}
\end{corollary}
\begin{proof}
As computed in the proof of the previous corollary, the components of $\II$ with respect to the orthonormal frame \eqref{eq:2.1} are given by
\begin{displaymath}
\II^\ast =\left ( \begin{array}{ccc} 0 & 0 & 0 \\ 0 & 0 & \frac{K}{4\rho}\sqrt{\frac K {r^2+1}} \\ 0 &\frac{K}{4\rho}\sqrt{\frac K {r^2+1}}  & \HM
\end{array} \right ).
\end{displaymath}
The roots of the characteristic polynomial$$\det (\II^\ast -\kappa \1)=-\kappa\big(-\kappa(\HM-\kappa)-K^3/16\rho^2(r^2+1)\big )$$ are then $\kappa_1=0$, $\kappa_{2,3}=(-\HM\pm \sqrt{\HM^2+K^3/16\rho^2(r^2+1)})/2$.
\end{proof}

Since ${u_1}_{|M^3_\Gamma}=\rho^2(r^2+1)^2$, the three elementary symmetric functions associated with the second fundamental form, i.e., essentially its  trace $\HM$ and the scalar curvature $S$, are manifestly invariant under the action of the isometry group $\U(2)$. 
The fact that the mean curvature of the hypersurfaces $M^3_\Gamma$ is given in terms of the geodesic curvature of $\Gamma$ in $S^2$ appears to be natural, since the geometry of the vector bundle $T^\ast\P^1(\C)$ is determined by the elliptic geometry of $\P^1(\C)\simeq S^2$. Note that  $\dot u \ddot v -\dot v \ddot u$ is the geodesic curvature of $\Gamma$ as a curve in $\C$ with respect to the Euclidean metric.

As an immediate consequence one obtains the following statement.

\begin{corollary}
\label{cor:3.3}
Let $\Gamma$ be a curve in $S^2$ of bounded geodesic curvature. Then the functionals 
\begin{equation*}
\int \HM^\alpha  \,d M^3_\Gamma \qquad \text{and} \qquad \int (\kappa_1 \kappa_2+\kappa_1 \kappa_3+\kappa_2 \kappa_3)^\beta  \,d M^3_\Gamma
\end{equation*}
stay bounded for $\alpha>3$ and $\beta>3/2$, respectively.
\end{corollary}
Consequently, the Willmore functional $\int \HM^3\,d M^3_\Gamma$ remains unbounded, the hypersurfaces $M^3_\Gamma$ thus being not accessible to integral geometry.

\section{On the geodesic flow of the hypersurfaces $M^3_\Gamma$}
\label{sec:8}

In this section we will study the structure of the geodesic flow  of the hypersurfaces $(M^3_\Gamma,h_t)$ and compute the exponential growth $\mu_\infty(M^3_\Gamma)$ explicitly,  at least in the case where $\Gamma$ is  a generalized circle in $\C$ that arises by a M\"{o}bius transform from a circle in $\C$ with center at the origin.  In general, \emph{the exponential growth} of an open, complete Riemannian manifold $(M^n,g)$ is defined as
\begin{displaymath}
\mu_\infty:= \limsup_{R \to \infty} \frac 1 R \log \vol(B_R(q_0)),
\end{displaymath}
where $q_0$ is a point in $M^n$ and  $\vol(B_R(q_0))$ denotes the volume of the ball of radius $R$ with center at $q_0$. If $\mu_\infty=0$, one says that $M^n$ has \emph{subexponential growth}. In case $M^n$ has finite volume, this quantity is not interesting, since then one always has  $\mu_\infty=0$, but for $\vol(M^n)=\infty$ the exponential growth is directly related to the infimum of the essential spectrum of the Laplace operator on $M^n$. We will return to this point in section \ref{sec:7}. There we will be able to calculate the exponential growth of $(M^3_\Gamma)$ for arbitrary closed curves.  

Let $\gamma(\tau)=\Psi(s(\tau),\rho(\tau),\phi(\tau))$ be a smooth curve in $M^3_\Gamma$ and $X(\tau)=\sum_{j=1}^3 X^j(\tau)$ $\gd/\gd_{\eta_j}$ a vector field along $\gamma(\tau)$. Its covariant derivative with respect to $\gamma$ is given by the formula
\begin{equation*}
\frac{\nabla X}{dt}=\sum\limits_{k=1}^3\Big(\frac d{dt} X^k(\tau)+\sum\limits_{i,j=1}^3 \Gamma^i_{jk} X^i(\tau)\dot \gamma^j(\tau) \Big) \gd_{\eta_k},
\end{equation*}
where $\Gamma^i_{jk}=h_t(\nabla_{\gd_{\eta_i}}\gd_{\eta_j},\gd_{\eta_k})$ are the components of the Levi--Civita connection of $M^3_\Gamma$ with respect to the coordinate frame $\mklm{\gd_s,\gd_\rho,\gd_\phi}$. For a geodesic  it holds that $\nabla \dot \gamma(\tau)/dt=0$  and one obtains the system of differential equations
\begin{equation}
\label{eq:8.0}
\ddot \gamma^k(\tau)+\sum\limits_{i,j} \Gamma^k_{ij}(\gamma(\tau)) \dot \gamma^i(\tau) \dot \gamma^j(\tau)=0, \qquad k=1,2,3.
\end{equation}
However, it turns out to be more convenient to determine the geodesic lines of the hypersurfaces $M^3_\Gamma$ by considering the first integrals of the geodesic flow.
Let us consider therefore the geodesic system $(TM^3_\Gamma,\E)$ of $M^3_\Gamma$, where the Lagrangian $\E$ is given by the metric,
$$\E:TM^3_\Gamma \rightarrow \R, \quad X \mapsto \frac 1 2 h_t (X,X).$$
The function $\E$ is a first integral of the geodesic flow, i.\ e.\ with respect to the coordinate reper $\mklm{\gd_s,\gd_\rho,\gd_\phi}$ one has that 
\begin{align*}
\label{eq:7.1}
\begin{split}
\E(\dot\gamma(\tau))&=\frac 1 2 \big [h_{11} \dot s^2(\tau) +2 (h_{12} \dot s(\tau) \dot \rho(\tau) +h_{13} \dot s(\tau) \dot \phi(\tau))+h_{22}\dot \rho^2(\tau) +h_{33} \dot \phi^2(\tau)\big ]\equiv \E
\end{split}
\end{align*}
is constant for any geodesic line. Let now $p=\Psi(s,\rho,\phi)$ and $z \in e^{i\tau} \in S^1$.  Since the coefficients of the metric $h_t$ do not depend on the angle variable $\phi$, the map
\begin{equation*}
\kappa_z:M^3_\Gamma\cap M^4 \longrightarrow M^3_\Gamma\cap M^4, \kappa_z(p)=\Psi(s,\rho,(\phi+\tau)\,\text{mod}\, 2\pi),
\end{equation*}
represents a one--parameter family of isometries.
Consequently, using Noether's theorem, the function
\begin{equation*}
\M_1:T(M^3_\Gamma\cap M^4) \rightarrow \R,\quad \M_1(X):=h_t\Big(\frac d{d\tau} \kappa_\tau(p)_{|\tau=0}\Big)=h_t(\gd_\phi,X)_{|p},
\end{equation*}
is a second first integral of the geodesic flow and a computation yields the formula
\begin{equation}
\label{eq:8.1}
\M_1(\dot \gamma(\tau))=h_{13} \dot s(\tau)+h_{33} \dot \phi(\tau)\equiv \M_1.
\end{equation}
For $\dot s=0$ and $\dot \phi=0$ it can be  seen immediately from the equations \eqref{eq:8.0} for a geodesic or the relation $\E(\dot \gamma(\tau))\equiv \E$ that, for 
\begin{equation}
\label{eq:8.1a}
\dot \rho^2=\frac {2\E}{(r^2+1)K}= \E \frac{\sqrt{\rho^4(r^2+1)^2+t^4}}{\rho^2(r^2+1)^2}, \qquad s,\, \phi \quad \text{constant},
\end{equation}
the curve $\gamma(\tau)=\Psi(s,\rho(\tau),\phi)$ must be a geodesic in $M^3_\Gamma$. 

We will assume from now on that $r\equiv r_0$ is constant and, in this case, determine the distance of a point $p=\Psi(s,\rho,\phi)$ to the set $\Gamma\equiv \mklm{[0,0]}\subset M^3_\Gamma$. Since $r(s)\equiv r_0$ and $\dot \phi_\Gamma(s)=\epsilon /r_0$ are constant, $\epsilon=\pm 1$, the coefficients $h_{ij}$ do also not depend on $s$ so that
\begin{equation*}
\mu_z:M^3_\Gamma\cap M^4 \longrightarrow M^3_\Gamma\cap M^4, \mu_z(p)=\Psi((s+\tau r) \, \text{mod} \, 2\pi r,\rho,\phi),
\end{equation*}
is an additional one--parameter group of isometries and Noether's Theorem gives a third first integral,
\begin{equation*}
\M_2:T(M^3_\Gamma\cap M^4) \rightarrow \R,\quad \M_2(X):=h_t\Big(\frac d{d\tau} \mu_\tau(p)_{|\tau=0}\Big)=h_t(\gd_s,X)_{|p},
\end{equation*}
i.\ e.,
\begin{equation}
\label{eq:8.2}
\M_2(\dot \gamma(\tau))=h_{11} \dot s(\tau)+h_{13} \dot \phi(\tau)\equiv \M_2
\end{equation}
is constant for any geodesic line as well. From the equations  \eqref{eq:8.1}, \eqref{eq:8.2} one obtains
\begin{equation}
\label{eq:8.3}
2\E=K(r^2+1)\dot \rho^2+\M_1 \dot \phi + \M_2 \dot s
\end{equation}
and
\begin{align*}
\dot s&=\epsilon \frac{\M_1-(r^2+1)K\rho^2\dot \phi}{r\rho^2 K}, & \dot \phi&=\epsilon\frac{\M_2-(K+H\rho^2)\rho^2\dot s}{r\rho^2K}. 
\end{align*}
Solving the latter two equations with respect to $\dot \phi$ and $\dot s$ yields
\begin{align*}
\dot \phi&=\left (\epsilon \M_2-\frac{(K+H\rho^2)\rho^2 \M_1}{rK\rho^2}\right ) (r\rho^2 K-(K+H\rho^2) \rho^2 (r^2+1)/r)^{-1}\\
&=\frac{(4+GH\rho^2) \M_1-4\epsilon\M_2 r}{4 \rho^2 G},
\end{align*}
as well as
\begin{align*}
\dot s&=\frac \epsilon{4 KG\rho^2r} [4 G\M_1-(r^2+1) (4 (\M_1-\epsilon\M_2r)K+ 4 \rho^2 H\M_1 )]\\
&=\frac \epsilon {16 \rho^2 r}[4(G-K(r^2+1)-H\rho^2(r^2+1)) \M_1 +4\epsilon(r^2+1)r K \M_2]\\
&=\frac{-\epsilon r\M_1+(r^2+1) \M_2}{4\rho^2}K;
\end{align*}
thus the functions $s=s(\tau)$ and $\phi=\phi(\tau)$ are determined by the function $\rho=\rho(\tau)$. Equation \eqref{eq:8.3} now reads
\begin{align*}
2\E&=\frac 1 G \left [ 4(r^2+1) \dot \rho^2+ \frac 1 {4 \rho^2} \Big ( 4 (\M_1^2-2 \epsilon \M_1\M_2 r+(r^2+1) \M_2^2) +GH \rho^2 \M_1^2 \Big ) \right ].
\end{align*}
Note that $\M_1^2-2 \epsilon \M_1\M_2 r+(r^2+1) \M_2^2=(\M_1-\epsilon \M_2 r)^2 +\M_2^2$ is non--negative. By inserting the expressions for $G$ and $H$ into the previous equation one finally obtains  the following ordinary differential equation for $\rho=\rho(\tau)$:
\begin{align}
\label{eq:8.4}
\dot \rho^2&=\frac{\sqrt{\rho^4(r^2+1)^2+t^4}}{\rho^2(r^2+1)^2} \E - \frac 1 {4\rho^2(r^2+1)} \left ( \frac{t^4 \M_1^2}{\rho^4(r^2+1)^3}+(\M_1-\epsilon \M_2 r)^2 +\M_2^2 \right ).
\end{align}
Thus, for $r(s)=r_0$, all geodesics $\gamma(\tau)=\Psi(s(\tau), \rho(\tau),\phi(\tau))$ in $M^3_{\Gamma}\cap M^4$ are parametrized by the three parameters $\E,\M_1,\M_2$. 
We are now able to compute the distance of a point $p \in M^3_{\Gamma(r= r_0)}$ to the curve $\Gamma$.
\begin{proposition}
\label{prop:8.1}
Let $\Gamma$ be a circle in $\C$ of radius $r(s)= r_0$. The distance of a point $p_0=\Psi(s_0,\rho_0,\phi_0)\in M^3_{\Gamma}\cap M^4$ to the curve $\Gamma\subset M^3_{\Gamma}$  is given by
\begin{equation}
\label{eq:8.4a}
\dist(p_0,\Gamma)=\frac 1 {t\sqrt{2}} \rho_0^2(r_0^2+1) \F(1/ 2, 1/4,3/ 2, -\frac{\rho_0^4(r_0^2+1)^2}{t^4}), 
\end{equation}
where $\F$ denotes the hypergeometric function, which  is defined for $z\in \C, \, |z|<1,$ by the series
\begin{equation*}
\F(\alpha,\beta,\gamma,z)=1+\frac{\alpha \, \beta}{\gamma \cdot 1} z + \frac{\alpha(\alpha+1) \beta(\beta+1)}{\gamma(\gamma+1) \cdot 1 \cdot 2 }z^2+\dots ,
\end{equation*}
the parameters $\alpha,\beta,\gamma$ being arbitrary complex numbers, $\gamma \not=0,\pm 1,\pm 2,\dots$.
\end{proposition}
\begin{proof}
Let $\gamma(\E,\M_1,\M_2):(0,\tau_0] \rightarrow M^3_{\Gamma}$ be a geodesic of positive energy $\E$ from the curve $\Gamma$  to the point $p_0$ with coordinates $s_0,\rho_0,\phi_0$. For $\M_1,\,\M_2=0$ the geodesic $\gamma(\E,\M_1,\M_2)$ is precisely the geodesic line \eqref{eq:8.1a} already described. If $\M_1$  were not equal zero, at least $\dot \phi$ would be  different from zero almost everywhere; then equation \eqref{eq:8.4} would imply that there exists a critical value $\rho_{\mathrm{crit}}>0$ for which
\begin{equation}
\label{eq:8.5}
\E\sqrt{\rho^4(r^2+1)^2+t^4}=\frac{r^2+1}4\left (\frac{t^4\M_1^2}{\rho^4(r^2+1)^3} +(\M_1-\epsilon \M_2 r )^2 + \M_2^2\right).
\end{equation}
For smaller values of $\rho$ the right--hand side of \eqref{eq:8.4} would become negative, implying that $\rho(\tau) \geq \rho_{\mathrm{crit}} >0$ must hold for all $\tau \in (0,\tau_0]$. This means that for $\M_1\not=0$ the geodesic $\gamma(\E,\M_1,\M_2)$ can never reach the curve $\Gamma$. Assume therefore $\M_1=0$, $\M_2$ being arbitrary. By \eqref{eq:8.4} we have
\begin{displaymath}
\dot \rho^2=\frac{4\sqrt{\rho^4(r^2+1)^2+t^4}\E-(r^2+1)^2 \M_2^2}{4\rho^2(r^2+1)^2}.
\end{displaymath}
In case that $4t^2 \E-(r^2+1)^2 \M_2^2<0$, this expression becomes negative for small $\rho$ so that $\gamma(\E,M_1=0, \,4t^2\E/(r^2+1)^2 < M^2_2)$ can never reach the set $\Gamma$. However, for $4t^2 \E-(r^2+1)^2 \M_2^2\geq 0$ we have  that $\dot \rho^2$ is non--negative for all $\tau$, as well as
\begin{displaymath}
\dot \phi=-\epsilon r \frac {r^2+1}{2\sqrt{\rho^4(r^2+1)^2+t^4}} \M_2, \qquad \dot s= \frac{(r^2+1)^2}{2\sqrt{\rho^4(r^2+1)^2+t^4}} \M_2,
\end{displaymath}
so there are infinitely many geodesic lines $\gamma(\E,M_1=0, \,4t^2\E/(r^2+1)^2 \leq M^2_2 )$ reaching the set $\Gamma$ in $M^3_\Gamma$ in a spiral motion. In this case, equation \eqref{eq:8.4} implies for $u_1(\tau)=\rho^2(\tau)(r^2+1)$ the relation
\begin{equation}
\label{eq:8.5b}
u_{1,\tau}=2 \rho \dot \rho(r^2+1)=\sqrt{4\E \sqrt{u^2_{1}+t^4}-(r^2+1)^2M_2^2}>0,
\end{equation}
i.e.,\ $u_{1,\tau}$ as well as $u_1$ are strictly monotone increasing as functions in $\tau$ and the point $p_0$ is reached earliest, that is, for smallest $\tau_0$, in case that $\M_2$ is also zero. Since the length of a geodesic is given by
\begin{align*}
L_{\gamma(\E,\M_1,\M_2)}=\int \limits_0^{\tau_0} \big|\dot \gamma(\E,\M_1,\M_2)\big| d\tau
=\int \limits_0^{\tau_0} \sqrt{h_t(\dot \gamma,\dot \gamma)} d\tau=\sqrt{2\E} \tau_0,
\end{align*}
the distance of the point $p_0$ to the set $\Gamma\subset M^3_\Gamma$ must be given by the length of the geodesic $\gamma(\E,\M_1=\M_2=0)$.

The integral $\int 1/\sqrt[4]{ax^2+t^4} dx$ cannot be represented by elementary functions and one has
\begin{equation}
\label{eq:8.5a}
\int \frac{dx}{\sqrt[4]{ax^2+t^4}}=\frac{x}{t} \F(1/2,1/4,3/2,-ax^2/t^4),
\end{equation}
where $\F(\alpha,\beta,\gamma,z)$ is the hypergeometric function introduced above. For $\Re(\alpha+\beta-\gamma)<0$ the defining series  converges even in $|z|\leq 1$; the hypergeometric function has an analytic continuation for $|z|>1$ and under the assumption that $\Re\gamma>\Re \beta >0$ it can be written for all $z$ as the integral 
\begin{displaymath}
\F(\alpha,\beta,\gamma,z)= \frac{\Gamma(\gamma)}{\Gamma(\beta)\Gamma(\gamma-\beta)} \int\limits_0^1 (-\zeta)^{\beta-1} (1-\zeta)^{\gamma-\beta-1}(1-\zeta z)^{-\alpha} d\zeta,
\end{displaymath}
where  $\Gamma$ denotes the Gamma function and $|\mathrm{arg}(-z)|<\pi$ is assumed in order to make the integrand uniquely defined. If  $z$ is real, differentiation under the integral with respect to $z$ gives the stated equality \eqref{eq:8.5a} if one takes the relation $\F(m,\beta,\beta,z)=(1-z)^{-m}$, $m \in \R$, $\beta$ arbitrary, into account additionally.
For $\M_1=M_2=0$ we finally deduce from \eqref{eq:8.5b} 
\begin{align*}
\tau_0&=\frac 1 {2\sqrt{\E}}\int\limits_0^{\tau_0} \frac{u_{1,\tau}}{\sqrt[4]{u_1^2+t^4}} d\tau=\frac 1 {2\sqrt{\E}}\int\limits_0^{u_1(\tau_0)} \frac{du_1}{\sqrt[4]{u_1^2+t^4}}\\
&=\frac{1}{2t\sqrt{\E}} u_1(\tau_0) \F(1/2,1/4,3/2,-u_1^2(\tau_0)/t^4),
\end{align*}
and thus 
\begin{displaymath}
\dist(p_0,\Gamma)=\frac 1 {t\sqrt{2}} \rho_0^2(r^2+1) \F(1/ 2, 1/4,3/ 2, -\frac{\rho_0^4(r^2+1)^2}{t^4}), 
\end{displaymath}
finishing the proof.
\end{proof}
We are now in a position to compute the exponential growth of the hypersurface $M^3_{\Gamma}$ in case that $\Gamma=\gd B(0,r_0)$ is a circle in $\C$. Note that we can estimate the volume of the ball with radius $R$ around a point $q_0 \in \Gamma \subset M^3_\Gamma$ by the volume of the union over all $R$--balls around points of $\Gamma$, thus obtaining 
\begin{gather*}
\vol(B_R(q_0))\leq \vol\Big(\bigcup\limits_{q \in \Gamma} B_R(q)\Big)
=\int\limits_0^{2\pi}\int\limits_0^{\rho_R}\int\limits_0^{2\pi r} \sqrt{\det h_t}\, ds \wedge d\rho \wedge d\phi\\
=2\pi\sqrt{8}\int\limits_0^{\rho_R}\int\limits_0^{2\pi r} \frac{\rho^3(r^2+1)}{\sqrt[4]{\rho^4(r^2+1)^2+t^4}} ds \wedge d\rho=\frac{4\pi^2 r \sqrt {8}}{3(r^2+1)} \big[ (\rho_R^4(r^2+1)^2+t^4)^{3/4}-t^3],
\end{gather*}
since by our previous considerations $$\bigcup\limits_{q \in \Gamma} B_R(q)=\mklm{p=\Psi(s,\rho_R,\phi) \in M^3_{\Gamma}:
 s \in [0,2\pi r),\, \phi \in [0,2\pi)},$$ where $\rho_R$ is  given by the expression \eqref{eq:8.4a} for $R=\dist(p,\Gamma)$. 
The analytic continuation of $\F(\alpha,\beta,\gamma,z)$ for $|z|>1$ is given by the formula
\begin{align*}
\F(\alpha,\beta,\gamma,z)&=\frac{\Gamma(\gamma)\Gamma(\beta-\alpha)}{\Gamma(\beta)\Gamma(\gamma-\alpha)} (-z)^{-\alpha} \F(\alpha,\alpha+1-\gamma,\alpha+1-\beta),1/z)+\\
&+\frac{\Gamma(\gamma)\Gamma(\alpha-\beta)}{\Gamma(\alpha)\Gamma(\gamma-\beta)} (-z)^{-\beta} \F(\beta,\beta+1-\gamma,\beta+1-\alpha,1/z),
\end{align*}
so that for $\rho_0$ being big enough the distance of $p_0=\Psi(s_0,\rho_0,\phi_0)$ to the set $\Gamma$ is given by
\begin{align*}
&\dist(p_0,\Gamma)=\frac 1 {t\sqrt{2}} \rho_0^2(r^2+1) \left [\frac{\Gamma(3/2)\Gamma(-1/4)}{\Gamma(1/4)\Gamma(1)} \frac{t^2}{\rho_0^2(r^2+1)} \F\left(\frac 1 2,0,\frac 5 4,\frac{-t^4}{\rho_0^4(r^2+1)^2}\right)\right.\\&+\left.\frac{\Gamma(3/2)\Gamma(1/4)}{\Gamma(1/2)\Gamma(5/4)} \frac{t}{\rho_0\sqrt{r^2+1}}\F\left(\frac 1 4,-\frac{1}{4},\frac3 4,\frac{-t^4}{\rho_0^4(r^2+1)^2}\right)\right ]\\
&=\frac{\Gamma(3/2)\Gamma(-1/4)}{\Gamma(1/4)\Gamma(1)}\frac{t}{\sqrt{2}}+ \frac{\Gamma(3/2)\Gamma(1/4)}{\Gamma(1/2)\Gamma(5/4)} \sqrt{\frac{r^2+1}{2}}\rho_0 \left (1+\frac 1 {12} \frac {t^4}{\rho_0^4(r^2+1)^2}+\dots \right ),
\end{align*}
implying that $\dist(p_0,\Gamma)$ is proportional to $\rho_0\sqrt{r^2+1}$ for $1\ll\rho_0$.
We obtain for $q_0\in \Gamma$ that 
\begin{gather*}
\lim_{R \to \infty} \frac 1 R \log \vol (B_R(q_0))\leq\\
\lim_{\rho \to \infty} \left [ \log \frac{4\sqrt{8} \pi^2r}{3(r^2+1)}+\log \big ((\rho^4(r^2+1)^2+t^4 ) ^{3/4}-t^3\big ) \right ]\cdot \left [\frac{\Gamma(3/2)\Gamma(-1/4)}{\Gamma(1/4)\Gamma(1)}\frac{t}{\sqrt{2}}+\right.\\
+\left.\frac{\Gamma(3/2)\Gamma(1/4)}{\Gamma(1/2)\Gamma(5/4)} \sqrt{\frac{r^2+1}{2}}\rho \left (1+\frac 1 {12} \frac {t^4}{\rho^4(r^2+1)^2}+\dots \right )\right]^{-1}\\
=\lim_{\rho \to \infty}\frac{3 \rho^3(r^2+1)^2}{\rho^4(r^2+1)^2+t^4-t^3\sqrt[4]{\rho^4(r^2+1)^2+t^4}}\left [\frac{\Gamma(3/2)\Gamma(1/4)}{\Gamma(1/2)\Gamma(5/4)}\sqrt{\frac{r^2+1}{2}}\right.\\\left. \cdot \left(1-\frac 1 4 \frac{t^4}{\rho^4(r^2+1)^2}+ \dots \right ) \right ]^{-1}=0,
\end{gather*}
the corresponding limes superior therefore being zero, too. By isometry arguments we thus obtain the following proposition.

\begin{proposition}
Let $\Gamma=\gd B(0,r_0)$ be a circle in $\C$ with center at the origin and radius $r_0$. Then $\mu_\infty(M^3_{A\Gamma})=0$ for all $A \in \U(2)$.
\endproof
\end{proposition}

\medskip

We want to finish this section with some remarks concerning closed geodesics in $M^3_\Gamma$, where we assume again that $\Gamma$ is  a circle in $\C$ of radius $r(s)= r_0$. In this case $M^3_{\Gamma}$ is foliated by the two--dimensional tori $T^2_{\rho_0,r_0}$, $\rho_0>0$ being constant. Let $\gamma(\tau)=\Psi(s(\tau),\rho_0,\phi(\tau)):[0,L_\gamma] \rightarrow T^2_{\rho_0,r_0}\subset M^3_{\Gamma}$ be a geodesic line parametrized by arc length. Since $\dot \rho=0$, $\ddot s$ and $\ddot \phi$ are also zero. Relation \eqref{eq:8.3} then reads $\M_2 \dot s= 2\E - \M_1 \dot \phi$ and equation \eqref{eq:8.5} must hold, representing  a condition on $\E=\E(\rho_0,r_0)$ for  given values of  $\M_1,\M_2$. Now, by the $S^1\times S^1$--symmetry of $T^2_{\rho_0,r_0}$,
\begin{displaymath}
\Psi(s_0,\rho_0,\phi_0)=\Psi(s_0+2\pi r_0 \cdot n, \,\rho_0,\,\phi_0+2\pi \cdot m), \qquad n,m \in \Z.
\end{displaymath}
Writing $s(\tau)=s_0+\dot s \tau$, $\phi(\tau)=\phi_0+\dot \phi \tau$ we see that $\gamma(\tau)$ is a closed geodesic if and only if $\dot s L_\gamma=2\pi r_0 \cdot n$ and $ \dot \phi L_\gamma=2\pi\cdot m$ are satisfied, i.\ e., if
\begin{displaymath}
\frac{\dot s}{\dot \phi}=r_0 \cdot \frac n m \in r_0 \cdot \Q^\ast, \qquad n,m \not=0;
\end{displaymath}
of course, if  $n$ or $m$ are zero, $\gamma(\tau)$ is also a closed geodesic. Inserting the expressions for $\dot s$ and $\dot \phi$ computed above we obtain for the previous condition
\begin{equation*}
\frac{-\epsilon r_0\M_1+(r^2_0+1)\M_2}{4(\M_1-\epsilon r_0 \M_2)+GH\rho_0^2\M_1} \in r_0 \cdot \Q^\ast.
\end{equation*}
Note that $GH=4t^4/\rho^6(r^2+1)^3$. Taking all together we find as solutions for $\M_1$ and $\M_2$
\begin{align*}
\M_1&=\frac{m(r_0^2+1)/4+\epsilon n r_0^2}{\rho_0^4(r_0^2+1)^2+t^4} \rho_0^4(r_0^2+1)^2,\\
\M_2&=\frac{r_0}{r_0^2+1}(n+\epsilon \M_1),
\end{align*}
where $n,m$ are integers. The curve $\gamma(\E,\M_1,\M_2)$ is then a closed geodesic in $T^2_{\rho_0,r_0}$, where $\E,\M_1,\M_2$ depend on $\rho_0,r_0,n,m$ as explained above. In particular, there must be at least countably many closed geodesics in $T^2_{\rho_0,r_0}\subset M^3_{\Gamma}$.

\section{Integrals of subharmonic functions on the hypersurfaces $M^3_\Gamma$}
\label{sec:4}

%Since in case of a closed curve $\Gamma$ and $t\not=0$ the hypersurfaces $(M^3_\Gamma,h_t)$ are open complete manifolds, every $\L^2$--harmonic function must be constant (see e.g. \cite{eichhorn}) and thus vanish since the volume of the hypersurfaces $\vol(M^3_\Gamma)$ is infinite, so that the $\L^p$--kernel of the Laplacian on $M^3_\Gamma$ becomes trivial for all $p \geq 1$.

In this section we will show that the $\L^p$--kernel of the Laplacian on the hypersurfaces $(M^3_\Gamma,h_t)$ becomes trivial for all $p\geq1$, where $t\geq 0$ and $\Gamma$ are arbitrary.
We will base our considerations on the much more general work of Greene and Wu \cite{greene-wu}, who studied integrals of certain generalized subharmonic functions on connected non--compact Riemannian manifolds admitting a canonical exhaustion function and showed that these integrals cannot be bounded. More precisely, they showed that the following theorem holds.

\begin{theorem}[Greene and Wu]
\label{thm:greene-wu}
Let $M$ be a connected non--compact oriented $\Cinft$ Riemannian manifold. Suppose that there exists a continuous proper function $\phi:M\rightarrow \R$ and a compact set $K_\phi\subset M$ such that

a) $\phi|_{M\setminus K_\phi}$ is ${\rm{C}}^2$.

b) $\phi|_{M\setminus K_\phi}$ is uniformly Lipschitz continuous.

c) $\phi|_{M\setminus K_\phi}$ is subharmonic.

Denote by $\Sigma(M)$ the closure of the set of all $\Cinft$ subharmonic functions in ${\rm{C}}^0(M)$. Then, if  $f$ is a nonnegative function in $\Sigma(M)$ such that 
$$\{p \in M: f(p) >0,\, \phi(p)>\max_{K_\phi} \phi,\, \grad \phi(p) \not=0\}\not=\emptyset,$$ there exist constants $A_f>0$ and $\tau_0$ such that 
\begin{displaymath}
\int\limits_{M^\phi_\tau} f dM \geq A_f(\tau-\tau_0)
\end{displaymath}
for all $\tau\geq \tau_0$, where $M^\phi_\tau$ denotes the set of all $p \in M$ such that $\phi(p)\leq \tau$; in particular, $\int_M f dM=+\infty$. 
\end{theorem}

A description of the set $\Sigma(M)$ is given by the following proposition.

\begin{proposition}[Greene and Wu]
\label{prop:greene-wu}
Let $M$ be a non--compact $\Cinft$ Riemannian manifold. Then the following functions are in $\Sigma(M)$:

1) Any function $f:M\rightarrow \R$ that is the limit uniformly on compact subsets of $M$ of a sequence of functions in $\Sigma(M)$,

2) ${\rm{C}}^2$ subharmonic functions,

3) $u^p$ where $u$ is a ${\rm{C}}^2$ nonnegative subharmonic function and $p\geq1$,

4) $|u|^p$ where $u$ is a harmonic function and $p\geq1$,

5) any geodesically convex function.
\end{proposition}

In general the scalar Laplacian on a Riemannian manifold $(M^n,g)$, acting on $\Cinft$ functions,  is given by $\Delta f=-\div(\grad f)$, where for a vector field $X \in \X(M^n)$ its divergence with respect to an orthonormal frame $\mklm{Y_1,\dots,Y_n}$ is given by
\begin{displaymath}
\div(X)= \sum \limits_{i=1}^n  g(\nabla _{Y_i}X,Y_i)=\sum \limits_{i=1}^n Y_i(X^i)+\sum\limits_{i,j=1}^n X^j \omega_{ji}(Y_i).
\end{displaymath}
Here the $X^i$ denote the components of  $X$ and the $\omega_{ij}$ the connection forms of the   Levi--Civita connection $\nabla$ of  $M^n$.
In the following we will show that the above results also apply for the considered hypersurfaces $M^3_\Gamma$, $\Gamma$ being arbitrary, obtaining in particular the vanishing of the $\L^p$--kernel of the Laplacian even in case $M^3_\Gamma$ is not complete.
Let us first start considering the function
\begin{displaymath}
\phi^\ast:=\rho \sqrt{r^2+1},
\end{displaymath}
which is $\Cinft$ on $M^3_\Gamma \cap M^4$. One calculates with respect to the orthonormal frame \eqref{eq:2.1}
\begin{displaymath}
Y_1(\phi^\ast)=\frac 1{\sqrt{h_{22}}} \sqrt{r^2+1}, \qquad Y_1Y_1(\phi^\ast)=-\frac  1{2{h_{22}}}  (\log K)_{,\rho} \sqrt{r^2+1},
\end{displaymath}
and thus
\begin{align*}
\Delta \phi^\ast &= -Y_1Y_1 (\phi^\ast)-\sum\limits_{i=1}^3Y_1(\phi^\ast) \omega_{1i}(Y_i)=\frac  1{{h_{22}}}\Big ( \frac{1}2 (\log K)_{,\rho} -\frac 2{\rho} \Big ) \sqrt{r^2+1}=\\
&=\frac{\sqrt{\rho^4(r^2+1)^2+t^4}}{2\rho^3(r^2+1)^{3/2}} \Big (\frac{t^4}{\rho^4(r^2+1)^2+t^4} -2 \Big).
\end{align*}
Because of $\sup_\rho t^4(\rho^4(r^2+1)^2+t^4)^{-1}=1$ it follows that $\phi^\ast$ is subharmonic and one computes further that
\begin{equation*}
|\grad \phi^\ast|^2=Y_1^2(\phi^\ast)=K^{-1}=\frac{\sqrt{\rho^4(r^2+1)^2+t^4}}{2\rho^2(r^2+1)}.
\end{equation*}
We define now the $\Cinft$ function $\lambda: \R \rightarrow [0,1)$ by  $\lambda(x)=e^{-1/x^{2}}$ for $x>0$ and $\lambda(x)=0$ for $x\leq0$ and put
\begin{equation}
\label{eq:4.2}
\mu(x):=\int\limits^x_{-\infty} \lambda(y)\lambda(1-y) dy \Big / \int\limits^{+\infty}_{-\infty} \lambda(y)\lambda(1-y) dy.
\end{equation}
The function $\mu:\R \rightarrow [0,1]$  is $\Cinft$ too, monotone, equal to zero for $x\leq 0$ and one for $x\geq1$. Let $0 < s_0 < L_\Gamma$. Then
\begin{equation*}
\phi:= \rho\sqrt{r^2+1} \cdot \mu(\rho)
\end{equation*}
is $\Cinft$ on $M^3_\Gamma$ and subharmonic on $M^3_\Gamma\setminus K_\phi$ where
\begin{equation*}
K_{\phi}:=\mklm{p=p(s,\rho,\phi) \in M^3_\Gamma: \rho \leq 1, s \leq s_0}.
\end{equation*}
Note that $K_\phi$ is compact and that $\phi$ is proper, i.\ e.,
\begin{displaymath}
\phi^{-1}[0,\kappa]=\mklm{p \in M^3_\Gamma: \, \rho\sqrt{r^2+1} \cdot \mu(\rho) \leq \kappa}
\end{displaymath}
is compact for all $\kappa \in \R$. 
We show that $\phi$ is globally Lipschitz. In order to do so, let us first remark that $|\grad \phi|\leq B_\phi$ on $M^3_\Gamma$, where $B_\phi$ is a constant, since $|\grad \phi|^2$ tends asymptotically to $1/2$ on $M^3_\Gamma \setminus K_\phi$ and, as a smooth function, remains bounded on $K_\phi$. Now let $p$ and $q$ be two points on $M^3_\Gamma$, and $\gamma(\tau)$ the shortest geodesic between them so that $\text{dist}(p,q)=L_\gamma$; we assume that $\gamma$ is parametrized by arc length. Since
\begin{displaymath}
h_t(\grad \phi(\gamma(\tau_0)),\dot \gamma(\tau_0))=d\phi(\gamma(\tau_0))(\dot \gamma(\tau_0))=\frac d{d\tau} \phi(\gamma(\tau))|_{\tau=\tau_0},
\end{displaymath}
one has by Cauchy--Schwarz that
\begin{align*}
|\phi(p)-\phi(q)|&=\Big|\int_0^{L_\gamma} \frac d{d\tau} \phi(\gamma(\tau)) d\tau\Big|=\Big|\int_0^{L_\gamma} h_t(\grad \phi(\gamma(\tau)),\dot \gamma(\tau)) d \tau\Big|\\
&\leq \int_0^{L_\gamma} |\grad \phi(\gamma(\tau)|\cdot |\dot \gamma(\tau)| d\tau \leq B_\phi \, \text{dist}\,(p,q),
\end{align*}
i.e., $\phi$ is uniformly Lipschitz continuous on $M^3_\Gamma$. Summing up we obtain the following proposition.

\begin{proposition}
\label{prop:4.1}
On the connected non-compact oriented $\Cinft$ Riemannian manifolds $(M^3_\Gamma, h_t)$ there exists, for every $t\not=0$ and every curve $\Gamma$, a proper continuous function $\phi:M^3_\Gamma \rightarrow \R$ and a compact set $K_\phi \subset M^3_\Gamma$ such that 

a) $\phi$ is $\Cinft$,

b) $\phi$ is uniformly Lipschitz,

c) $\phi|_{M^3_\Gamma\setminus K_\phi}$ is subharmonic. 

In particular the conclusions of Theorem \ref{thm:greene-wu} hold.
\endproof
\end{proposition}

Note that the above proposition is also true in case $t=0$, i.e., for the non--complete $\Cinft$ Riemannian manifolds $(M^3_\Gamma\cap M^4, h_0)$. As a consequence of the proposition we obtain the following vanishing theorem.

\begin{corollary}
\label{cor:4.1}
Let $p \geq 1$. There exist no $\L^p$--harmonic functions,  on the hypersurfaces $(M^3_\Gamma,h_t)$ for arbitrary $t\in\R$ and curves $\Gamma$. 
\end{corollary}
\begin{proof}
Let $\phi$ and $K_\phi$ be given as in the previous proposition and let $u$ be a harmonic function on $M^3_\Gamma$. By Proposition \ref{prop:greene-wu} one has $|u|^p \in \Sigma(M^3_\Gamma)$ for all $p \geq 1$. Now, by the Aronszajn--Cordes uniqueness theorem for second order differential operators $\sum_{|\alpha|\leq 2} a_\alpha(x)\gd ^\alpha$ with elliptic metric principal symbol \cite{aronszajn} $u$ cannot vanish identically on $M^3_\Gamma\setminus K_\phi$ unless it vanishes everywhere. Therefore, for $u$ not being trivial, the set $\mklm{m \in M^3_\Gamma: |u|^p(m) > 0,\, \phi(m)>\max_{K_\phi} \phi,\, \grad \phi(m)\not=0}$ is not empty, and by Theorem \ref{thm:greene-wu} there exist constants $A$ and $\tau_0$ such that
\begin{displaymath}
\int_{M^{\phi}_\tau} |u|^p dM^3_\Gamma \geq A(\tau-\tau_0)
\end{displaymath}
for all $\tau\geq \tau_0$. In particular $\Ker_{\L^p}(\Delta)=\mklm{0}$ for all $p\geq1$. 
\end{proof}

\section{Einstein and  $T$--Killing spinors on the hypersurfaces $M^3_\Gamma$}%{The Einstein--Dirac system and $T$--Killing spinors on the hypersurfaces $M^3_\Gamma$}
\label{sec:5}

In the sequel we will consider  the Dirac operator  $\Dirac$ on the hypersurfaces $M^3_\Gamma$, whose geometry has been studied in the previous sections. For $t>0$, the homotopy type of $M^3_\Gamma$ is given by $\R^2 \times \Gamma/\mklm{\pm 1}$. If the curve $\Gamma$ is not closed, $M^3_\Gamma$ cannot be complete and admits only one spin structure. Otherwise $M^3_\Gamma$ has the same homotopy type as the circle $S^1$ and, consequently, admits two spin structures. The trivial spin structure is characterized by the fact that there exists a global trivialization of the $\Spin(3)$--principal bundle covering an arbitrary orthonormal frame bundle, while the non--trivial spin structure admits a trivialization of this kind only locally. On the other hand, the unique spin structure of the Eguchi--Hanson space $H^2$ induces a spin structure on the hypersurface $M^3_\Gamma\subset H^2$ by reduction of the former  with respect to the normal vector field of $M^3_\Gamma$. It turns out that the induced spin structure is the trivial one if and only if the winding number of the closed curve $\Gamma$ is even. In the following most of the results will be derived for the induced spin structure, though some of them that follow from purely geometric arguments hold for both spin structures.

 First we will try to determine solutions to the Dirac equation that are also solutions to the Einstein equation and we will show that the aforementioned hypersurfaces do not admit such solutions in case $t\not=0$. Nevertheless, it is possible to construct such solutions explicitly by deformation into the singular situation, though these solutions are no longer complete. In the complete case and if $M^3_\Gamma$ is a minimal surface, one can further show the existence of a spinor field satisfying a generalized Killing equation for spinors. 

Let $e_1, \dots, e_n$ denote the standard basis of the Euclidean space $\R^n$ and introduce the complex two-dimensional matrices
\begin{equation*}
g_1=\left ( \begin{array}{cc} i & 0 \\ 0 & -i \end{array} \right ), \quad g_2=\left ( \begin{array}{cc} 0 & i \\ i & 0 \end{array} \right ), \quad E=\left ( \begin{array}{cc} 1 & 0 \\ 0 & 1 \end{array} \right ), \quad T=\left ( \begin{array}{cc} 0 & -i \\ i & 0 \end{array} \right ).
\end{equation*}
In case $n=2m$, the spin representation of the $n$--dimensional complex Clifford--algebra $C^c_n$ is given by the isomorphism
\begin{gather*}
\kappa_{2m}: C^c_{2m} \simeq End(\Delta_{2m}), \quad 
\kappa_{2m}(e_j):= E \otimes \dots \otimes E \otimes g_{\alpha(j)}\otimes T \otimes \dots \otimes T,
\end{gather*}
where $ j=1,\dots,2m$ and $\alpha(j)$ is equal to $1$ and $2$ for $j$ odd and even, respectively. For $n=2m+1$ one has the representation%s
\begin{gather*}
\kappa_{2m+1}: C^c_{2m+1} \simeq End(\Delta_{2m+1})\oplus End(\hat \Delta_{2m+1})\stackrel{\,pr_1\,}{\longrightarrow} End(\Delta_{2m+1}), \\
\kappa_{2m+1}(e_j):=\kappa_{2m}(e_j), \qquad \kappa_{2m+1}(e_{2m+1}):=-iT\otimes \dots \otimes T,
\end{gather*}
%as well as 
%\begin{gather*}
%\hat\kappa_{2m+1}: C^c_{2m+1} \simeq End(\Delta_{2m+1})\oplus End(\hat \Delta_{2m+1})\stackrel{\,pr_2\,}{\longrightarrow} End( \Delta_{2m+1}), \\
%\kappa_{2m+1}(e_j):=\kappa_{2m}(e_j), \qquad \kappa_{2m+1}(e_{2m+1}):=-iT\otimes \dots \otimes T,
%\end{gather*}
where $\Delta_{2m}=\Delta_{2m+1}=\hat \Delta_{2m+1}=\C^{2^m}$ denote the corresponding representation spaces as well as the representations itself. The induced representations of $\Spin(n)\subset C^c_n$ will be denoted by the same symbols.

We denote by  $\Sigma(M^2_\Gamma)$ or simply $\Sigma$ the spinor bundle considered in each case of $M^3_\Gamma$, by $\eklm{\cdot,\cdot}$ its hermitean inner product and by $\Gamma(\Sigma)$ the space of smooth sections in $\Sigma$. Further we identify  the tangent bundle $TM^3_\Gamma$ and the cotangent bundle $T^\ast M^3_\Gamma$ with the aid of  $h_t$. The Clifford multiplication $TM^3_\Gamma \otimes_\R \Sigma(M^3_\Gamma)\rightarrow 
\Sigma(M^3_\Gamma)$ of a spinor and a vector can then be extended naturally to a multiplication $\Lambda(M^3_\Gamma) \otimes_\R \Sigma(M^3_\Gamma)\rightarrow  \Sigma(M^3_\Gamma)$ of a spinor and a form.
The Levi--Civita connection $\nabla$ of $(M^3_\Gamma,g_t)$ induces a covariant derivative in $\Sigma(M^2_\Gamma)$, which we will denote by $\nabla$, too. With respect to an orthonormal frame $\mklm{Y_1,\, Y_2, \, Y_3}$ one has for $\nabla$ the local representation
\begin{displaymath}
\nabla:\Gamma(\Sigma) \longrightarrow \Gamma(T^\ast(M^3_\Gamma) \otimes \Sigma), \qquad \nabla \psi=d\psi+\frac{1}{2} \sum\limits_{i<j}^3\omega _{ij} Y_i \cdot Y_j \cdot \psi,
\end{displaymath}
where the $\omega_{ij}$ are the connection forms of the Levi--Civita connection $\nabla$. The Dirac operator $\Dirac :\Gamma({\Sigma})\rightarrow\Gamma({\Sigma})$ on $M_\Gamma^3$  is then locally given by
\begin{displaymath}
\Dirac\psi= \sum\limits_{i=1}^3 Y_i \cdot \nabla_{Y_i} \psi,
\end{displaymath}
where $X \cdot \psi$ denotes the Clifford multiplication of a vector field with a spinor; in the realization of  the complex Clifford algebra $C_3^c\simeq M(2,\C)\oplus M(2,\C)$ given above, the vectors $Y_1,\, Y_2, \, Y_3$ are represented by the matrices $g_1,\,  g_2,\, -iT$, respectively.
Note that   in the three--dimensional Clifford algebra it hold that $e_i=\epsilon _{ijk} \,e_j e_k$, where $\epsilon_{ijk}$ denotes the totally skew symmetric tensor. With respect to the global trivialization \eqref{eq:2.1} the $1$--forms $\omega_{ij}$ have been computed in Proposition \ref{prop:2.2}. Let us now introduce the following definitions.
\begin{definition}
A non--trivial spinor field $\psi$ on a Riemannian spin manifold $(M^n,g)$ with $n \geq 3$ is called a \emph{positive} resp.\ \emph{negative Einstein spinor} with eigenvalue $\lambda \in \R$ if it is a solution of the Dirac equation and the Einstein equation 
\begin{equation*}
\Dirac\psi=\lambda \psi, \qquad \Ric-\frac{1}2 S g=\pm \frac{1}4 T_\psi,
\end{equation*}
where $T_\psi(X,Y)=\Re\eklm{X\cdot \nabla_Y\psi+Y \cdot \nabla_X \psi,\psi}$ is the symmetric $(0,2)$--tensor field defined by $\psi$, the energy momentum tensor of $\psi$.
\end{definition}
As shown in \cite{friedrich-kim}, in dimension $n=3$ and in case that the scalar curvature does not vanish, the existence of an Einstein spinor is equivalent to the existence of a so--called WK spinor:
\begin{definition} Let $(M^n,g)$ be  a Riemannian spin manifold whose scalar curvature $S$ does not vanish anywhere. A non--trivial spinor field on $M$ satisfying  the field equation
\begin{align}
\label{eq:5.1}
\begin{split}
2(n-1)S\,\nabla_X\psi&={n} \,dS(X) \psi +{2\lambda}\frac{n-1}{n-2}\big(2 \Ric(X)-S \,X\big) \cdot \psi+  X \cdot dS \cdot \psi
\end{split}
\end{align}
is called a \emph{weak Killing spinor} or  \emph{WK spinor} with \emph{WK number} $\lambda\in \R$. 
\end{definition}

For general $n$ each solution $\psi$ of the field equation \eqref{eq:5.1} with $\lambda S <0$ and $\lambda S>0$  corresponds to a positive and  negative Einstein spinor with eigenvalue $\lambda$, respectively.
For the existence of a WK spinor the following necessary condition is known \cite{friedrich-kim}: 

\begin{proposition}[Friedrich and Kim]
\label{prop:5.1}
Let $(M^n,g)$ be a Riemannian spin manifold with non--vanishing  scalar curvature and $\psi$ a WK spinor on $(M^n,g)$ with WK number $\lambda$. Then
\begin{align}
\label{eq:5.2}
\begin{split}
4(n-1)\lambda^2[(n^2-5n&+8)S^2-4\modulus{\Ric}^2]\\&=(n-2)^2[(n-1)S^3+n\modulus{dS}^2+2(n-1)S(\Delta S)].
\end{split}
\end{align}
\end{proposition}
%As previously remarked, $dS\in \Gamma(T^\ast(M^3_\Gamma))$ can be understood as an element in $\Gamma(TM^3_\Gamma)$, identifying it with its gradient:
%\begin{align*}
%dS&=Y_1(S) \omega^1+Y_2(S) \omega^2+Y_3(S) \omega^3 \equiv Y_1(S) Y_1+Y_2(S)Y_2+Y_3(S) Y_3=\grad S.
%\end{align*}
We show in the following that, for $t\not=0$, the condition \eqref{eq:5.2} cannot be fulfilled on $M^3_\Gamma$ for any  choice of the curve $\Gamma$.

\begin{proposition}
\label{prop:5.2}
For $t\not=0$ and for any spin structure the hypersurfaces $(M^3_\Gamma,h_t)$ do not admit solutions of the WK equation and, hence, there can be no solution to the Dirac--Einstein system.  
\end{proposition}

\begin{proof}
Assume that a WK spinor with WK number  $\lambda$ is given  on  $(M^3_\Gamma,h_t)$. Then, by Proposition \ref{prop:5.1}
\begin{equation}
\label{eq:5.3}
8 \lambda ^2(2S^2-4\modulus{\Ric}^2)=2S^3+3 \modulus{dS}^2+4S(\Delta S)
\end{equation}
must hold, where $S$ has been computed in Theorem \ref{thm:2.1} Using the relation $(r^2+1)S_{,s}=r\dot r \rho S_{,\rho}$ one computes with respect to the trivialization \eqref{eq:2.1} 
\begin{align*}
Y_1(S)&=\frac{1}{\sqrt{h_{22}}} \frac{\gd}{\gd \rho} S=-\frac{1}{\sqrt{h_{22}}} \frac{\rho^3(r^2+1)^2}{(\rho^4(r^2+1)^2+t^4)^{5/2}}(4t^4-2\rho^4(r^2+1)^2),\\
Y_2(S)&=\frac{1}{\sqrt{h_{33}}}\frac{\gd}{\gd \phi} S= 0,\\
Y_3(S)&=D \frac{\gd}{\gd s}S+E \frac{\gd}{\gd \rho}S+F \frac{\gd}{\gd \phi}=0,
\end{align*}
thus obtaining for the Laplacian of $S$ that 
\begin{align*}
-\Delta S&= \div \grad S=Y_1Y_1(S)+\sum\limits_{i=1}^3 Y_1(S) \omega_{1i}(Y_i)\\
&=Y_1Y_1(S)+Y_1(S) \frac{1}{\sqrt{h_{22}}} \left ( \frac{1}2 (\log h_{33})_{,\rho}-(\log D)_{,\rho} \right )=Y_1Y_1(S)+Y_1(S) \frac{2}{\sqrt{h_{33}}},
\end{align*}
compare Proposition \ref{prop:2.2}. One computes further that 
\begin{align*}
Y_1Y_1(S)&=\frac{1}{\sqrt{h_{22}}} \left ( \frac{1}{\sqrt{h_{22}}} S_{,\rho}\right)_{,\rho}=\frac{1}{h_{22}}S_{,\rho^2}-\frac{1}{2(h_{22})^2} h_{22,\rho} S_{,\rho}\\
%&=-3\frac{\rho^8(r^2+1)^4-7\rho^4t^4(r^2+1)^2+2t^8}{(\rho^4(r^2+1)^2+t^4)^3}+t^4\frac{2t^4-\rho^4(r^2+1)^2}{(\rho^4(r^2+1)^2+t^4)^3}
&=-\frac{2(4t^8-20(1+r^2)^2t^4\rho^4+3(1+r^2)^4\rho^8)}{(\rho^4(r^2+1)^2+t^4)^3}\\
\intertext{as well as}
Y_1(S)\frac{2}{\sqrt{h_{33}}}&= -2 \frac{2t^4-\rho^4(r^2+1)^2}{(\rho^4(r^2+1)^2+t^4)^2},\\
\intertext{obtaining for $\Delta S$ the expression}
\Delta S&=\frac{8t^8-18\rho^4 t ^4(r^2+1)^2+\rho^8(r^2+1)^4}{(\rho^4(r^2+1)^2+t^4)^3}.
\end{align*}
Since $\modulus{dS}^2=Y_1(S)^2$, $\modulus{\Ric}^2=\tr \Ric^2=R^2_{11}+R^2_{22}+R^2_{33}$, one obtains 
\begin{gather*}
%8 \lambda^2\left ( 2 \frac{ \rho^8(r^2+1)^4}{(\rho^4(r^2+1)^2+t^4)^3}- \frac{4t^8+(2t^4-\rho^4(r^2+1)^2)^2+(4t^4+\rho^4(r^2+1)^2)^2}{(\rho^4(r^2+1)^2+t^4)^3}\right)\\=
8 \lambda^2\frac{1}{(\rho^4(r^2+1)^2+t^4)^3}(-24 t ^8 -4 \rho^4 t^4(r^2+1)^2)
\end{gather*}
for the left--hand side of \eqref{eq:5.3} and 
\begin{align*}
%&-2 \frac{\rho^{12}(r^2+1)^6}{(\rho^4(r^2+1)^2+t^4)^{9/2}}+3 \frac{2\rho^4(r^2+1)^2(2t^4-\rho^4(r^2+1)^2)^2}{(\rho^4(r^2+1)^2+t^4)^{9/2}}\\
%&-4 \frac{ \rho ^4 (r^2+1)^2}{(\rho^4(r^2+1)^2+t^4)^{3/2}}\cdot \frac{8t^8-18\rho^4 t ^4(r^2+1)^2+\rho^8(r^2+1)^4}{(\rho^4(r^2+1)^2+t^4)^3}\\&=
\frac{-2}{(\rho^4(r^2+1)^2+t^4)^3} S[-8t^8+48 \rho ^4(r^2+1)^2 t ^4 ],
\end{align*}
for the right--hand side, so that the condition \eqref{eq:5.3} reads 
\begin{equation*}
\lambda^2(-12 t ^8-2 \rho ^4 (r^2+1)^2 t^4)= S( t^8-6 \rho ^4 ( r^2+1)^2 t ^4),
\end{equation*}
and one sees that in case $t\not=0$, it cannot be satisfied for any choice of  the curve $\Gamma$. Note that since the integrability condition \eqref{eq:5.3} is purely geometric, the assertion of the proposition holds for any spin structure.
\end{proof}

Only for $t=0$ the condition \eqref{eq:5.3} is  fulfilled for arbitrary values of $\lambda$, since then both sides vanish.
In this case the hypersurfaces $(M^3_\Gamma,h_0)$ are no longer complete for any curve, the metric becoming degenerate along the exceptional curve; one finds that $K=G=2,\, H=0$ and the Ricci tensor and the scalar curvature are
 \begin{gather}
\label{eq:5.3a}
\begin{split}
 \Ric&=\frac{1}{2 \rho^6(r^2+1)^3} \left ( \begin{array}{ccc} 0 & 0 & 0  \\ 0 & -\rho^4(r^2+1)^2 & 0 \\ 0 & 0 & -\rho^4(r^2+1)^2 \end{array} \right ),\\
  S&=-\frac{1}{\rho^2(r^2+1)}.
\end{split} 
\end{gather}
In the following we show that, in this case, solutions of the Dirac--Einstein system  can be constructed explicitly on $(M^3_\Gamma\cap M^4,h_0)$ for an arbitrary choice of the curve $\Gamma$.

In order to do so let $\psi$ be a non--trivial spinor field on $M^3_\Gamma$ that satisfies the spinor equation \eqref{eq:5.1} for $n=3$,
\begin{gather*}
\nabla_X\psi=\frac{3}{4S} dS(X) \psi +\frac{2\lambda}{S} \Ric(X) \cdot \psi-\lambda X \cdot \psi+ \frac{1}{4S} X \cdot dS \cdot \psi.
\end{gather*}
Putting
$\psi=\sqrt{-S} \,\chi$, the above equation can be reformulated into an equation for $\chi$. Using $\nabla(f\psi)=d\,f\otimes \psi +f \nabla \psi$ for a function $f$ and the relation
$X\cdot dS=-dS(X)+X \times dS$ in the 3--dimensional Clifford algebra yields
\begin{equation}
\label{eq:5.4}
\nabla _X\,\chi=\lambda\left( \frac{2}{S} \Ric(X)-X \right ) \cdot \chi +\frac{1}{4S}(X \times dS) \cdot \chi.
\end{equation}
As already shown, with respect to the base \eqref{eq:2.1} only $Y_1(S)$ is different from zero and one obtains
\begin{displaymath}
X \times dS=\omega^3(X)Y_1(S) Y_2-\omega^2(X) Y_1(S)Y_3.
\end{displaymath}
Further, one has
\begin{displaymath}
\frac{2}{S} \Ric(X)-X=\frac{2R_{11}-S}S\omega^1(X) Y_1+\frac{2R_{22}-S}S\omega^2(X) Y_2+\frac{2R_{33}-S}S\omega^3(X) Y_3.
\end{displaymath}
In the realization of the complex Clifford algebra given above one then obtains due to Proposition \ref{prop:2.2} that 
\begin{align*}
\nabla \chi&=d\chi+\frac{1}2 \sum \limits_{i<j} \omega_{ij} Y_i \cdot Y_j \cdot \chi=
d\chi+\frac{1}2  \left ( \begin{array}{cc} i \omega_{23} & -\omega_{12}-i\omega_{13} \\ \omega_{12}-i\omega_{13} & -i\omega_{23} \end{array} \right )  \left (\mnd \begin{array}{c} \chi_1 \\ \chi_2 \end{array}\mnd \right )\\
&=d\chi+\frac{1}{2 \sqrt{h_{22}}}\left [\frac{(\log h_{33})_{,\rho}}2\left (\begin{array}{cc} 0 & -\omega^2 \\
 \omega^2  &0 \end{array} \right )+(\log D)_{,\rho}\left (\begin{array}{cc} 0 & i \omega^3 \\ i \omega^3 &0 \end{array} \right )\right ]\left ( \mnd \begin{array}{c} \chi_1 \\ \chi_2 \end{array} \mnd \right ).
\end{align*}
Now, if $t=0$,
\begin{align*}
\frac{2 R_{11}-S}S=-1, \qquad \frac{2 R_{22}-S}S=0,\qquad
\frac{2 R_{33}-S}S=0
\end{align*}
as well as
\begin{align*}
(\log D)_{,\rho}=-\frac{1}{\rho},\qquad (\log h_{33})_{,\rho}=\frac{2}{\rho},\qquad(\log S)_{,\rho}= -\frac 2{\rho}.
\end{align*}
Summing up, \eqref{eq:5.4} now reads
\begin{align*}
\left (\begin{array}{c} d \chi_1 \\ d\chi_2 \end{array} \right )&=\left \{ \lambda \left ( \begin{array}{cc} -i\omega^1 & 0 \\ 0 & i\omega^1 
\end{array} \right ) -\frac{1}{2\rho\sqrt{h_{22}}} \left ( \begin{array}{cc} 0 & i\omega^3+\omega^2 \\ i \omega^3 - \omega^2 & 0 \end{array} \right ) \right.\\
&\left.-\frac{1}{2\rho\sqrt{h_{22}}}\left ( \begin{array}{cc}  0 & -\omega^2-i\omega^3 \\ \omega^2 -i\omega^3 & 0 \end{array} \right ) \right\}\left (\mnd \begin{array}{c} \chi_1 \\ \chi_2 \end{array} \mnd \right )\\
&=\lambda \left ( \begin{array}{cc} -i\omega^1 & 0 \\ 0 & i\omega^1 
\end{array} \right ) \left (\mnd\begin{array}{c} \chi_1 \\ \chi_2 \end{array} \mnd\right ),
\end{align*}
the summands $(X\times dS)/4S$ and $\sum_{i<j} \omega_{ij}(X)\, e_i\cdot e_j/2$ cancelling out each other. Since $d\omega_1=0$, the system above can be integrated.
Taking into account the equality $d \chi_j=\sum_{i=1}^3Y_i(\chi_j)\omega^i=\chi_{j,s}\, ds+ \chi_{j,\rho}\, d\rho+\chi_{j,\phi}\,d\phi$ and the expressions for the $\omega^i$ one derives the system of partial differential equations
\begin{align*}
\frac{\gd}{\gd s}\left (\mnd\begin{array}{c} \chi_1 \\ \chi_2 \end{array}\mnd\right )&=\left( \!\begin{array}{cc} f_0 & 0 \\ 0 & \overline{f_0} \end{array} \!\right )  \left (\mnd\begin{array}{c} \chi_1 \\ \chi_2 \end{array} \mnd\right ), 
& \frac{\gd}{\gd \rho}\left (\mnd\begin{array}{c} \chi_1 \\ \chi_2 \end{array} \mnd\right )&=\left(\! \begin{array}{cc} f_2 & 0 \\ 0 & \overline{f_2} \end{array}\! \right )  \left (\mnd\begin{array}{c} \chi_1 \\ \chi_2 \end{array}\mnd\right ), 
\end{align*}
where 
\begin{align*}
f_0&= -i\lambda \frac{r\dot r\rho}{r^2+1} \sqrt{(r^2+1)K}, & f_2&=-i\lambda \sqrt{(r^2+1)K}.
\end{align*}
are functions in the variables $\rho$ and $s$. Note that $(r^2+1)f_0=r \dot r \rho f_2$.
Further one has 
\begin{displaymath}
\rho \frac{\gd}{\gd s} f_2=\rho \left (-\lambda i \frac{1}2 \sqrt{\frac{K}{r^2+1}} 2 r \dot r\right )=f_0,
\end{displaymath}
showing that 
\begin{displaymath}
\chi_1=e^{f_2(s) \rho}, \qquad \chi_2=e^{-f_2(s)\rho}
\end{displaymath}
is a solution of the system above. Transforming  back to the original WK equation yields the following proposition.

\begin{proposition}
 Consider the family of hypersurfaces $(M^3_\Gamma\cap M^4,h_0)$, where $\Gamma$ is an arbitrary curve. Then
\begin{equation*}
\psi=\frac 1{\rho\sqrt{r^2+1}} \left (\begin{array}{c} e^{-\lambda\sqrt{2(r^2+1)}\rho i} \\
e^{\lambda\sqrt{2(r^2+1)}\rho i}
\end{array}\right )
\end{equation*}
is a WK spinor of length $\modulus{\psi}^2=-S\modulus{\chi}^2=-S(\modulus{\chi_1}^2+\modulus{\chi_2}^2)=-2S$ and WK number $\lambda\in \R$. Thus, the normalized spinor
\begin{equation*}
\sqrt{\frac{-S}{\modulus{\lambda} \modulus{\psi}^2}} \psi=\frac 1{\rho \sqrt{2(r^2+1)|\lambda|}}\left (\begin{array}{c} e^{-\lambda\sqrt{2(r^2+1)}\rho i} \\
e^{\lambda\sqrt{2(r^2+1)}\rho i}\end{array}\right )
\end{equation*}
is an Einstein spinor on $M^3_\Gamma\cap M^4$ with eigenvalue $\lambda$.
\end{proposition}
The homotopy type of $M^3_\Gamma\cap M^4$ is given by $\R^\ast_+ \times S^1\times \Gamma$; therefore it has at least two spin structures, the one involved here being determined by the global trivialization \eqref{eq:2.1}.
Recall  that $M^3_\Gamma\cap M^4$ is parametrized by the length parameter $s$ of the curve $\Gamma(s)=r(s)e^{i\phi_\Gamma(s)} \subset \C$ and the fiber parameters $0<\rho <\infty,\, 0\leq \phi <2\pi$. The metric $h_0$ is then given by the formula
\begin{equation*}
h_0=2\big(\rho^2 ds^2+(r^2+1)(d\rho^2+\rho^2 d\phi^2)\big)+r \dot r \rho\, ds\, d\rho+r^2\dot \phi_\Gamma \rho^2 ds \,d\phi,
\end{equation*}
and the Ricci tensor has rank two, see equation \eqref{eq:5.3a}. Similar examples of WK spinors on  a $3$--dimensional non--complete Riemannian manifold with negative scalar curvature have been constructed in \cite{friedrich-kim}.

We introduce now the notion of a T--Killing spinor \cite{friedrich-kim2}.  
\begin{definition}
Let $(M^n,g)$ be a Riemannian spin manifold. A spinor field $\psi$ without zeros will be called a \emph{T--Killing spinor} if the trace $\Tr(\hat T_\psi)=\frac 1{\norm{\psi}^2} \Tr(T_\psi)$ is constant and $\psi$ is a solution of the field equation
\begin{equation*}
\nabla_X \psi=-\frac 1{2} \hat T_\psi(X) \cdot \psi, \qquad X \in \X(M^n).
\end{equation*}
Here $\hat T_\psi(X,Y)=\frac 1{\norm{\psi}^2} T_\psi(X,Y)$ is the energy momentum tensor of the normalized spinor $\psi/\norm{\psi}$.
\end{definition}

As remarked at the  beginning, $(H^2,g_t)$ is endowed with a hyperk\"ahler structure and therefore Ricci--flat and self--dual. Due to this, there is a parallel spinor on $H^2$, and the study of its restriction to $M^3_\Gamma$ will enable us to construct a $T$--Killing spinor explicitly. There we follow a similar construction carried out  in \cite{friedrich98, friedrich-kim2}, where the restriction of a parallel spinor on the Euclidean space $\R^3$ to an isometrically immersed closed $2$--surface of constant mean curvature is considered, yielding examples of $T$--Killing spinors on any surface of constant mean curvature in $\R^3$. 

We consider first the restriction of the spinor bundle of $H^2$ to the submanifold $M^3_\Gamma$  (compare \cite{baum-friedrich-grunewald-kath}). Note that the Clifford representation $\Delta_{2k+2}$ can be constructed directly from the Clifford representation $\Delta_{2k+1}$ by setting
\begin{gather*}
\Delta_{2k+2}:= \Delta_{2k+1}\oplus\Delta_{2k+1}
\end{gather*}
and defining the Clifford multiplication  in $\Delta_{2k+2}$ by means of the Clifford multiplication in $\Delta_{2k+1}$,
\begin{gather*}
e_i \cdot  (\psi_1\oplus \psi_2):= e_i\cdot  \psi_1 \oplus (- e_i\cdot \psi_2), \qquad 1 \leq i\leq 2k+1,\\
e_{2k+2} \cdot (\psi_1\oplus \psi_2):=\psi_2\oplus (-\psi_1).
\end{gather*}
The mapping 
\begin{displaymath}
f:=i^{k+1} (e_1  \dots  e_{2k+2}): \Delta_{2k+2}\longrightarrow \Delta_{2k+2}
\end{displaymath}
is an automorphism of the corresponding $\Spin(2k+2)$--representation, and because of $(e_1\dots e_{2k+2})^2=(-1)^{k+1}$  it turns out to be an involution. Thus the spin representation $\Delta_{2k+2}$ decomposes into the eigensubspaces of $f$, and we denote them by $\Delta_{2k+2}^{\pm}$. Explicitly one has 
\begin{displaymath}
f(\psi_1\oplus \psi_2)=i^{k+1}  (e_1\dots e_{2k+1}\cdot \psi_2\oplus e_1\dots e_{2k+1}\cdot \psi_1),
\end{displaymath}
yielding in particular for $k=1$ the relation
\begin{displaymath}
f(\psi_1\oplus \psi_2)=-  \big (e_1 e_2 e_{3}\cdot \psi_2\oplus e_1 e_2 e_{3}\cdot \psi_1\big )=\psi_2\oplus \psi_1,
\end{displaymath}
since $e_1e_2e_3=-1$ in the three--dimensional Clifford algebra. In this way one  obtains
\begin{align*}
\Delta_{4}^{\pm}=\mklm{\psi_1\oplus \psi_2 \in \Delta_4: \psi_2=\pm\psi_1},
\end{align*}
i.e., a spinor in $\Delta_{4}^+$ or $\Delta_{4}^-$ uniquely defines a spinor in $\Delta_3$ and vice versa. Thus we have defined  two isomorphisms of $\Spin(3)$ representations,
\begin{equation}
\label{eq:5.3b}
\Delta_3 \simeq \Delta_{4}^{\pm}: \phi_1 \longmapsto \phi_1\oplus (\pm \phi_1).
\end{equation}

Since the four--dimensional spin manifold $(H^2,g_t)$ is simply connected, it has only one spin structure, and we denote the corresponding spinor bundle by $\Sigma_{H^2}$. It splits into the subbundles $\Sigma^+_{H^2}$ and $\Sigma^-_{H^2}$, according to the above decomposition of $\Delta_4$, and as a consequence of  $\Delta_{4}= \Delta_{3}\oplus\Delta_{3}$ and \eqref{eq:5.3b} we have the identifications
\begin{equation*}
{\Sigma_{H^2}}_{|M^3_\Gamma}\simeq\Sigma\oplus\Sigma,\qquad
\Sigma \simeq {\Sigma^\pm_{H^2}}_{|M^3_\Gamma}, 
\end{equation*}
where $\Sigma$ is the induced spinor bundle on $M^3_\Gamma$. Consider now a spinor field $\phi^+\in \Gamma(\Sigma^+_{H^2})$ and its restriction $\phi^+_{|M^3_\Gamma}=\phi_1\oplus \phi_1$  to $M^3_\Gamma$, where $\phi_1 \in \Gamma(\Sigma)$ is a three--dimensional spinor field. In particular note that for a field of unit normal vectors on $M^3_\Gamma$ the relation $N \cdot (\phi_1\oplus \phi_1)=
\phi_1\oplus(-\phi_1)$ holds, according to the realization of $\Delta_4$ given above.
 By using the local formulas for the different covariant derivatives one obtains for the spinorial derivative of $\phi^+$ on $M^3_\Gamma$ the relation
\begin{align*}
\nabla^{\Sigma_{H^2}}_X \phi^+&=d\phi^+(X)+\frac 1 2 \sum_{1\leq i<j \leq 3} \omega_{ij}(X)\, Y_i \cdot Y_j \cdot (\phi_1\oplus \phi_1)\\&+\frac 1 2 \sum_{1\leq i<4} \omega_{i4}(X) \,Y_i \cdot N \cdot  (\phi_1\oplus \phi_1)\\&=(\nabla^\Sigma_X \phi_1\oplus \nabla^\Sigma_X\phi_1)- \frac 1 2 (\nabla ^{H^2}_X N \cdot \phi_1 \oplus \nabla ^{H^2}_X N \cdot \phi_1)
\end{align*}
for every vector field $X \in TM^3_\Gamma$, since $\omega_{ij}(X)=g_t(\nabla^{H^2}_X Y_i,Y_j)=h_t(\nabla^{M^3_\Gamma}_X Y_i,Y_j)$ and $\omega_{i4}(X)=g_t(\nabla^{H^2}_X Y_i,N)=-h_t(Y_i, \nabla^{H^2}_X N)$. Here and until the end of this section $\mklm{Y_1,Y_2,Y_3}$ denotes an arbitrary section in the frame bundle of $M^3_\Gamma$.
% (note that  $\nabla^{H^2}_X N \in TM^4_\Gamma$).
Since one part of the Weyl tensor of the Eguchi--Hanson space $H^2$ vanishes, we can assume that the parallel spinor on $H^2$ is contained in   $\Gamma(\Sigma^+_{H^2})$ and given by $\phi^+$. Hence $\nabla^{\Sigma_{H^2}}_X \phi^+=0$,  and with $\II(X)= \nabla ^{H^2}_X N$ we obtain  the equation
\begin{equation*}
\nabla^\Sigma_X \phi_1=\frac 1 2 \II(X) \cdot \phi_1
\end{equation*}
for the corresponding three--dimensional spinor $\phi_1$.
Further, since $\II$ is a symmetric bilinear form, $\sum_{i=1}^3 Y_i \cdot \II(Y_i)=-\HM$ is a scalar and one obtains 
\begin{equation*}
\Dirac^\Sigma \phi_1=\sum Y_i \cdot \nabla^\Sigma_{Y_i} \phi_1=-\frac \HM {2}  \phi_1;
 \end{equation*}
moreover, $\phi_1$ has constant length because it is given by the  restriction of a parallel spinor.
 We summarize these results in the following lemma.
\begin{lemma}
\label{lemma:5.1}Let $\Sigma$ denote the induced spinor bundle of $M^3_\Gamma$. Then there exists a spinor $\psi^\ast\in\Gamma(\Sigma)$ on $M^3_\Gamma$ with 
$$\nabla ^\Sigma_X \psi^\ast=-\frac{1}2 \II(X) \psi^\ast,\qquad \Dirac\psi^\ast=\\\frac \HM{2} \psi^\ast,\\\qquad \norm{\psi^\ast}=1.$$
\endproof
\end{lemma}
Let now $\psi^\ast$ be given as in the previous lemma. Then $\nabla_X \psi^\ast=-
(\II(X)/2) \cdot  \,\psi^\ast$, so that  
\begin{align*}
\hat T_{\psi^\ast}(X,Z)&=-\frac1{2\norm{\psi^\ast}^2}\Re \eklm{X \cdot \II(Z) \cdot \psi^\ast+ Z \cdot \II(X) \cdot \psi^\ast, \psi^\ast}, \qquad X,Z \in \X(M^3_\Gamma).
\end{align*}
Making use of the relation $\Re \eklm{X \cdot \psi,\psi}=0$, which holds for an arbitrary vector field $X$ and spinor $\psi$, one computes in the base of the $Y_i$
\begin{align*}
\hat T_{\psi^\ast}(X,Z)&=-\frac1{2\norm{\psi^\ast}^2}\Re \eklm{\sum \limits_{i,j,k=1}^3(X^iZ^j
+Z^iX^j) \II^\ast_{jk} Y_i\cdot Y_k \cdot \psi^\ast,\psi^\ast}\\&=-\frac1{2\norm{\psi^\ast}^2}\Re \eklm{-2\sum \limits_{i,j=1}^3X^iZ^j
\II^\ast_{ij} \psi^\ast,\psi^\ast}=\II(X,Z),
\end{align*}
since only the summands with $i=k$ are different from zero.   In particular, one has
\begin{align*}
\Tr(\hat T_{\psi^\ast})=-\frac1{\norm{\psi^\ast}^2} \sum\limits_{i=1}^3 \Re \eklm{Y_i \cdot \II^\ast(Y_i) \cdot \psi^\ast,\psi^\ast}=\Tr \II=\HM,
\end{align*}
and it follows that $\Tr(\hat T_{\psi^\ast})$ is constant if $\HM$ is constant. Since the latter only occurs if $\HM$ vanishes identically, we deduce the following proposition.

\begin{proposition}
\label{prop:5.3}
Denote by $\Sigma$ the induced spinor bundle of $M^3_\Gamma$ and let $M^3_\Gamma$ be a minimal surface, i.\ e., $\Gamma$ a great circle in $S^2$. Then there exists a T--Killing spinor $\psi^\ast\in \Gamma(\Sigma)$ with $\Tr(\hat T_{\psi^\ast})=0$ satisfying the field equation
\begin{equation*}
\nabla_X \psi^\ast=-\frac1{2} \hat T_{\psi^\ast}(X) \cdot \psi^\ast=-\frac1{2} \II(X) \cdot \psi^\ast.
\end{equation*}
For any other choice of the curve $\Gamma$ there are no T--Killing spinors.
\end{proposition}

\section{The spectrum of the Dirac operator}
\label{sec:6}

In this section we will study some properties of  the spectrum $\sigma(\Dirac)$ of the Dirac operator on the hypersurfaces $M^3_\Gamma$, $\Gamma$ being a closed curve, so that $M^3_\Gamma$ is complete.
In general, the Dirac operator $\Dirac$ on a Riemannian spin manifold $(M^n,g)$ is an elliptic formally selfadjoint differential operator of first order and, as a differential operator, closable. If $M$ is complete, $\Dirac$ is essentially selfadjoint as an unbounded operator in $\L^2(\Sigma)$ with domain $\Ctest(M^n,\Sigma)$ and the kernels of $\Dirac$ and $\Dirac^2$ coincide, see e.g. \cite{friedrich}. Here $\L^2(\Sigma)$ is defined as the completion of $\Ctest(M^n,\Sigma)$, the space of sections in $\Sigma$ with compact support, with respect to the norm induced by the scalar product
\begin{displaymath}
(\psi_1,\psi_2)=\int\limits_{M^n} \langle \psi_1(x),\psi_2(x)\rangle dM^n, \qquad \psi_i \in \Ctest(M^n,\Sigma).
\end{displaymath}
One has $\sigma(\Dirac)=\sigma(\overline\Dirac)$. If $M^n$ is complete, $\sigma(\overline{\Dirac})$ is real and consists only of the approximation spectrum since, in this case, 
$\overline{\Dirac}$ has no residual spectrum.  If, additionally,  $M^n$ is non--compact, one has to expect point spectrum as well as continuous spectrum; in particular, we are interested in the essential spectrum of $\overline{\Dirac}$, which is defined by
\begin{displaymath}
\sigma_{\mathrm{ess}}(\overline{\Dirac}):=\mklm{\lambda \in \C: \text{ there is a Weyl sequence for $\lambda$ and $\overline{\Dirac}$}},
\end{displaymath}
and represents the continuous spectrum together with the eigenvalues of infinite multiplicity.
The main result of this section will consist in showing that the infimum of $\sigma (\overline{\Dirac}^2)$ on $(M^3_{A\Gamma},h_t)$, where $\Gamma$ is a closed curve and $A \in \U(2)$, becomes arbitrarily small for arbitrary values of the parameter $t$, and that $0 \in \sigma(\overline \Dirac)$; for $\Gamma$ arising  by a M\"{o}bius transform from a circle in $\C$ with center at the origin, we also show that the $\L^2$--kernel of $\Dirac$ and $\overline \Dirac$ are trivial, thus obtaining $0 \in \sigma_{\mathrm{ess}}(\overline \Dirac)$ in this case. As we use the global trivialization \eqref{eq:2.1}, these results hold for the trivial spin structure. 
\begin{theorem}
\label{thm:6.1}
Let $\Gamma$ be a closed curve and $\overline{\Dirac}$ the closure of the Dirac operator on the hypersurfaces $(M^3_{A\Gamma},h_t)$, endowed with the trivial spin structure, where $t\not=0$ and $A \in \U(2)$. Then, for arbitrary $\delta > 0$, 
$$\inf \{\lambda:\lambda \in \sigma(\overline{\Dirac}^2)\}<\delta, $$ 
and $0 \in \sigma(\overline \Dirac)$. 
\end{theorem}
We will prove these statements by using the $\min$--$\max$ principle. For this, we need the following lemmas.

\begin{lemma}
The $\L^2_{\mathrm{loc}}$--kernel of $\Dirac$ on $(M^3_\Gamma,h_t)$ is non--trivial for arbitrary $t$ and $\Gamma$.
\end{lemma}
\begin{proof}
With respect to the realization of the previous section one has for $\psi \in \Gamma(\Sigma)$ that
\begin{align*}
\nabla_{Y_k} \psi &= %d\psi(Y_k)
%+\frac{1}2  \left ( \begin{array}{cc} i \omega_{23}(Y_k) & -\omega_{12}(Y_k)-i\omega_{13}(Y_k) \\ \omega_{12}(Y_k)-i\omega_{13}(Y_k) & -i\omega_{23}(Y_k) \end{array} \right )\left (\begin{array}{c} \psi_1 \\ \psi_2 \end{array} \right )\\
\left(\mnd\begin{array}{c} d\psi_1(Y_k) \\ d\psi_2(Y_k) \end{array} \mnd\right ) +\frac{1}{4 \sqrt{h_{22}}}\left [(\log h_{33})_{,\rho}\left (\begin{array}{cc} 0 & -\omega^2(Y_k) \\
 \omega^2(Y_k)  &0 \end{array} \right )\right.\\&\left.+2(\log D)_{,\rho}\left (\begin{array}{cc} 0 & i \omega^3(Y_k) \\ i \omega^3(Y_k) &0 \end{array} \right )\right ]\left (\mnd\begin{array}{c} \chi_1 \\ \chi_2 \end{array} \mnd\right ),
\end{align*}
i.e., $\nabla_{Y_1}\psi =d\psi(Y_1)$ for every $\psi$.
Then  $d\psi_j=\sum_{i=1}^3 Y_i(\psi_j)\omega^i$ implies that  the Dirac operator on $M^3_\Gamma\cap M^4$ is given by
\begin{align*}
\Dirac\psi&=\left ( \mnd\begin{array}{c} i Y_1(\psi_1)+iY_2(\psi_2)-Y_3(\psi_2)
 \\ -i Y_1(\psi_2)+i Y_2 ( \psi_1) +Y_3(\psi_1) \end{array} \mnd\right )+\frac{i}{\rho\sqrt{h_{22}}} \left ( \mnd\begin{array}{cc} 1 & 0 \\ 0 & -1 \end{array} \mnd\right )  \left ( \mnd\begin{array}{c} \psi_1  \\ \psi_2  \end{array} \mnd\right ),
\end{align*}
since
\begin{displaymath}
\frac{1}{\sqrt{h_{22}}} \Big( \frac 1{\rho}+\frac1 {2} (\log K)_{,\rho}-(\log D)_{,\rho}\Big)=\frac{2}{\rho\sqrt{h_{22}}}.
\end{displaymath}
By taking $\Sigma = \rho /\sqrt{\det h_t}=1/{2\rho\sqrt{h_{22}}}$ into account one obtains on $M^3_\Gamma\cap M^4$, for the Dirac operator, the system of partial differential equations
\begin{align}
\label{eq:6.1}
\frac1{\sqrt{h_{22}}}\left [i\left ( \frac \gd{\gd \rho} +\frac{1}\rho\right) \psi_1 - \frac1{2\rho } \Omega_I \psi_2\right ]= \lambda \psi_1,\\
\frac1{\sqrt{h_{22}}} \left [-i\left( \frac \gd{\gd \rho} +\frac{1}\rho\right) \psi_2 + \frac1{2\rho } \Omega_{II} \psi_1\right]= \lambda \psi_2,
\end{align}
where
\begin{align*}
\Omega_1&= (r^2+1) K \frac{\gd}{\gd s} - r \dot r \rho K \frac \gd{\gd \rho} -(r^2 \dot \phi_\Gamma K+2i) \frac \gd {\gd \phi}\\ 
\Omega_{II}&=(r^2+1) K \frac{\gd}{\gd s} - r \dot r \rho K \frac \gd{\gd \rho} -(r^2 \dot \phi_\Gamma K-2i) \frac \gd {\gd \phi}.
\end{align*}
Let now $\lambda=0$ and $\psi$ be of the form $
\psi=\psi(\rho,s)=\rho^{-1} \gamma(s)$. Clearly one has then
\begin{displaymath}
i \left( \frac \gd{\gd \rho} + \frac 1 { \rho}\right ) \psi_j(\rho,s)=i\gamma_j(s) \left (-\frac 1 {\rho^2}+\frac 1 {\rho^2} \right ) =0
\end{displaymath}
as well as
\begin{align*}
r \dot r \rho \frac \gd { \gd \rho} \psi_j(\rho,s) &= - \frac {r\dot r} \rho \gamma_j(s),&
(r^2+1) \frac \gd {\gd s} \psi_j(\rho,s)&=\frac{r^2+1}\rho \frac \gd{\gd s} \gamma_j(s).
\end{align*}
Equating these expressions yields the relation
\begin{align*}
\frac \gd{\gd s} (\log \gamma_i)&=-\frac1{2} \frac \gd {\gd s}( \log (r^2+1))
\end{align*}
 for $\gamma$, so that by integration
\begin{displaymath}
\log \gamma_j=-\frac1{2} \log (r^2+1) + \log C_j.
\end{displaymath}
Putting  $\gamma_j=(r^2+1)^{-1/2}C_j$ one sees that all spinors of the form
\begin{displaymath}
\psi=\frac1{\rho\sqrt{r^2+1}} \left ( \begin{array}{c} C_1 \\ C_2 \end{array} \right ) 
\end{displaymath}
are harmonic on $M^3_\Gamma\cap M^4$, where $C_i \in \C$ are constants. Since, further 
\begin{displaymath}
\int _{U} |\psi|^2 dM^3_\Gamma=\int _{U} |\psi|^2 \sqrt{\det h_t}\, ds \wedge d\rho \wedge d\phi,
\end{displaymath}
for an open region  $U \subset M^3_\Gamma$ and $\det h_t=4 \rho^4 h_{22}$, the harmonic spinors $\psi$ are in $\L^2_{\mathrm{loc}}(\Sigma)$.
\end{proof}

\begin{lemma}
\label{lemma:6.3}
Let $\Gamma$ and $t$ be arbitrary. Then there exists a $\L^2_{\mathrm{loc}}$--harmonic spinor $\psi_0$ on $M^3_\Gamma \cap M^4$  which can be approximated pointwise by spinors $\psi_\epsilon \in \L^2(\Sigma)\cap \Gamma(\Sigma)$ depending on a parameter $\epsilon >0$ such that $\Dirac\psi_\epsilon\in \L^2(\Sigma)$.
\end{lemma}

\begin{proof}
To begin with, note that $\sqrt{-S}$ converges pointwise to $1/{\rho\sqrt{r^2+1}}$ as $t\to 0$, and we therefore introduce the function
\begin{equation*}
S_\epsilon:=-\frac{\rho^4 (r^2+1)^2}{(\rho^4(r^2+1)^2+\epsilon^4)^{3/2}}, \qquad \epsilon >0,
\end{equation*}
replacing in $S$ the parameter $t$ of the K\"{a}hler potential by the new parameter $\epsilon$.
One computes
\begin{align*}
\frac \gd {\gd \rho} \sqrt{-S_\epsilon}&= \frac 1 \rho \sqrt{-S_\epsilon} \frac{2\epsilon^4-\rho^4(r^2+1)^2}{(\rho^4(r^2+1)^2+\epsilon^4)},\\
 \left ( \frac \gd {\gd \rho}+ \frac 1 {\rho}\right ) \sqrt{-S_\epsilon}&=\frac 1{\rho} \sqrt{-S_\epsilon} \frac{3\epsilon^4}{\rho^4(r^2+1)^2+\epsilon^4},
\end{align*}
as well as 
\begin{displaymath}
\Omega_I \sqrt{-S_\epsilon}=0, \qquad \Omega_{II}\sqrt{-S_\epsilon}=0,
\end{displaymath}
since $(r^2+1) S_{\epsilon,s}=r \dot r\rho S_{\epsilon,\rho}$. Each other function in $\rho$ and $s$ of the functional dependence $\rho \sqrt{r^2+1}$ is also harmonic with respect to $\Omega_I$ and $\Omega_{II}$. We put 
\begin{equation}
\label{eq:6.1a}
\psi_\epsilon:= \sqrt{-S_\epsilon} e^{-3\epsilon^4\rho\sqrt{r^2+1}} \left ( \mnd\begin{array}{c} C_1 \\ C_2 \end{array} \mnd \right ),\qquad C_i \in \C \quad \text{constant}.
\end{equation}
For  $\epsilon \to 0$ one has then
\begin{displaymath}
\psi_\epsilon \longrightarrow \frac1{\rho\sqrt{r^2+1}} \left ( \mnd\begin{array}{c} C_1 \\ C_2 \end{array} \mnd\right ).
\end{displaymath}
As remarked, $\Omega_I\psi_\epsilon=0,\,\Omega_{II}\psi_\epsilon=0$, so that 
\begin{align*}
\Dirac\psi_\epsilon&=\frac{i}{\sqrt{h_{22}}} \left ( \frac \gd {\gd \rho}+ \frac 1 {\rho}\right )\left [ \sqrt{-S_\epsilon} e^{-3\epsilon^4 \rho\sqrt{r^2+1}} \right ] \left ( \mnd\begin{array}{c} C_1 \\ -C_2 \end{array} \mnd\right )\\
&=\frac{3\epsilon^4i}{\rho \sqrt{h_{22}}}\left (\frac{1}{\rho^4(r^2+1)^2+\epsilon^4}-\rho\sqrt{r^2+1} \right ) \sqrt{-S_\epsilon} e^{-3\epsilon^4 \rho \sqrt{r^2+1}}\left ( \mnd\begin{array}{c} C_1 \\ -C_2 \end{array} \mnd\right ).
\end{align*}
Then one computes, since $\epsilon>0$,
\begin{align*}
\norm{\psi_\epsilon}_{\L^2}^2&=\int (-S_\epsilon) e^{-6\epsilon^4 \rho\sqrt{r^2+1}}(|C_1|^2+|C_2|^2) dM^3_\Gamma\\
&=2\pi C  \int\limits_0^\infty\int\limits_0^{L_\Gamma}(-S_\epsilon) \sqrt{\det h_t}e^{-6\epsilon^4 \rho\sqrt{r^2+1}} ds \wedge d\rho<\infty,
\end{align*}
where $C=|C_1|^2+|C_2|^2$, as well as
\begin{gather*}
\norm {\Dirac\psi_\epsilon}^2_{\L^2}\\
=\int (-S_\epsilon) e^{-6\epsilon^4 \rho\sqrt{r^2+1}}(|C_1|^2+|C_2|^2) \frac{9\epsilon^8}{h_{22}\rho^2} \left (\frac{1}{\rho^4(r^2+1)^2+\epsilon^4}-\rho\sqrt{r^2+1} \right )^2dM^6_\Gamma\\
=2\pi C  \int\limits_0^\infty\int\limits_0^{L_\Gamma} (-S_\epsilon) \sqrt{\det h_t}e^{-6\epsilon^4 \rho\sqrt{r^2+1}}\frac{9\epsilon^8}{h_{22}\rho^2} \left (\frac{1}{\rho^4(r^2+1)^2+\epsilon^4}-\rho\sqrt{r^2+1} \right )^2 ds \wedge d\rho
\end{gather*}
i.\ e., the $\psi_\epsilon$ are $\L^2$--approximations of  $\L^2_{\mathrm{loc}}$--harmonic spinors,
\begin{gather*}
\L^2(\Sigma)\cap \Gamma(\Sigma) \ni \psi_\epsilon \to \psi_0 \in \L^2_{\mathrm{loc}}(\Sigma),
\end{gather*}
$\Dirac \psi_\epsilon$ being in $\L^2$, too.
\end{proof}

In the following we will use the abbreviations
\begin{align*}
p_\epsilon&:=\sqrt{\det h_t} e^{-6\epsilon^4 \rho\sqrt{r^2+1}},& q_\epsilon&:=\frac{9\epsilon^8}{h_{22}\rho^2} \left (\frac{1}{\rho^4(r^2+1)^2+\epsilon^4}-\rho\sqrt{r^2+1} \right )^2;
\end{align*}
for $\epsilon \to 0$ we then have that 
\begin{align*}
-S_\epsilon\, p_\epsilon \rightarrow \frac{\sqrt{\det h_t}}{\rho^2(r^2+1)}, \qquad S_\epsilon \, p_\epsilon\, q_\epsilon \rightarrow 0
\end{align*}
pointwise. Let $\psi_\epsilon$ be as in \eqref{eq:6.1a}.While $\norm{\psi_\epsilon}_{\L^2}$ becomes unbounded for $\epsilon \to 0$, it does not follow that $\norm {\Dirac\psi_\epsilon}_{\L^2} \to 0$  for small $\epsilon$. Nevertheless, we will show that for given $\delta>0$ and $\epsilon$ small enough, $\norm{\Dirac\psi_\epsilon}_{\L^2}/\norm{\psi_\epsilon}_{\L^2}<\delta$, thus proving Theorem \ref{thm:6.1}. For this we have to determine precise estimates for the Rayleigh quotient $\norm{\Dirac\psi_\epsilon}_{\L^2}^2/\norm{\psi_\epsilon}_{\L^2}^2$ from above, where the point is to find bounds not depending on  $\rho$. 

\begin{proof}[Proof of Theorem $\ref{thm:6.1}$]
Let $\psi_\epsilon$ be as in Lemma \ref{lemma:6.3}, equation \eqref{eq:6.1a}.  One has
\begin{gather*}
-S_\epsilon\sqrt{\det h_t}=\frac{2 \sqrt{2}\rho^7 (r^2+1)^3}{(\rho^4 (r^2+1)^2+\epsilon^4)^{3/2}(\rho^4 (r^2+1)^2+t^4)^{1/4}},\\
\frac \gd {\gd \rho}\Big (-S_\epsilon\sqrt{\det h_t}\Big )= \frac{2\sqrt{2}(r^2+1)^3 \rho^6\big[\rho^4(r^2+1)^2(t^4+6\epsilon^4)+7 \epsilon^4t^4\big] }{(\rho^4 (r^2+1)^2+\epsilon^4)^{5/2}(\rho^4 (r^2+1)^2+t^4)^{5/4}}>0,\\
\sup_{\rho >0}(-S_\epsilon\sqrt{\det h_t})=\frac{2\sqrt{2}}{\sqrt{r^2+1}}; 
\end{gather*}
therefore $-S_\epsilon\sqrt{\det h_t}$ is strictly increasing and tends to 
${2\sqrt{2}/\sqrt{r^2+1}}$  as $\rho \to \infty$, so that it seems natural to estimate $\norm {\psi_\epsilon}^2_{\L^2}$ from below according to
\begin{equation*}
\norm {\psi_\epsilon}^2_{\L^2} =2\pi C \int \limits_0^{L_\Gamma}\int\limits_0^\infty (-S_\epsilon  \, p_\epsilon)\, d\rho \wedge ds \geq 2\pi C \int\limits _0^L \inf_{\rho \geq P} (-S_\epsilon  \sqrt{\det h_t}) \int \limits_P^\infty e^{-6\epsilon^4\rho \sqrt{r^2+1}} d\rho \wedge ds.
\end{equation*}
Here $P>0$ is a cutting point to be determined in a convenient manner, such that the resulting lower bound for  $\norm {\psi_\epsilon}^2_{\L^2}$ is as great as possible. A possible choice would be the turning point  $P_w$ of $-S_\epsilon\sqrt{\det h_t}$, which can be calculated by the condition $(-S_\epsilon\sqrt{\det h_t})_{,\rho^2}=0$ by solving the equation of third degree in $u_1^2=\rho^4(r^2+1)^2$ 
\begin{displaymath}
(u_1^2+\epsilon^4)(u_1^2+t^4)[30 u_1^2\epsilon^4 +42 \epsilon^4 t^4]=5 u_1^2[t^4(u_1^2+\epsilon^4)^2+12 \epsilon^4(u_1^2+t^4)^2 ].
\end{displaymath}
Since this turns out to be a little bit involved and does not necessarily lead to optimal estimates, we look for a condition for  $\rho$ instead such that
\begin{equation}
\label{eq:6.2}
0 < \frac \gd {\gd \rho}\Big(-S_\epsilon\sqrt{\det h_t}\Big)|_\rho \leq a \ll 1,
\end{equation}
i.\ e.,
\begin{gather*}
\big[\rho^4(r^2+1)^2  \big] ^{6/5}\, \big[6\epsilon^4(\rho^4(r^2+1)^2+t^4)+t^4(\rho^4(r^2+1)^2+\epsilon^4)\big]^{4/5} \\\leq \left ( \frac a{2\sqrt{2}} \right ) ^{4/5} (\rho^4 (r^2+1)^2 +\epsilon^4)^2(\rho^4 (r^2+1)^2 +t^4).
\end{gather*}
This is fulfilled if
\begin{gather*}
\big [(\rho^4(r^2+1)^2+t^4)(6\epsilon^4+t^4) \big ]^{4/5} \leq
\left ( \frac a{2\sqrt{2}} \right ) ^{4/5} \big [\rho^4 (r^2+1)^2\big ]^{4/5} (\rho^4(r^2+1)^2+t^4),
\end{gather*}
where we assumed $\epsilon < t$.
For small  $\epsilon$ and $t$ this does not represent a much  stronger condition. Again this is assured if
\begin{displaymath}
\left (\frac{2 \sqrt{2}(6\epsilon^4+t^4)}{a}\right ) ^4  \leq [\rho^4(r^2+1)^2]^5,
\end{displaymath}
and we put 
\begin{equation*}
P_a:= \mu(\epsilon,t,a)\, \frac{1}{\sqrt{r^2+1}}, \qquad \text{where} \qquad \mu(\epsilon,t,a):=\sqrt[5]{\frac{2 \sqrt{2}(6\epsilon^4+t^4)}{a}}.
\end{equation*}
Then one calculates 
\begin{equation*}
\inf _{\rho \geq P_a}(-S_\epsilon\sqrt{\det h_t})=-S_\epsilon\sqrt{\det h_t}|_{P_a}=M(\epsilon,t,a) \, \frac1{\sqrt{r^2+1}},
\end{equation*}
the function $M$ being given by
\begin{displaymath}
M(\epsilon,t,a):=\frac{2 \sqrt{2}\mu^7}{(\mu^4+\epsilon^4)^{3/2}(\mu^4+t^4)^{1/4}}.
\end{displaymath}
We remark that, as $\epsilon \to 0$, the functions $\mu$ and $M$  tend to a finite value  that is independent of $\rho$, namely
\begin{displaymath}
\lim _{\epsilon \to 0} \mu(\epsilon,t,a)=\sqrt[5]{\frac{2\sqrt{2}t^4}{a}}, \qquad \lim _{\epsilon \to 0} M(\epsilon,t,a)=2\sqrt{2}\left (\frac {\sqrt[5]{\left(\frac{2 \sqrt{2}}{at} \right ) ^4 }}{\sqrt[5]{\left(\frac{2 \sqrt{2}}{at} \right ) ^4}+1} \right )^{1/4},
\end{displaymath}
the cutting point $P_a$ also remaining finite.
We finally  obtain an estimate for $\norm{\psi_\epsilon}_{\L^2}$ of the form
\begin{align*}
\norm{\psi_\epsilon}^2_{\L^2} & \geq 2\pi C \int\limits_0^{L_\Gamma} \inf _{\rho \geq P_a}(-S_\epsilon  \sqrt{\det h_t}) \int \limits_{P_a}^\infty e^{-6\epsilon^4\rho \sqrt{r^2+1}} d\rho \wedge ds \\ 
&= 2\pi C M(\epsilon,t,a) \int \frac 1{\sqrt{r^2+1}} \,\frac{ e^{-6\epsilon^4 P_a \sqrt{r^2+1}}}{6\epsilon^4 \sqrt{r^2+1}}  ds.
\end{align*}
Note that $M$ tends  to $2\sqrt{2}$ if, additionally, $a\cdot t \to 0$  so that the value of $-S_\epsilon \sqrt{\det h_t}$ at the point $P_a$ becomes arbitrarily close to $\sup_\rho(-S_\epsilon \sqrt{\det h_t})=2\sqrt{2/r^2+1}$. This can always be achieved by choosing $a$ small enough, though for big $t$ the cutting point $P_a$ becomes big, too. Nevertheless, we will see that this is of no relevance for later arguments. For small $t$ we do not lose too much by the above estimate, since then $P_a$ is also small.

We turn now to estimating $\norm{\Dirac\psi_\epsilon}^2_{\L^2}$. First, one has
\begin{gather*}
-S_\epsilon  \, q_\epsilon \, \sqrt{\det{h_t}}\\=9\sqrt{2}\epsilon^8 \rho^3(r^2+1) \frac{\sqrt[4]{\rho^4(r^2+1)^2+t^4}}{(\rho^4(r^2+1)^2+\epsilon^4)^{3/2}}
\left (\frac{1}{\rho^4(r^2+1)^2+\epsilon^4}-\rho\sqrt{r^2+1} \right )^2 < \infty, 
\end{gather*}
and we set
\begin{align*}
\Delta^2&:=\frac{\rho^2}{\rho^4(r^2+1)^2+\epsilon^4} \left (\frac{1}{\rho^4(r^2+1)^2+\epsilon^4}-\rho\sqrt{r^2+1} \right )^2,\\
\Lambda&:=\rho(r^2+1)  \frac{\sqrt[4]{\rho^4(r^2+1)^2+t^4}}{\sqrt{\rho^4(r^2+1)^2+\epsilon^4}},
\end{align*}
which yields $-S_\epsilon  \, q_\epsilon  \, \sqrt{\det h_t}=9\sqrt{2}\epsilon^8 \Lambda\Delta^2$. The function $\Lambda$ vanishes only for $\rho=0$. The zeros of $\Delta$ are $\rho=0$  and the solutions of the equation of fifth degree in $\rho$,
\begin{equation}
\label{eq:6.4}
\rho\sqrt{r^2+1}(\rho^4(r^2+1)^2+\epsilon^4)-1=0.
\end{equation}
Now $\rho\sqrt{r^2+1}$ becomes zero for $\rho=0$ and is strictly increasing;
$(\rho^4(r^2+1)^2+t^4)^{-1}$ is equal  to $t^{-4}$ for $\rho=0$ and strictly decreasing. The equation \eqref{eq:6.4} has therefore exactly one real solution; it is positive and will be denoted in the following by $Q$. Note that $Q$ is greater than 0 and bounded from above by $1/\sqrt{r^2+1}$.
Since  $-S_\epsilon q_\epsilon\sqrt{\det h_t}$ is non--negative and $(-S_\epsilon q_\epsilon\sqrt{\det h_t} )_{,\rho}=9\sqrt{2}\epsilon^8  \Delta (2\Lambda\Delta_{,\rho}+\Delta \Lambda_{,\rho})$,  the numbers $0$ and $Q$ are the only  absolute minima of $-S_\epsilon q\sqrt{\det h_t}$. The absolute value of 
$\Delta$ can then be estimated according to 
\begin{align*}
|\Delta|&=\frac{\rho}{\sqrt{\rho^4(r^2+1)^2+\epsilon^4}}\modulus{\frac 1{\rho^4(r^2+1)^2+\epsilon^4}-\rho \sqrt{r^2+1}} \leq \\
&\leq \left \{ \begin{array}{cc} \frac{\rho^2 \sqrt{r^2+1}}{\sqrt{\rho^4(r^2+1)^2+\epsilon^4}} &\text{for} \, \rho \geq Q, \\
\frac{\rho}{(\rho^4(r^2+1)^2+\epsilon^4)^{3/2}} & \text{for} \, \rho \leq Q. \end{array} \right. 
\end{align*}
The relation 
\begin{displaymath}
\frac{\gd}{\gd \rho} \left (\frac{\rho^2 \sqrt{r^2+1}}{\sqrt{\rho^4(r^2+1)^2+\epsilon^4}} \right )= \frac{2\rho \epsilon^4\sqrt{r^2+1}}{(\rho^4(r^2+1)^2+\epsilon^4)^{3/2}} >0
\end{displaymath}
as well as
\begin{displaymath}
\sup_\rho \frac{\rho^2 \sqrt{r^2+1}}{\sqrt{\rho^4(r^2+1)^2+\epsilon^4}} =\frac 1{\sqrt{r^2+1}}
\end{displaymath}
imply the estimate
\begin{displaymath}
|\Delta| \leq \frac1{\sqrt{r^2+1}} \qquad \text{for} \, \rho \geq Q.
\end{displaymath}
In a similar way one sees by 
\begin{displaymath}
\frac{\gd}{\gd \rho} \left ( \frac \rho{(\rho^4(r^2+1)^2+\epsilon^4)^{3/2}} \right )=\frac{\epsilon^4-5\rho^4(r^2+1)^2}{(\rho^4(r^2+1)^2+\epsilon^4)^{5/2}},
\end{displaymath}
that $\rho{(\rho^4(r^2+1)^2+\epsilon^4)^{-3/2}}$ has a maximum at $\rho_{max}=\frac \epsilon{\sqrt[4]{5}\sqrt{r^2+1}}$ with
\begin{displaymath}
\frac\rho{(\rho^4(r^2+1)^2+\epsilon^4)^{3/2}}_{\big|\rho_{max}}=\frac 1{\epsilon^5\sqrt[4]{5}(6/5)^{3/2}\sqrt{r^2+1}},
\end{displaymath}
and we obtain the estimate
\begin{equation*}
|\Delta| \leq \frac 1 {\epsilon^5\sqrt[4]{6^6/5^5} \sqrt{r^2+1} }\qquad \text{for} \, \rho \leq Q.
\end{equation*}
Now $\Lambda$ tends asymptotically to $\sqrt{r^2+1}$ as $\rho \to \infty$ and one computes
\begin{equation*}
\frac{\gd}{\gd \rho} \Lambda= \frac{(r^2+1)\big [\epsilon^4 t^4+(2 \epsilon^4-t^4)\rho^4(r^2+1)^2 \big ]}{(\rho^4(r^2+1)^2+\epsilon^4)^{3/2}(\rho^4(r^2+1)^2+t^4)^{3/4}},
\end{equation*}
so that for $2 \epsilon^4<t^4$ one sees that $\Lambda$ has a maximum at
$
\rho_{\mathrm{max}}'=\frac{\epsilon t}{\sqrt[4]{t^4-2\epsilon^4}} \frac 1 {\sqrt{r^2+1}};
$
otherwise it is strictly increasing. Inserting $\rho_{\mathrm{max}}'$ in $\Lambda$ we obtain
\begin{displaymath}
\Lambda_{|\rho_{\mathrm{max}}'}=\sqrt{r^2+1}\frac t \epsilon N(\epsilon,t), 
\end{displaymath}
where 
$N(\epsilon,t):=t/ (\sqrt{2}\sqrt[4]{t^4-\epsilon^4})$,
and thus, for $\Lambda$, the estimate
\begin{equation*}
\Lambda \leq \left \{ \begin{array}{cc}\sqrt{r^2+1}  &\text{for} \quad 2 \epsilon^4 \geq t^4, \\ \sqrt{r^2+1}  N(\epsilon,t) \frac t \epsilon & \text{for} \quad 2 \epsilon^4 < t^4. \end{array} \right.
\end{equation*}
As $\epsilon \to 0$, the function $N$ tends  to $1/\sqrt{2}$. Summarizing we find that, under the assumption  that $2 \epsilon^4 < t^4$, $-S_\epsilon q\sqrt{\det h_t}$ can be estimated from above according to
\begin{equation*}
-S_\epsilon q_\epsilon  \sqrt{\det h_t}=9\sqrt{2}\epsilon^8 \Lambda\Delta^2 \leq \left \{ \begin{array}{cc} \frac{9\sqrt{2}N(\epsilon,t)}{\sqrt{r^2+1}} t \epsilon^7 &\text{for} \quad \rho \geq Q, \\ \frac{9\sqrt{2}N(\epsilon,t)}{\kappa^2\sqrt{r^2+1}} \frac t {\epsilon^3} & \text{for} \quad \rho \leq Q,\end{array} \right.
\end{equation*}
where  $\kappa:=\sqrt[4]{6^6/5^5}$; finally we obtain for $\norm{\Dirac\psi_\epsilon}_{\L^2}$, assuming  $\epsilon$ to be small, that 
\begin{align*}
\norm{\Dirac\psi_\epsilon}^2_{\L^2} &\leq 2 \pi C \int\limits_0^{L_\Gamma} \sup_{\rho \leq Q} (-S_\epsilon q_\epsilon  \sqrt{\det h_t}) \int \limits _0^Q e^{-6\rho \epsilon^4 \sqrt{r^2+1}} \, d\rho\wedge ds\\
&+2 \pi C \int\limits_0^{L_\Gamma} \sup_{\rho \geq Q} (-S_\epsilon q_\epsilon  \sqrt{\det h_t}) \int \limits _Q^\infty e^{-6\rho \epsilon^4 \sqrt{r^2+1}} \, d\rho\wedge ds\\
&=2\pi C \int\limits_0^{L_\Gamma} \frac{9\sqrt{2}N(\epsilon,t)t}{\sqrt{r^2+1}}  \left ( \frac 1 {\kappa^2 \epsilon^3}  \frac {1-e^{-6Q\epsilon^4\sqrt{r^2+1}}}{6 \epsilon^4\sqrt{r^2+1}}+\epsilon^7 \frac {e^{-6Q\epsilon^4\sqrt{r^2+1}}}{6 \epsilon^4 \sqrt{r^2+1}}\right ) ds. 
\end{align*}
Under the assumption that $2 \epsilon^4 < t^4$  we obtain  the estimate
\begin{align*}
\frac{\norm{\Dirac\psi_\epsilon}^2_{\L^2}}{\norm{\psi_\epsilon}^2_{\L^2}}&\leq \frac{9\sqrt{2}N(\epsilon,t)t}{M(\epsilon,t,a)}\frac{\int \big ( \epsilon^{-3}\kappa^{-2} (1-e^{-6Q\epsilon^4\sqrt{r^2+1}})+\epsilon^7 \,e^{-6Q\epsilon^4 \sqrt{r^2+1}} \big ) \, ds}{\int  e^{-6P_a \epsilon^4 \sqrt{r^2+1}} \, ds}
\end{align*}
for the Rayleigh quotient. The expression
\begin{align*}
\frac{e^{-6Q\epsilon^4\sqrt{r^2+1}}-1}{\epsilon^3}&= \sum \limits_{k=1}^\infty \frac 1 {k!}(-6Q \sqrt{r^2+1})^k\, \epsilon^{4k-3}
\end{align*}
tends to zero as $\epsilon \to0$, so that the  Rayleigh quotient itself becomes arbitrarily small for $\epsilon \to 0$. Since for closed  curves  $\Gamma$ the 
hypersurfaces $M^3_\Gamma$ are complete, both $\overline{\Dirac}$ and 
$\overline{\Dirac}^2$  are self--adjoint, and by the 
$\min$--$\max$ principle, see e.\ g.\ \cite{reed-simon}, one has
\begin{equation*}
\inf\mklm{\lambda:\lambda \in \sigma(\overline{\Dirac}^2)}=\inf \limits_{0 \not=\psi \in \D(\overline{\Dirac}^2)} \frac{\norm{\overline\Dirac\psi}^2_{\L^2}}{\norm{\psi}^2_{\L^2}},
\end{equation*}
since $\overline \Dirac^2$ is bounded from below. The domain of definition of the closure $\overline{\Dirac}$ of the Dirac operator is given by
\begin{align*}
\D(\overline{\Dirac})=\Big\{&\psi \in \L^2(\Sigma): \exists \, \text{a series} \, \psi_n \in \D(\Dirac) : \psi_n \to \psi\, \text{and}\,\\&  \Dirac \psi_n \, \text{converges in $\L^2(\Sigma)$}\Big\},
\end{align*}
and in case $\psi \in \D(\overline \Dirac) \cap \Gamma(\Sigma)$, one has  $\overline \Dirac \psi=\Dirac\psi$ . The first assertion of the theorem then follows by noting that the inequalities $\int \norm {\psi_\epsilon}^2 dM^3_\Gamma < \infty$, $\int \norm{\Dirac \psi_\epsilon}^2 dM^3_\Gamma < \infty$ and $\int \norm{\Dirac^2 \psi_\epsilon}^2 dM^3_\Gamma < \infty$ imply that $\psi_\epsilon$ lies in $\D(\overline{\Dirac})$ and  $\D(\overline{\Dirac}^2)$, respectively,  since $M^3_\Gamma$ is assumed to be complete.
To see this, let $p_0\in M^3_\Gamma$ be fixed and $\mu(x):\R \rightarrow[0,1]$ be the function defined in \eqref{eq:4.2}. Following \cite{friedrich} we put
\begin{displaymath}
b_n(p):=\mu\Big(2-\frac{\dist(p,p_0)} n\Big), \qquad n=1,2,\dots, \quad p \in M^3_\Gamma.
\end{displaymath}
 Then $b_n\equiv 1$ on  $B_n(p_0)$ and $\supp b_n \subset B_{2n}(p_0)$. Further one sees that  $b_n$ is Lipschitz--continuous and, hence, almost everywhere differentiable with $|\grad b_n|^2 \leq M/n^2$, where $M$ is a constant. Since $M^3_\Gamma$ is complete, the closed envelopes of the geodesic balls $B_n(p_0)$ are compact in $M^3_\Gamma$ and therefore
\begin{displaymath}
\psi_n:= b_n\, \psi_\epsilon \in \D(\Dirac)=\Ctest(M^3_\Gamma,\Sigma).
\end{displaymath}
Since $\norm{\psi_\epsilon}^2_{L^2}<\infty$, one has  $\psi_n \to \psi_\epsilon$ in $\L^2(\Sigma)$. In the same way $\norm{\Dirac \psi_\epsilon}^2_{L^2}<\infty$ implies with the relation
$ \Dirac \psi_n=b_n \Dirac \psi_\epsilon+\grad b_n \cdot \psi_\epsilon$
that $\Dirac\psi_n \to \Dirac\psi_\epsilon$ in $\L^2(\Sigma)$. Consequently, one obtains $\psi_\epsilon \in \D(\overline \Dirac)$, and in a similar way $\Dirac \psi_\epsilon \in \D(\overline \Dirac)$.Finally, by setting $\tilde\psi_\epsilon:=\psi_\epsilon/\norm{\psi_\epsilon}_{\L^2}$, we obtain a sequence of elements in $\D(\overline \Dirac)\cap \Gamma(\Sigma)$ of unit length for which $\|\Dirac \tilde\psi_\epsilon\|_{\L^2} \to 0$ as $\epsilon \to 0$, which implies that $0 \in \sigma_{\mathrm{approx}}(\overline D)=\sigma(\overline D)$.
\end{proof}

In the following we will study the $\L^2$--kernel of the Dirac operator in case that $\Gamma$ is a circle in $\C$ with center at the origin and radius $r\equiv r_0$.  Let  $p=\Psi(s,\rho,\phi)\in M^3_\Gamma\cap M^4$ and $z=e^{i\tau} \in S^1$. As explained in section \ref{sec:8}, in this case 
\begin{align*}
\kappa_z&:M^3_\Gamma\cap M^4 \longrightarrow M^3_\Gamma\cap M^4, \kappa_z(p)=\Psi(s,\rho,(\phi+\tau)\, \text{mod} \,2\pi),\\
\mu_z&:M^3_\Gamma\cap M^4 \longrightarrow M^3_\Gamma\cap M^4, \,\mu_z(p)=\Psi((s+\tau)\, \text{mod} \, 2\pi r,\rho,\phi)
\end{align*}
represent two isometric $S^1$--actions on $M^3_\Gamma\cap M^4$. Putting 
$$(\kappa_z \psi)(p):=\psi( \kappa_{z^{-1}}(p)) \qquad \text{und} \qquad (\mu_z \psi)(p):=\psi( \mu_{z^{-1}}(p)) $$ 
one obtains two continuous unitary $S^1$--representations in $\L^2(\Sigma)$, since by the invariance of the volume form under $\kappa_z$ and $\mu_z$ the equality  
\begin{displaymath}
\int \big\langle\psi(\kappa_{z^{-1}}(p)),\phi(\kappa_{z^{-1}}(p))\big\rangle \,dM^3_\Gamma(p)=\int \big\langle\psi(p),\phi(p)\big\rangle \, (\kappa_z)^\ast (dM^3_\Gamma)(p), \qquad \phi,\psi \in L^2(\Sigma),
\end{displaymath}
and a similar one for $\mu_z$ hold. 
Then,  by the theorem of Stone, there exist uniquely determined self--adjoint operators $M$ and $M_1$ such that $\kappa_{e^{i\tau}}=e^{i\tau M}$, $\mu_{e^{i\tau}}=e^{i\tau M_{1}}$. They are given by $M=i\gd_\phi$, $M_1=i\gd_s$, while the corresponding eigenfunctions are determined by
\begin{displaymath}
i\frac \gd {\gd \phi} e^{i\alpha \phi}=-\alpha  e^{i\alpha \phi} , \qquad i\frac \gd {\gd s} e^{i\alpha \phi_\Gamma}=-\beta \dot \phi_\Gamma e^{i\beta \phi_\Gamma},
\end{displaymath}
where $\alpha$ and $\beta$ are integers and $\dot \phi_\Gamma=\epsilon/r,\, \epsilon=\pm 1$. 
Because of  $\Dirac_{|p}=\Dirac_{|\kappa_z(p)}=\Dirac_{|\mu_z(p)}$, the operator $\Dirac$ commutes with  $\kappa_z$ and $\mu_z$, so that each of the eigensubspaces  $E_\lambda$ of  $\Dirac$ and  $\overline \Dirac$ corresponding to the eigenvalue $\lambda$ decomposes into the eigensubspaces of the unitary $S^1\times S^1$--action according to
\begin{equation*}
E_\lambda=\bigoplus_{\alpha,\beta} \H^\lambda_\alpha\otimes \H^\lambda_\beta,
\end{equation*} 
in concordance with the spectral decomposition of the operators $M$ and $M_1$; in particular, one has
$
\Ker_{L^2}(\overline\Dirac)=\bigoplus_{\alpha,\beta} \H_\alpha\otimes H_\beta.
$ 
A general solution of the Dirac equation $\Dirac \psi=\lambda \psi$ on  $M^3_\Gamma\cap M^4$ can then be written as a product of the form
\begin{equation}
\label{eq:6.5}
\psi(s,\rho,\phi)=e^{i\alpha \phi} e^{i\beta \phi_\Gamma(s)} R(\rho),
\end{equation}
where $R$ is a function of $\rho$.
Thus, the system of partial differential equations  \eqref{eq:6.1} leads to a system of ordinary differential equations
\begin{align}
\label{eq:6.6}
\begin{split}
\frac1{\sqrt{h_{22}}}\left [i\left ( \frac \gd{\gd \rho} +\frac{1}\rho\right) R_1 - \frac i{2\rho } \Big ( (r^2+1)K  \beta \dot \phi_\Gamma-(r^2 \dot \phi_\Gamma K+2i)\alpha\Big)  R_2\right ]= \lambda R_1,\\
\frac1{\sqrt{h_{22}}} \left [-i\left( \frac \gd{\gd \rho} +\frac{1}\rho\right) R_2 + \frac i {2\rho }   \Big ( (r^2+1)K  \beta \dot \phi_\Gamma-(r^2 \dot \phi_\Gamma K-2i)\alpha\Big)R_1\right]= \lambda R_2
\end{split}
\end{align}
for the radial function $R(\rho)$. Introducing $\delta:=\big ((r^2+1)\beta-r^2 \alpha\big) \dot \phi_\Gamma/2$, we put
\begin{align*}
f&:=\big ( \delta K-i\alpha \big)/
\rho, & g&:=\big ( \delta K+i\alpha \big )/\rho,
\end{align*}
and make the substitution   
\begin{displaymath}
\chi:=C \, \rho  \left ( \mnd\begin{array}{c} e^{i\lambda \int \sqrt{h_{22}} d\rho}R_1(\rho) \\ e^{-i\lambda \int \sqrt{h_{22}} d\rho}R_2(\rho) \end{array} \mnd\right ),
\end{displaymath}
so that one obtains for $\chi$ the system of differential equations
\begin{gather*}
\label{eq:6.7}
\frac d {d\rho}\left ( \mnd\begin{array}{c} \chi_1 \\ \chi_2 \end{array} \mnd\right )=\frac 1 \rho \left ( \mnd\begin{array}{c} \chi_1 \\ \chi_2 \end{array} \mnd\right )+
\left (\begin{array}{cc} i\lambda \sqrt {h_{22}} & 0 \\ 0 & - i\lambda \sqrt {h_{22}} \end{array} \right ) \left ( \mnd\begin{array}{c} \chi_1 \\ \chi_2 \end{array} \mnd\right )\\+C \,\rho\left ( \begin{array}{cc}  e^{i\lambda \int \sqrt{h_{22}} d\rho} & 0 \\ 0 &  e^{-i\lambda \int \sqrt{h_{22}} d\rho}\end{array} \right ) \left ( \begin{array}{cc} -1/\rho-i\lambda \sqrt{h_{22}} & f  \\ g  & -1/\rho +i\lambda \sqrt{h_{22}}  \end{array} \right ) \left ( \mnd\begin{array}{c} R_1 \\ R_2 \end{array} \mnd\right )
\\=\left ( \begin{array}{cc} 0 & \tilde f  \\ \tilde g  & 0  \end{array} \right ) \left ( \mnd\begin{array}{c} \chi_1 \\ \chi_2 \end{array} \mnd\right )
\end{gather*}
with $\tilde f:=e^{2i\lambda\int \sqrt{h_{22}} d\rho} f$, $\tilde g:=e^{-2i\lambda\int \sqrt{h_{22}} d\rho} g$. Note that
\begin{align*}
\int \sqrt{h_{22}} d\rho&=\frac{(r^2+1) \rho^2}{\sqrt2(\rho^4(r^2+1)^2 +t^4)^{1/4}} \left ( 1 + \frac {\rho^4(r^2+1)^2}{t^4} \right ) ^{1/4} \F\left (\frac 1 2, \frac 1 4, \frac 4 3, -\frac{\rho^4(r^2+1)^2}{t^4} \right). 
\end{align*}
If $\alpha$ or $\beta$ are different from zero, neither $f$ nor $g$ vanish; differentiating again gives
\begin{align*}
(\chi_1)_{,\rho^2}&=\tilde f_{,\rho} \,\chi_2+\tilde f \,(\chi_2)_{,\rho} =(\log \tilde f)_{,\rho} \,(\chi_1)_{,\rho}+\tilde f \, \tilde g \,\chi_1,\\
(\chi_2)_{,\rho^2}&=\tilde g_{,\rho} \,\chi_1+\tilde g \,(\chi_1)_{,\rho} =(\log \tilde g)_{,\rho} \,(\chi_2)_{,\rho}+\tilde f \, \tilde g \,\chi_2,
\end{align*}
and one obtains the differential equations of second order    
\begin{align}
\label{eq:6.8a}
\frac {d^2\chi_1}{d\rho^2}(\rho)+p(\rho)\,\frac {d\chi_1}{d\rho}(\rho)+q(\rho)\,\chi_1(\rho) = 0, \\
\label{eq:6.8b}
\frac {d^2\chi_2}{d\rho^2}(\rho)+\bar p(\rho)\,\frac {d\chi_2}{d\rho}(\rho)+q(\rho)\,\chi_2(\rho) = 0,
\end{align}
where
\begin{align*}
p(\rho)&=\frac 1 \rho - \delta  \frac{K_{,\rho}}{\delta K-i\alpha}-2i\lambda \sqrt{h_{22}},&
q(\rho)&=-\frac 1 {\rho^2} (\delta^2 K^2+\alpha^2) \leq 0.
\end{align*}
If one puts $\chi_2:=\tilde f^{-1} (\chi_1)_{,\rho}$ and   $\chi_1:=\tilde g^{-1} (\chi_2)_{,\rho}$,respectively, each solution of \eqref{eq:6.8a} or   \eqref{eq:6.8b} corresponds to a solution of the above  system of differential equations for $\chi$, i.\ e., solving  the latter system of two differential equations of first order is equivalent to finding a solution of the differential equation of second order \eqref{eq:6.8a}  or  \eqref{eq:6.8b}. The latter are differential equations of Sturm--Liouville type and our next goal will consist in showing that, for $\lambda=0$ and $\alpha \not=0$, they cannot have any bounded solutions and, in particular, that they do not lead to $\L^2$--integrable solutions $\Psi$ of  the Dirac equation. For this purpose we will make use of the following theorem proved by Hartman \cite{hartman}.

\begin{theorem}[Hartman]
Let $I$ be an interval in $\R$ and $w(x)$ a solution of the differential equation
\begin{displaymath}
\ddot w(x)+p(x) \dot w(x)+q(x) w(x)=0, \qquad x \in I,
\end{displaymath}
with continuous complex valued coefficients $p$ and $q$. If
\begin{equation}
\label{eq:6.9}
\Re \Big[ -q(x)-\frac 1 4 |p(x)|^2\Big] \geq 0,
\end{equation}
then $r(x)=|w(x)|^2$ is concave, i.\ e., $\ddot r(x) \geq 0$.
\end{theorem}
Now, in our case one computes
\begin{align*}
\Re p(\rho) &= \frac 1 2 (p+\bar p)=\frac 1 \rho - \frac {\delta K_{,\rho}} 2 \left (\frac 1 { \delta K-i\alpha}+ \frac 1 { \delta K +i\alpha} \right )\\
&=\frac 1 \rho - \frac { \delta ^2 K K_{,\rho}}{ \delta ^2 K^2+\alpha^2},\\
\Im p(\rho)&= \frac 1 {2i} ( p -\bar p)=-\frac {\delta K_{,\rho}}{2i} \left (\frac 1 { \delta K-i\alpha}- \frac 1 { \delta K +i\alpha} \right )-2\lambda \sqrt{h_{22}}\\
&=-\frac {\alpha \delta K_{,\rho}}{\delta ^2 K^2+\alpha^2} -2\lambda \sqrt{h_{22}}
\end{align*}
and thus
\begin{align*}
-q(\rho)-\frac 1 4 |p(\rho)|^2 &=\frac 1 {\rho^2} \left ( \delta ^2 K^2+\alpha^2-\frac 1 4 \right )+\frac {\delta ^2 K^2 (\log K)_{,\rho}}{\delta ^2 K^2+\alpha^2}\\&\cdot \left ( \frac 1 {2\rho} -\frac{(\log K)_{,\rho}} 4 -\frac {\alpha \lambda \sqrt{ (r^2+1) K}}{\delta K} \right )-\lambda^2 (r^2+1) K.
\end{align*}
Because of   
\begin{displaymath}
\frac 1 {2\rho}-\frac{(\log K)_{,\rho}} 4=\frac 1 {2\rho} \left ( 1 - \frac { t^4}{\rho^4(r^2+1)^2+t^4} \right) > 0
\end{displaymath}
one recognizes that, for $\lambda=0$ and $\alpha \not=0$, the condition \eqref{eq:6.9} is fulfilled for the differential equations \eqref{eq:6.8a} and  \eqref{eq:6.8b}, while for $\lambda \not=0$ the expression $-q(\rho)-|p(\rho)|^2/4$ tends asymptotically to  $-2\lambda^2(r^2+1)$ for $\rho \to \infty$. For $\lambda=0$ and $\alpha =0$ it becomes also negative as $\rho \to 0$. As a consequence of the preceeding theorem we obtain the following lemma.
\begin{lemma}
\label{lemma:6.4}
Assume that $\lambda=0$ and $\alpha\not=0$, and let $\chi_1$, $\chi_2$ be solutions  of the differential equations \eqref{eq:6.8a} and   \eqref{eq:6.8b}, respectively.  Then $|\chi_1|^2$ and $|\chi_2|^2$ are concave.
\endproof
\end{lemma}
We are now in a position to prove the announced theorem.
\begin{theorem}
\label{thm:6.2}
Let $\Gamma$ be a circle in $\C$ with center at the origin and radius $r\equiv r_0$, and $\psi$ a spinor on $(M^3_\Gamma\cap M^4,h_t)$ of the form  \eqref{eq:6.5}. If  $\psi$ is a solution of the Dirac  equation  with respect to the trivial spin structure corresponding  to the eigenvalue $\lambda=0$ and if $\alpha\not=0$, then $\norm{\psi}^2_{\L^2}=\infty$.
\end{theorem}
\begin{proof}
Let  $\Dirac \psi=0$ and $\alpha \not=0$. By our previous considerations $\chi_1=C\rho R_1(\rho)$
satisfies the differential equation \eqref{eq:6.8a} and we consider its continuation  
\begin{equation}
\label{eq:6.10}
\frac {d^2\chi_1}{dz^2}(z)+p(z)\,\frac {d\chi_1}{dz}(z)+q(z)\,\chi_1(z) = 0, \qquad z \in \C,
\end{equation}
to the whole complex domain. For $\alpha \not=0$ both $p(z)$ and $q(z)$, $z\in \C$, are meromorphic functions with poles of first and second order at zero, respectively. The differential equation \eqref{eq:6.10} is therefore of Fuchssian type and zero is a regular singular point. Let $\chi_{1,1}$, $\chi_{1,2}$ form a fundamental system of solutions  of \eqref{eq:6.10}; they can be expanded  around the origin  into the uniformly convergent series 
\begin{align*}
\chi_{1,1}(z)&=z^{\epsilon_1} \Big(1+\sum\limits_{n=1}^{\infty} a_n z^n\Big), & \chi_{1,2}(z)&=z^{\epsilon_2} \Big(1+\sum\limits_{n=1}^{\infty} b_n z^n\Big),
\end{align*}
where $a_n$, $b_n$ are constants and $\epsilon_1$, $\epsilon_2$ are the roots of the equation
\begin{displaymath}
\epsilon^2+(p^0-1)\epsilon +q^0=0
\end{displaymath}
with
\begin{displaymath}
q^0=\lim_{z\to0} z^2 \, q(z), \qquad p^0=\lim_{z\to 0} z \, p(z),
\end{displaymath}
see e.\ g.\ \cite{whittaker-watson}. One obtains $\epsilon_1+\epsilon_2=1-p^0$, $\epsilon_1\epsilon_2=q^0$, which yields  in our case that $\epsilon_1+\epsilon_2=1-1=0,\,  \epsilon_1\epsilon_2=-\alpha^2$, and hence $\epsilon_1=\alpha$, $\epsilon_2=-\alpha$. Evidently, analogous considerations hold for $\chi_2=C\rho R_2(\rho)$, too. Now,
\begin{equation*}
\norm{\psi}^2_{\L^2}=\int \norm{\psi(s,\rho,\phi)}^2 dM^3_\Gamma=\int\limits_0^{2\pi}\int\limits_0^\infty\int\limits_0^{2\pi r} \frac 1 {\rho^2}\big(|\chi_1(\rho)|^2+|\chi_2(\rho)|^2 \big ) \sqrt{\det h_t} \, ds \wedge d\rho \wedge d\phi.
\end{equation*}
In order that the above integral remains bounded it is necessary that $|\chi_1(\rho)|^2$ and $|\chi_2(\rho)|^2$ decrease with order greater than one for $\rho\to \infty$, since 
\begin{displaymath}
\frac 1 {\rho^2} \sqrt{\det h_t}=\frac{\sqrt 8 \rho(r^2+1)}{\sqrt[4]{\rho^4(r^2+1)^2+t^4}} \sim \text{constant};
\end{displaymath}
therefore  $|\chi_i(\rho)|^2 < 1/\rho$, $i=1,2$,  must hold for large $\rho$. As, moreover, $|\chi_i(\rho)|^2$ is smooth, there exists  a $\rho_0$ such that $(|\chi_i(\rho_0)|^2)_{,\rho} \leq -1/\rho_0^2<0$. However, by Lemma \eqref{lemma:6.4} one has that $|\chi_i(\rho)|^2_{,\rho}$ is monotone increasing so that 
\begin{displaymath}
\frac d {d\rho} |\chi_i(\rho)|^2\leq -\frac 1 {\rho_0^2}<0 \qquad \text {for all} \quad \rho \leq \rho_0
\end{displaymath}
must hold. Consequently, $|\chi_i(\rho)|^2$ is monotone decreasing and strictly monotone decreasing for  $\rho \leq \rho_0$. Let us now assume that $\alpha=1,\,2,\, \dots$ without loss of generality. If $\chi(\rho)$ is not identically zero, it follows that, in a neighbourhood of the origin, its components $\chi_1$ and $\chi_2$ must have the developments 
\begin{equation*}
\chi_i(\rho)=A_i\rho^{-\alpha} \Big[1+\sum\limits_{n=1}^\infty c_n^i\rho^n\Big],
\end{equation*}
where $A_i,\, c_n^i$ are constants; otherwise  one would have  $(|\chi_i(0)|^2)_{,\rho}\geq 0$.
Let now $\rho_1$ be sufficiently small so that  $\chi_1$ and $\chi_2$ can be developed as above and, in particular, 
\begin{displaymath}
\Big | \sum\limits_{n=0}^\infty \Re c_n^i \rho^n\Big|<\frac 1 2 \qquad \text{for all} \quad \rho<\rho_1.
\end{displaymath}
Then
\begin{align*}
\int\limits_0^{2\pi}\int\limits_0^{\rho_1}\int\limits_0^{2\pi r} &\norm{\psi(s,\rho,\phi)}^2\, dM^3_\Gamma=4\pi^2 r \int\limits_0^{\rho_1} \rho^{-2(\alpha+1)} \sum\limits_{i=1,2} A_i^2\Big [ \Big (1+ \sum\limits_{n=0}^\infty \Re c_n^i \rho^n\Big)^2\\&+\Big(\sum\limits_{n=0}^\infty \Im c_n^i\rho^n\Big)^2 \Big]\sqrt{\det h_t} \,d\rho\geq \pi^2 r (A_1^2+A_2^2) \int \limits_0^{\rho_1} \rho^{-2(\alpha+1)}  \sqrt{\det h_t}\,d\rho\\
&\geq \pi^2 r (A_1^2+A_2^2) \frac {\sqrt 8 (r^2+1)}{\sqrt[4]{\rho^2_1(r^2+1)^2+t^4}} \int \limits_0^{\rho_1} \rho^{-2\alpha+1} \,d\rho\\
&=\left\{ \begin{array}{cccc} \log \rho\big|^{\rho_1}_0&=&\infty, & \alpha=1,\hspace{.9cm}\\
\frac{\rho^{-2(\alpha-1)}}{-2(\alpha-1)} \Big |^{\rho_1}_0&=&\infty, & \alpha=2,3,\dots 
\end{array}\right.
\end{align*}
and hence $\norm{\psi}^2_{\L^2}=\infty$.
\end{proof}
We turn now to the remaining case of $\alpha=0$. If  $\psi(s,\rho,\phi)=e^{i\beta\phi_\Gamma} R(\rho)$ is a harmonic spinor on $M^3_\Gamma\cap M^4$, the components of $\chi=\rho \,R(\rho)$ satisfy the differential equations  \eqref{eq:6.8a} and  \eqref{eq:6.8b}, respectively, where
\begin{equation*}
p=\frac 1 \rho -(\log K)_{,\rho}, \qquad q =-f^2=-g^2=-\left ( \frac{\delta K}{\rho} \right)^2,
\end{equation*}
i.\ e., for  $\chi_1$ and $\chi_2$ one obtains  the differential equations
\begin{equation*}
\frac {d^2\chi_i}{d\rho^2}(\rho)+p(\rho)\,\frac {d\chi_i}{d\rho}(\rho)-f^2(\rho)\,\chi_i(\rho) = 0 
\end{equation*}
and these can be integrated explicitly. Indeed, putting
\begin{align*}
\chi_i(\rho)=B_i e^{\delta \, \arsinh \big (\frac {\rho^2(r^2+1)}{t^2}\big)}=\frac {B_i}{t^{2\delta}} \Big (\rho^2(r^2+1)+\sqrt{\rho^4(r^2+1)^2+t^4}\Big )^\delta
\end{align*}
one verifies that
\begin{align*}
\frac {d\chi_i}{d\rho}(\rho)&=\frac{\delta \chi_i(\rho)}{\rho^2(r^2+1)+\sqrt{\rho^4(r^2+1)^2+t^4}}\left ( 2 \rho (r^2+1)+ \frac 1 2 \frac{4\rho^3(r^2+1)^2}{\sqrt{\rho^4(r^2+1)^2+t^4}} \right )\\&=\frac {2\rho(r^2+1)\delta}{\sqrt{\rho^4(r^2+1)^2+t^4}}\chi(\rho)=f(\rho) \,\chi(\rho),\\
\frac {d^2\chi_i}{d\rho^2}(\rho)&= \left (f^2(\rho)+f \left(-\frac 1 \rho+(\log K)_{,\rho} \right ) \right ) \chi_i(\rho)=\big(f^2(\rho)-f(\rho)\, p(\rho)\big )\, \chi_i(\rho). 
\end{align*}
We continue $\psi$ to a spinor on $M^3_\Gamma$ by setting ${\psi}_{|\Gamma}\equiv 0$.
Let now $\dot \phi_\Gamma=1/r$ and $\beta=-1,-2,\dots$, so that $\delta = (r^2+1) \beta/2r \leq \beta<0$. Then one computes 
\begin{gather*}
\int \norm{\psi(s,\rho,\phi)}^2 \, dM^3_\Gamma\\=4 \pi^2 r \frac {B_1^2+B_2^2}{t^{4\delta}}\int \frac{\sqrt 8 \rho(r^2+1)}{\sqrt[4]{\rho^4(r^2+1)^2+t^4}}\Big (\rho^2(r^2+1)+\sqrt{\rho^4(r^2+1)^2+t^4}\Big )^{2\delta} \, d\rho <\infty,
\end{gather*}
so that $\psi \in \L^2(M^3_\Gamma)$. Nevertheless, $\psi$ is not smooth at $\rho=0$, so that $\psi  \notin \D(\overline \Dirac)$.
Thus we have completely determined the $\L^2$--kernel of the Dirac operator in case that $\Gamma$ is a circle in  $\C$ with center at the origin and obtain the following theorem. 

\begin{theorem}
\label{thm:6.3}
Let $\Gamma$ be a generalized circle in $\C$ that arises by a M\"{o}bius  transform from a circle in $\C$ with center at the origin, and $\Dirac$ the Dirac operator on   $(M^3_\Gamma,h_t)$ with respect to the trivial spin structure. Then   
\begin{equation}
\label{eq:6.11}
\Ker_{L^2}\big(\Dirac_{|(M^3_\Gamma \cap M^4)}\big)=\bigoplus_{\beta=-1,-2,\dots} \H_0 \otimes \H_\beta,
\end{equation}
while the $\L^2$--kernel of the Dirac operator and its closure are trivial. In particular, $0 \in \sigma^{\L^2}_{ess}(\overline\Dirac)$.
\end{theorem}
\begin{proof}
Let $\Gamma$ be a circle in $\C$ with center around the origin and radius $r\equiv r_0$.
Without loss of generality we can assume that $\dot \phi_\Gamma=1/r$. For  $\beta=-1,-2,-3,\dots$ and by the previous considerations     
\begin{align*}
\psi_\beta(s,\rho,\phi)&=\frac {e^{i\beta \phi_\Gamma}} \rho e^{\delta \, \arsinh \big(\frac {\rho^2(r^2+1)}{t^2}\big)} \left ( \mnd\begin{array}{c} B_1 \\ B_2 \end{array} \mnd\right )\\
&=\frac {e^{i\beta \phi_\Gamma}} {\rho\, t^{2\delta}}  \Big (\rho^2(r^2+1)+\sqrt{\rho^4(r^2+1)^2+t^4}\Big )^\delta \left ( \mnd\begin{array}{c} B_1 \\ B_2 \end{array} \mnd\right )
\end{align*} 
are harmonic $\L^2$--spinors on $M^3_\Gamma\cap M^4$, where $\delta=(r^2+1)\beta/2 r$, $B_i$ are constants and ${\psi_\beta}_{|\Gamma}\equiv 0$. By Theorem \ref{thm:6.2}, apart from  the trivial representation  no other representations of the $S^1$--action  $\kappa_z$ can occur in the  $\L^2$--kernel of the Dirac  operator and we obtain \eqref{eq:6.11} in case that $\Gamma= \gd B(0,r)$. The general statement then follows from the fact that $M^3_\Gamma$ and $M^3_{A\Gamma}$ are isometric for  $A \in \U(2)$. If, further, $\psi\in \L^2(\Sigma)$ is a harmonic spinor with respect to $\Dirac$, then the regularity theorem for solutions of elliptic differential equations implies that  $\psi\in \Gamma(\Sigma)$. However, since all $\L^2$--harmonic spinors have to be linear combinations of the $\psi_\beta$,  which, nevertheless, are not regular at $\rho=0$, the $\L^2$--kernel of the Dirac operator and its closure turn out to be trivial. Since, by theorem \ref{thm:6.1},  zero  belongs to the spectrum of $\overline \Dirac$, it follows that $0 \in \sigma^{\L^2}_{ess}(\overline\Dirac)$.  
\end{proof}

\section{On the spectrum of the Laplacian}
\label{sec:7}

 In this section we will continue the study of the Laplacian on the hypersurfaces $M^3_\Gamma$, which we began in Section \ref{sec:4}. Unlike the Dirac operator, the spectrum of the Laplacian on an open complete manifold is related to the underlying geometry in a much more intrinsic way. Thus, lower bounds for the Ricci tensor imply upper bounds for its smallest spectral value, and by studying the geodesic flow and the exponential growth of the manifold one obtains statements about the infimum of the essential spectrum of the Laplace operator and vice versa. 

Operating on functions, the Hodge--Laplace operator and the Bochner--Laplace operator coincide, and  we have $\Delta=\nabla^\ast \nabla:\Cinft(M^3_\Gamma) \rightarrow \Cinft(M^3_\Gamma)$ on the hypersurfaces $M^3_\Gamma$; further, since $M^3_\Gamma$ is complete for a closed curve $\Gamma$, $\Delta$ is essentially selfadjoint as an operator in $\L^2(M^3_\Gamma)$ with domain $\Ctest(M^3_\Gamma)$, where the domain of $\overline \Delta$ is given by the Sobolev space $^2\Omega^0(M^3_\Gamma)=\Sob^2(M^3_\Gamma)$. One has $\sigma(\Delta)=\sigma(\overline \Delta)$. Now, for the smallest spectral value of the Laplacian  the following proposition holds in general (see e.g. \cite{eichhorn}). 

\begin{proposition}
Let $(M^n,g)$ be an open complete Riemannian manifold, the components of the Ricci tensor being  bounded from below by $-(n-1) C$, where $C\geq 0$. Then the smallest spectral value of the Laplacian $\mu_0(M^n)$ satisfies
\begin{displaymath}
\mu_0(M^n) \leq \frac{(n-1)^2} 4 C.
\end{displaymath}
\end{proposition}

Hence, as an immediate consequence we obtain the following statement.

\begin{corollary}
\label{cor:7.1}
Let $\Gamma$ be a closed curve in $\C$. Then the smallest spectral value of the Laplacian on the hypersurfaces $(M^3_{A\Gamma},h_t)$ satisfies $\mu_0(M^3_{A\Gamma}) \leq t^{-2}$, where $t\not=0$ and $A\in \U(2)$.
\end{corollary}
\begin{proof} By Theorem \ref{thm:2.1}, $R_{11} \geq R_{22} \geq R_{33}$. Further, since $R_{33}$ is strictly increasing one has that $\inf_\rho R_{33}={R_{33}}_{|\rho=0}=-2/t^2$ so that $R_{ij}\geq -2 C$, where $C=t^{-2}$. The assertion then follows from the proposition above.
\end{proof}

In the sequel we will proceed to find estimates for the infimum of the spectrum of $\overline \Delta$ on the considered hypersurfaces by using again the $\min$--$\max$ principle, and show that it becomes arbitrarily close to zero for any closed curve $\Gamma$, so that $\mu_0(M^3_{A\Gamma})=0$, where $A\in \U(2)$. Since, by Corollary \ref{cor:4.1}, this estimate gives also an estimate for the infimum of  the essential spectrum, we are in position to compute the exponential growth of $M^3_\Gamma$ for an arbitrary closed curve, thus generalizing the results previously obtained in section \ref{sec:8},  since, as already mentioned, the infimum of the essential spectrum of the Laplacian is closely related to the exponential growth of the underlying manifold. More precisely the following theorem proved by Brooks \cite{brooks} holds.
\begin{theorem}[Brooks]
\label{thm:brooks}
Let $(M^n,g)$ be an open complete manifold of infinite volume. Then
\begin{equation*}
\inf \sigma_{\mathrm{ess}}(\overline\Delta)=\frac1{4}\mu_\infty^2.
\end{equation*}
\end{theorem}
Consequently, the exponential growth of the hypersurfaces $M^3_\Gamma$ must be zero for any closed curve $\Gamma$. Let us now prove these assertions. 

First note that for $\phi \in \Sob^2(M^3_\Gamma)$,
\begin{displaymath}
 \int (\phi, \Delta \phi) dM^3_\Gamma=\int (\nabla \phi, \nabla \phi) dM^3_\Gamma
\end{displaymath}
holds, where $(\cdot, \cdot)$ denotes the scalar product in $\Omega^0(M^3_\Gamma)$ and 
 \begin{displaymath}
 (\nabla \phi, \nabla \phi):=\sum(\nabla_{Y_i} \phi, \nabla_{Y_i}\phi) =\sum Y^2_i(\phi)=|\grad \phi|^2.
 \end{displaymath}
By the $\min$--$\max$ principle we have  
\begin{equation}
\label{eq:7.0}
\inf \sigma (\overline \Delta)=\inf_{0 \not=f \in \D(\overline \Delta)}\frac{\int |\grad f|^2 dM^3_\Gamma}{\int |f|^2 dM^3_\Gamma}.
\end{equation}

Now we consider the function
\begin{equation*}
\HM_\epsilon:=\frac{\sqrt{2}}{\sqrt[4]{\rho^4(r^2+1)^2+\epsilon^4}}, \qquad \epsilon >0,
\end{equation*}
which is derived from the trace $\HM$ of the second fundamental form,  and by means of this function we generate estimates for $\sigma_{\mathrm{ess}}(\overline \Delta)$.

\begin{theorem}
\label{thm:7.2}
Let $\Gamma$ be a closed curve in $\C$ and $\overline \Delta$ the closure of the scalar Laplacian on $(M^3_{A\Gamma},h_t)$, where $A \in U(2)$. Then, for arbitrary $\delta>0$,
\begin{equation*}
\inf \sigma_{\mathrm{ess}} (\overline \Delta) <\delta.
\end{equation*}
\end{theorem}
\begin{proof}
By Corollary \ref{cor:3.3}, $\HM^\alpha$ is $\L^2$--integrable over $M^3_\Gamma$ for $\alpha>3/2$. One computes further that 
\begin{align*}
Y_1(\HM_\epsilon^\alpha)&=\frac 1{\sqrt{h_{22}}}  \frac \gd{\gd \rho} \left (\frac 2{\sqrt{\rho^2(r^2+1)^2+\epsilon^4}}\right)^{\alpha/2}=-\frac 1 {\sqrt{h_{22}}}
\frac{\alpha\rho^3(r^2+1)^2}{\rho^4(r^2+1)^2+\epsilon^4}\HM_\epsilon^\alpha,
\end{align*}
the derivatives $Y_2(\HM_\epsilon^\alpha)$ and $Y_3(\HM_\epsilon^\alpha)$ being zero so that
\begin{align*}
|\grad \HM_\epsilon^\alpha|^2=Y^2_1(\HM_\epsilon^\alpha)
=\frac{\alpha^2\rho^4(r^2+1)^2}{2(\rho^4(r^2+1)^2+\epsilon^4)^2} \sqrt{\rho^4(r^2+1)^2+t^4}  \, \HM_\epsilon^{2\alpha}.
\end{align*}
For $\alpha=2$, and assuming $t\leq \epsilon$, the monotony of the integral implies 
\begin{gather*}
\int \HM_\epsilon^4 dM^3_\Gamma=8 \sqrt{8}\pi \int \frac{ \rho^3 (r^2+1)}{(\rho^4(r^2+1)^2+\epsilon^4)\sqrt[4]{\rho^4(r^2+1)^2+t^4}}  \,ds\wedge  d\rho\\ \geq 8 \sqrt{8}\pi \int \frac{ \rho^3 (r^2+1)}{(\rho^4(r^2+1)^2+\epsilon^4)^{5/4}}  \,ds\wedge  d\rho
\\=8 \sqrt{8} \pi \int\limits^L_0 \left [-\frac1{(1+r^2)(\rho^4(r^2+1)^2+\epsilon^4)^{1/4}}\right]_0^\infty \, ds= \frac {8\sqrt{8} \pi}\epsilon \int\limits_0^L \frac{1}{r^2+1} \, ds.  
\end{gather*}
Similarly, under the assumption that $t\leq \epsilon$ one computes
\begin{gather*}
\int |\grad \HM_\epsilon^2|^2 dM^3_\Gamma =16\sqrt{8} \pi \int \frac{\rho^7(r^2+1)^3}{(\rho^4(r^2+1)^2+\epsilon^4)^3} \sqrt[4]{\rho^4(r^2+1)^2+t^4}  \, d\rho \wedge ds\\
\leq 16\sqrt{8} \pi \int \frac{\rho^7(r^2+1)^3}{(\rho^4(r^2+1)^2+\epsilon^4)^{11/4}}   \, d\rho \wedge ds\\
=16 \sqrt{8} \pi \int \limits_0^L \left (  \left [ -\frac{7 \rho^4(r^2+1)^2 +4\epsilon^4}{21 (\rho^4(r^2+1)^2+\epsilon^4)^{7/4}(r^2+1)} \right ] ^\infty_0 \,\right )\,  ds
=\frac{64 \sqrt{8}\pi}{21\epsilon^3} \int \limits_0^L \frac 1{(r^2+1)} \, ds,
\end{gather*}
showing that $\HM_\epsilon^2 \in \D(\overline\Delta)$. Summing up we have 
\begin{displaymath}
\frac{\int |\grad \HM_\epsilon^2|^2 dM^3_\Gamma}{\int \HM_\epsilon^4 dM^3_\Gamma}\leq \frac 8 {21 \epsilon^2} \qquad \text{for all} \quad 0<t<\epsilon;
\end{displaymath}
using \eqref{eq:7.0} one then obtains the stated bound from above for the essential spectrum of the Laplacian since, by Corollary \ref{cor:4.1}, zero can be no $\L^2$--eigenvalue of $\Delta$ and, hence, of $\overline \Delta$.
\end{proof}

\begin{corollary}
Let $t>0$ be arbitrary and $A \in \U(2)$. Then for any closed curve $\Gamma$ in $\C$,  $(M^3_{A\Gamma},h_t)$ has subexponential growth.
\end{corollary}
\begin{proof}
This is a consequence of the theorems \ref{thm:brooks} and \ref{thm:7.2}.
\end{proof}

%%%%%%%%%%%%%%%%%%%%%%%%% Bibliography  %%%%%%%%%%%%%%%%%%%%%%%%%%%
%%%%%%%%%%%%%%%%%%%%%%%%%%%%%%%%%%%%%%%%%%%%%%%%%%%%%%%%%%%%%%%%%%%%

\providecommand{\bysame}{\leavevmode\hbox to3em{\hrulefill}\thinspace}

%%%%%%%%%%%%%%%%%%%%%%%%% End of document  %%%%%%%%%%%%%%%%%%%%%%%%%
%%%%%%%%%%%%%%%%%%%%%%%%%%%%%%%%%%%%%%%%%%%%%%%%%%%%%%%%%%%%%%%%%%%%

\end{document}